\documentclass[journal]{IEEEtran}

\usepackage{amsmath,amssymb,amsfonts,bm,mathtools,amsthm}
\usepackage{algorithm}
\usepackage{algorithmic}
\usepackage{booktabs}
\usepackage{cite}
\usepackage{enumitem}
\usepackage{microtype}
\usepackage{url}
\usepackage{hyperref}
\usepackage{graphicx,multirow,xcolor}

\newcommand{\res}[1]{\textcolor{black}{#1}} 
\newcommand{\rev}[1]{{\color{black}#1}}

\newtheorem{assumption}{Assumption}
\newtheorem{theorem}{Theorem}
\newtheorem{lemma}{Lemma}
\newtheorem{proposition}{Proposition}
\newtheorem{corollary}{Corollary}

\theoremstyle{remark}
\newtheorem{remark}{Remark}

\newcommand{\R}{\mathbb{R}}
\newcommand{\E}{\mathbb{E}}
\newcommand{\Prob}{\mathbb{P}}
\newcommand{\calA}{\mathcal{A}}
\newcommand{\calB}{\mathcal{B}}
\newcommand{\calC}{\mathcal{C}}
\newcommand{\calE}{\mathcal{E}}
\newcommand{\calF}{\mathcal{F}}
\newcommand{\calG}{\mathcal{G}}
\newcommand{\calH}{\mathcal{H}}
\newcommand{\calN}{\mathcal{N}}
\newcommand{\calO}{\mathcal{O}}

\newcommand{\calS}{\mathcal{S}}
\newcommand{\calT}{\mathcal{T}}
\newcommand{\calV}{\mathcal{V}}
\newcommand{\calW}{\mathcal{W}}

\newcommand{\diag}{\operatorname{diag}}
\newcommand{\supp}{\operatorname{supp}}
\newcommand{\tr}{\operatorname{tr}}
\newcommand{\Cov}{\operatorname{Cov}}
\newcommand{\Var}{\operatorname{Var}}
\newcommand{\Proj}{\operatorname{Proj}}
\newcommand{\KL}{\operatorname{KL}}
\newcommand{\TV}{\operatorname{TV}}
\newcommand{\relint}{\operatorname{relint}}
\newcommand{\one}{\bm 1}
\newcommand{\bA}{\bm A}
\newcommand{\bB}{\bm B}
\newcommand{\bG}{\bm G}
\newcommand{\bH}{\bm H}
\newcommand{\bI}{\bm I}
\newcommand{\bL}{\bm L}
\newcommand{\bM}{\bm M}
\newcommand{\bR}{\bm R}
\newcommand{\bV}{\bm V}
\newcommand{\bSigma}{\bm\Sigma}
\newcommand{\bGamma}{\bm\Gamma}
\newcommand{\bLambda}{\bm\Lambda}
\newcommand{\bOmega}{\bm\Omega}
\newcommand{\bXi}{\bm\Xi}
\newcommand{\bDelta}{\bm\Delta}
\newcommand{\bw}{\bm w}
\newcommand{\bx}{\bm x}
\newcommand{\by}{\bm y}
\newcommand{\br}{\bm r}
\newcommand{\bc}{\bm c}
\newcommand{\bg}{\bm g}
\newcommand{\bu}{\bm u}
\newcommand{\bv}{\bm v}
\newcommand{\bz}{\bm z}
\newcommand{\bs}{\bm s}
\newcommand{\bb}{\bm b}
\newcommand{\bh}{\bm h}
\newcommand{\bpsi}{\bm\psi}
\newcommand{\bphi}{\bm\phi}
\newcommand{\beps}{\bm\epsilon}
\newcommand{\bzero}{\bm 0}
\newcommand{\tI}{\widetilde{\bm I}}
\newcommand{\norm}[1]{\left\lVert #1\right\rVert}
\newcommand{\ip}[2]{\left\langle #1,#2\right\rangle}
\newcommand{\dd}{\mathfrak d}

\title{Diffusion-Score-Based Joint Topology Identification and Generative Particle Filtering}

\author{Chuansen Peng, Hanning Tang, Xiaowei Yang and Xiaojing Shen%
\thanks{The work was supported by the Sichuan Provincial Natural Science Foundation Project (Innovative Research Group) 2026NSFSCZY0056 and in part by the National Natural Science
 Foundation of China (NSFC) under Grant U2133208, 62203313. \textit{(Corresponding author: Xiaojing Shen)}\\
 \indent Chuansen Peng and Xiaojing Shen are with School of Mathematics, Sichuan University, Chengdu, Sichuan, 610064, China. (e-mail: pengchuansen@stu.scu.edu.cn; shenxj@scu.edu.cn).\\
 \indent Hanning Tang is with Science and Technology on Information Systems Engineering Laboratory, National University of Defense Technology, Changsha, Hunan, 410073, China. (e-mail: than.ning0@outlook.com).\\
 \indent Xiaowei Yang is with College of Electronic Engineering, National University of Defense Technology, Hefei 230009, China. (e-mail: yxw8290@163.com). }}

\begin{document}
\maketitle

\begin{abstract}
\rev{Inferring the interaction graph of a networked system while tracking its latent state is central to infrastructure, environmental, and brain networks. Joint identification of interaction topology and latent network states is challenging when observations are intermittent and transition noise is unknown and non-Gaussian: a flexible generative transition can absorb graph-induced changes of the conditional mean, whereas a Gaussian transition discards the information carried by the residual shape. We develop a diffusion-score-based framework for joint topology identification and generative particle filtering. The graph governs a leaky-diffusion conditional mean and a frozen diffusion emulator represents the conditionally centered residual. First, a score--topology duality expresses edge-wise transition likelihood gradients as edge projections of the residual diffusion score, so topology learning requires neither transition densities nor differentiation through the sampler. Second, we prove identifiability from intermittent observations, invariance of the population minimizer across noise levels, and global strong convexity above an explicit smoothing threshold for bounded-plus-Gaussian residuals. With known support and oracle state pairs, the asymptotic covariance decreases to the Cram\'er--Rao bound as smoothing vanishes and is strictly smaller than that of any Gaussian quasi-likelihood estimator under an explicit non-Gaussianity condition; for unknown sparse support we characterize the boundary limit, and we quantify the topology bias caused by residual-score mismatch. Third, we couple score-driven graph updates with an approximately fully adapted generative particle filter and, for bounded residuals under explicit excitation and forgetting conditions, prove a mean-square closed-loop bound with a step-size-independent small-gain condition. Synthetic and real-data experiments demonstrate the benefits under non-Gaussian forcing.}
\end{abstract}


\section{Introduction}
\label{sec:intro}

\IEEEPARstart{N}{etworked} dynamical systems arise in traffic and power grids, hydrological networks, sensor systems, and the brain~\cite{shuman2013emerging,ortega2018graph}. Analyzing these systems often requires identifying the underlying interaction graph while tracking latent network states from sparse observations~\cite{giannakis2018topology,mateos2019connecting,dong2019learning}. These tasks are inherently coupled: the graph governs state propagation, whereas topology identification relies on accurate state estimates. Joint inference becomes particularly challenging in high dimensions, where the graph is often unknown, only a time-varying subset of nodes is observed, and the process noise may be heavy-tailed or multimodal. Such uncertainties complicate transition modeling and posterior sampling, while Gaussian approximations may fail to capture important forcing events. Diffusion generative models provide a flexible means of representing complex non-Gaussian dynamics through learned scores. However, exploiting this capability for joint topology identification and state estimation remains challenging.

Three bodies of work each resolve part of this problem. First, \emph{network topology inference} has clarified when and how a graph can be recovered from nodal data. Signal-processing approaches learn Laplacians from smooth or stationary graph signals~\cite{dong2016learning,kalofolias2016graph,segarra2017network,pasdeloup2018characterization,shafipour2021identifying}; dynamic formulations based on graph vector autoregressions, structural equation models, and state-space subspace methods exploit temporal causality~\cite{mei2017signal,isufi2019forecasting,zaman2021online,coutino2020state}; and recent work copes with partial observations by alternating topology estimation with Kalman-type signal reconstruction~\cite{ioannidis2019semiblind,zaman2023online} or by modeling hidden nodes explicitly~\cite{navarro2024joint}. Second, \emph{graph-aware state estimation} has clarified how a topology should be exploited for tracking. Graph Kalman filters~\cite{isufi2020observing}, neural-aided Kalman filters~\cite{revach2022kalmannet,buchnik2024gsp}, and the broader model-based deep learning paradigm~\cite{shlezinger2023model} culminate in GKNet~\cite{sabbaqi2026gknet}, which couples a stochastic heat-diffusion state equation and edge-structured Gaussian uncertainty with a learned Kalman architecture. Third, \emph{generative data assimilation} has shown that diffusion models~\cite{ho2020denoising,song2021sde,karras2022edm} emulate complex non-Gaussian transitions with calibrated uncertainty, from global weather forecasting~\cite{price2025probabilistic} to high-dimensional filtering and smoothing~\cite{rozet2023sda,bao2024ensemble,chen2025flowdas}. Most relevant here, Savary~\emph{et al.}~\cite{savary2026trainingfree} showed that a frozen diffusion emulator, guided by moment-matching posterior sampling (MMPS)~\cite{rozet2024mmps}, implements without retraining \rev{an approximation of} the fully adapted auxiliary particle filter (FA-APF)~\cite{pitt1999filtering,petetin2013optimal}, whose proposal minimizes the conditional variance of the importance weights~\cite{doucet2000sequential,snyder2015performance}. \rev{Our filtering module builds directly on this construction.}

\rev{Despite these advances, a gap remains between structural graph identification and generative state estimation. The obstacle is not merely that the graph is unknown or that the dynamics are non-Gaussian. A sufficiently flexible generative transition can absorb graph-dependent changes of the conditional mean, so that structural identification becomes ill posed even from fully observed states. Conversely, a graph-structured Gaussian transition retains identifiability but discards the likelihood information carried by non-Gaussian residual shapes: graph learners and graph-aware filters are, explicitly or implicitly, second-order procedures, and robust variants require a heavy-tailed family to be fixed a priori~\cite{cardoso2021graphical}. Generative filters, in turn, represent complex transitions but treat the transition mechanism as known. Finally, in joint inference, topology errors propagate into the state estimates, which then bias subsequent graph updates.}

\rev{We address this tension by assigning the graph and the generative model complementary, non-overlapping roles and by converting the learned residual score into graph-parameter likelihood gradients. The graph specifies the conditional transition mean, and a diffusion model represents the conditionally centered residual law. Topology-conditioned particle filtering supplies the latent state pairs that the score-based graph update requires, and the updated graph in turn sharpens the proposal. That a score can deliver parameter gradients is, in itself, a consequence of the chain rule and is related to Fisher's identity for state-space models~\cite{cappe2005inference,poyiadjis2011particle,kantas2015particle} and to Fisher score matching for simulation-based inference~\cite{sui2025fisher}. What is specific to the present problem, and what we analyze, is (i) the exact edge-projected interface induced by the Laplacian parameterization, (ii) the fact that diffusion smoothing changes the curvature and the Fisher information of the topology objective but not its population minimizer, (iii) explicit convexification thresholds and oracle efficiency gains over Gaussian quasi-likelihood, together with the topology bias induced by an imperfect emulator, and (iv) an error analysis of the closed loop with an approximate generative proposal.}

\subsection{Contributions and Organization}

The main contributions are summarized as follows.

\begin{itemize}[leftmargin=1.2em]

\item \rev{\textbf{Identifiable structured residual model and density-free structural score.}
We show that unstructured generative transitions are topology-unidentifiable and adopt a model in which the graph enters only through the conditional mean. A score--topology duality expresses edge-wise transition likelihood gradients as edge projections of the residual diffusion score at every noise level, so topology gradients require only forward evaluations of a frozen emulator; Gaussian graph learning is recovered as generalized least squares. We also show that translation errors of the emulator, such as source-detrending errors, are inherited one-to-one by the topology, and we bound the topology bias caused by a general residual-score mismatch.}

\item \rev{\textbf{Diffusion-noise geometry and oracle efficiency.}
We establish identifiability from intermittent node samples, invariance of the population minimizer across noise levels, and global strong convexity above an explicit smoothing threshold for bounded-plus-Gaussian residuals. For the oracle problem with known support, known residual law, and observed state pairs, the asymptotic covariance decreases monotonically to the Cram\'er--Rao bound as the noise vanishes and strictly improves on every Gaussian quasi-likelihood estimator under an explicit edge-projected non-Gaussianity condition. For unknown sparse supports, we characterize the non-Gaussian boundary limit of the nonnegativity-constrained estimator.}

\item \rev{\textbf{Approximation-aware joint filtering with a conditional closed-loop guarantee.}
TopoScore-FAF couples score-driven topology updates with an approximately fully adapted, topology-conditioned generative particle filter. Proposal errors are measured in total variation, which remains meaningful under support mismatch such as residual clipping. For bounded residuals under explicit excitation, forgetting, and emulator-accuracy conditions, we derive mean-square closed-loop bounds with a step-size-independent small-gain condition and an error floor that separates score, proposal, particle, and stochastic-approximation errors.}

\end{itemize}

Synthetic and real-data experiments evaluate topology recovery, state estimation, and their coupled behavior under non-Gaussian dynamics.

The remainder of the paper is organized as follows. Section~\ref{sec:model} introduces the model and the identifiability trap. Section~\ref{sec:duality} develops the score--topology duality, and Section~\ref{sec:theory} presents the identification, convexification, \rev{efficiency, and robustness} theory. Section~\ref{sec:algo} describes TopoScore-FAF, and Section~\ref{sec:closed} analyzes its closed-loop behavior. Sections~\ref{sec:exp} and~\ref{sec:real} report synthetic and real-data experiments, and the paper closes with a discussion of limitations and conclusions.
 
\emph{Notation.} $\norm{\cdot}$ is the Euclidean norm or the induced operator norm; $\preceq$ is the Loewner order and also the componentwise order for vectors; $\lambda_{\min},\lambda_{\max}$ denote extreme eigenvalues. For a measure $\mu$ and an integrable $f$, $\mu(f)=\int f\,d\mu$. $\TV(\mu,\nu)=\sup_{A}|\mu(A)-\nu(A)|$. $\varphi_{\bSigma}$ is the $\calN(\bzero,\bSigma)$ density. $\calB(\bw,r)$ is a closed Euclidean ball. \rev{For a measurable map $T$, $T_\#\mu$ is the image measure. For $\calS\subseteq\calE_c$, $\bw_{\calS}$ collects the entries indexed by $\calS$, and $\relint$ denotes the relative interior.}

\section{Model and Problem Statement}
\label{sec:model}

\subsection{Graph Parameterization}

Let $\calG^\star=(\calV,\calE^\star,\bm W^\star)$ be an unknown static, undirected, weighted graph on $N$ nodes. Let $\calE_c=\{(i,j):1\le i<j\le N\}$ be the complete candidate edge set, $M=N(N-1)/2$, and fix an orientation with incidence matrix $\bB=[\bb_1,\ldots,\bb_M]\in\R^{N\times M}$, $\bb_e=\bm e_i-\bm e_j$ for $e=(i,j)$. Every candidate Laplacian is
\begin{equation}
\bL(\bw)=\bB\diag(\bw)\bB^\top=\sum_{e=1}^M w_e\bb_e\bb_e^\top,
\qquad \bw\in\R_+^M,
\label{eq:Lparam}
\end{equation}
and the true graph corresponds to $\bw^\star$ with support $\calS^\star$ and $w^\star_{\min}=\min_{e\in\calS^\star}w_e^\star$. \rev{Unless the graph is complete, $\bw^\star$ has zero entries and therefore lies on the boundary of $\R_+^M$.} 

\subsection{Leaky Diffusion Mean and Residual Law}

The latent state $\bx_t\in\R^N$ evolves according to
\begin{equation}
\bx_{t+1}
=
\bA(\bw^\star)\bx_t+\br_{t+1},
\qquad
\bA(\bw)=\beta\bI-\gamma\bL(\bw),
\label{eq:intro-model}
\end{equation}
where $\bw^\star$ denotes the true edge-weight vector,
$\br_{t+1}$ is the transition residual,
$\beta\in(0,1)$ is the known leakage coefficient,
and $\gamma>0$ is the known diffusion step size. The leakage models dissipation (linear drag, losses, outflow) and is essential: \rev{with $\beta=1$, the consensus mode obeys $\one^\top\bx_{t+1}=\one^\top\bx_t+\one^\top\br_{t+1}$, because $\one^\top\bL(\bw)=\bzero^\top$. Whenever the residual has a non-degenerate component along $\one$, this mode is a random walk and admits no stationary law with finite variance.} The admissible set is
\begin{multline}
\calW=\Big\{\bw\in\R_+^M:\ \one^\top\bw\le R_w,\\
\gamma\lambda_{\max}(\bL(\bw))\le 1+\beta-\delta_A\Big\},
\label{eq:W}
\end{multline}
where $R_w>0$ is a prescribed upper bound on the total
edge weight, and $\delta_A\in(0,1+\beta)$ is a stability
margin. Since $\bw\mapsto\lambda_{\max}(\bL(\bw))$ is convex, $\calW$ is convex and compact. For $\bw\in\calW$, the eigenvalues of $\bA(\bw)$ lie in $[-1+\delta_A,\beta]$, hence
\begin{equation}
\norm{\bA(\bw)}\le\bar a:=\max\{\beta,1-\delta_A\}<1 .
\label{eq:abar}
\end{equation}

Conditionally on $\bx_t=\bx$, the residual $\br_{t+1}$ has density $\rho^\star(\cdot\mid\bx)$, independent of the past given $\bx_t$, with
\begin{equation}
\E[\br_{t+1}\mid\bx_t]=\bzero,\qquad \bSigma_r(\bx):=\Cov(\br_{t+1}\mid\bx_t=\bx).
\label{eq:centered}
\end{equation}
The law $\rho^\star$ may be state dependent, heteroscedastic, skewed, or multimodal, but it does not depend on $\bw^\star$. \rev{Condition~\eqref{eq:centered} is a structural normalization of the model, not a property that preprocessing can enforce; its role is explained in Section~\ref{sec:centered}, and the consequences of its violation are quantified in Section~\ref{sec:mismatch}.}

For all $\bu\in\R^M$, the \emph{edge-excitation operator} satisfies
\begin{equation}
\bXi(\bx)=\gamma\bB\diag(\bB^\top\bx),
\quad
\bXi(\bx)\bu=\gamma\bL(\bu)\bx .
\label{eq:Xi}
\end{equation}
For a candidate $\bw$, the induced residual satisfies
\begin{equation}
\bx_{t+1}-\bA(\bw)\bx_t=\br_{t+1}+\bXi(\bx_t)(\bw-\bw^\star),
\label{eq:shift}
\end{equation}
so a topology error acts as a state-dependent translation of the residual. The following quantity governs all curvature constants.

\begin{lemma}
\label{lem:excitation}
Let $\pi$ be a probability measure on $\R^N$ with second-moment matrix $\bSigma_\pi=\int\bx\bx^\top\pi(d\bx)$. For every $\bu\in\R^M$,
\begin{align}
2\gamma^2\lambda_{\min}(\bSigma_\pi)\norm{\bu}^2
&\le
\int\norm{\bXi(\bx)\bu}^2\pi(d\bx)\nonumber\\
&\le
2N\gamma^2\lambda_{\max}(\bSigma_\pi)\norm{\bu}^2 .
\end{align}
Moreover, $\norm{\bXi(\bx)}\le\gamma\sqrt{2N}\norm{\bx}$.
\end{lemma}

\begin{proof}
By~\eqref{eq:Xi}, $\int\norm{\bXi\bu}^2d\pi=\gamma^2\tr(\bL(\bu)\bSigma_\pi\bL^T(\bu))$, which lies between $\gamma^2\lambda_{\min}(\bSigma_\pi)\norm{\bL(\bu)}_F^2$ and $\gamma^2\lambda_{\max}(\bSigma_\pi)\norm{\bL(\bu)}_F^2$. The off-diagonal entries of $\bL(\bu)$ are $-u_{ij}$, each appearing twice, so $\norm{\bL(\bu)}_F^2\ge 2\norm{\bu}^2$. The diagonal entries are $d_i=\sum_{j\ne i}u_{ij}$ with $d_i^2\le(N-1)\sum_{j\ne i}u_{ij}^2$, so $\norm{\bL(\bu)}_F^2\le 2(N-1)\norm{\bu}^2+2\norm{\bu}^2=2N\norm{\bu}^2$. Finally, $\norm{\bXi(\bx)\bu}\le\gamma\norm{\bL(\bu)}_F\norm{\bx}\le\gamma\sqrt{2N}\norm{\bu}\norm{\bx}$.
\end{proof}

\rev{Lemma~\ref{lem:excitation} is an averaged statement. For a single state, $\operatorname{rank}\bXi(\bx)\le\operatorname{rank}\bB=N-1<M$ whenever $N\ge3$, and $\operatorname{range}\bXi(\bx)\subseteq\operatorname{range}\bB=\one^\perp$. Hence no single state pair identifies all edge weights; identification rests on the accumulation of excitation over time, encoded by $\bSigma_\pi\succ0$.}

\subsection{Observation Model}

Observations are
\begin{equation}
\by_t=\bH_t\bx_t+\bv_t,\qquad \bv_t\sim\calN(\bzero,\bR_t),
\label{eq:obs}
\end{equation}
with known $(\bH_t,\bR_t)$, \rev{$\bR_t\succ0$,} and likelihood $g_t(\by_t\mid\bx_t)$. For intermittent node sampling, we use the zero-filled masked form
\begin{equation}
\by_t=\bM_t(\bx_t+\bv_t),\quad \bM_t=\diag(m_{t,1},\ldots,m_{t,N}),
\label{eq:maskobs}
\end{equation}
where $m_{t,i}\in\{0,1\}$ are known mask
indicators and $\bv_t\sim\calN(\bzero,\bR_v)$
denotes the pre-masking observation noise with
time-invariant covariance $\bR_v$. \rev{The zero-filled form is used only to state the moment identities of Theorem~\ref{thm:ident}(a); as a full-dimensional Gaussian model it would have a singular covariance. All likelihood computations use the observed-subspace form: with $\calO_t=\{i:m_{t,i}=1\}$, $d_t=|\calO_t|$, and the row-selection matrix $\bm S_t\in\{0,1\}^{d_t\times N}$, we use~\eqref{eq:obs} with $\by_t\leftarrow\bm S_t\by_t$, $\bH_t=\bm S_t$, and $\bR_t=\bm S_t\bR_v\bm S_t^\top\succ0$.} The targets are the filter $\pi_t^\star=p_{\bw^\star}(\bx_t\mid\by_{1:t})$, the pair posterior $\Gamma_{t+1}^\star=p_{\bw^\star}(\bx_t,\bx_{t+1}\mid\by_{1:t+1})$, and the topology $\bw^\star$.

\subsection{Challenges in Joint Topology and State Inference}

The preceding model and inference targets reveal three technical challenges.

\emph{(D1) Structural identifiability.}
An unrestricted generative residual law can absorb changes in the graph-dependent conditional mean. Different topologies may therefore induce identical transition distributions, even when the states are fully observed. Identifying the graph requires separating structural dynamics from stochastic residuals.

\emph{(D2) Density-free topology learning.}
Diffusion emulators provide residual scores or denoisers rather than tractable transition densities. Standard likelihood-based parameter estimation is therefore not directly applicable. The key challenge is to recover topology gradients from the information supplied by the emulator.

\emph{(D3) Efficient and stable joint inference.}
Non-Gaussian residuals can lead to nonconvex topology objectives, creating a trade-off between optimization and statistical efficiency. Moreover, topology errors affect filtering, while filtering errors bias subsequent graph updates. Joint inference must control this feedback without sacrificing the statistical information contained in non-Gaussian dynamics.

\subsection{Approach and Main Insight}

Our approach assigns distinct statistical roles to the graph and the generative model. The graph governs the conditional transition mean through $\bA(\bw)$, while a frozen diffusion emulator represents the centered residual law $\rho^\star(\cdot\mid\bx)$. This separation removes the shift ambiguity of unrestricted generative transitions. Together with sufficient excitation and observation coverage, it enables topology identification from intermittent node measurements. Since the residual law is independent of the graph parameters, the same emulator can be used for every candidate topology.

The central insight is a duality between topology parameters and residual scores. For consecutive states $\bx,\bx'\in\R^N$ and a candidate topology $\bw$, the transition density is
\begin{equation}
p_{\bw}(\bx'\mid\bx)
=
\rho^\star\bigl(
\bx'-\bA(\bw)\bx\mid\bx
\bigr).
\end{equation}
When $\rho^\star$ is differentiable, define the residual score by
\[
\bs_0(\br\mid\bx)
=
\nabla_{\br}\log\rho^\star(\br\mid\bx),
\qquad
\br=\bx'-\bA(\bw)\bx.
\]
For each candidate edge $(i,j)$, the transition parameter score satisfies
\begin{equation}
\frac{\partial}{\partial w_{ij}}
\log p_{\bw}(\bx'\mid\bx)
=
\gamma(x_i-x_j)
\bigl(s_{0,i}-s_{0,j}\bigr),
\label{eq:intro-duality}
\end{equation}
where $s_{0,i}$ and $s_{0,j}$ are the corresponding components of $\bs_0(\br\mid\bx)$. Thus, the topology score is the edge projection of the residual score. \rev{This identity is a chain rule and holds for any differentiable residual density; centering is not needed for it, but is needed for the graph that maximizes the resulting likelihood to be the structural graph.}

For a noise level $\kappa>0$, let $p_\kappa(\cdot\mid\bx)$ denote the conditional density of $\br+\kappa\beps$, where $\beps\sim\calN(\bzero,\bI)$, and define
\[
\bs_\kappa(\bu\mid\bx)
=
\nabla_{\bu}\log p_\kappa(\bu\mid\bx).
\]
This smoothed residual score is learned by denoising score matching. The score--topology duality extends exactly to the corresponding smoothed transition likelihood. Consequently, topology gradients can be evaluated using a frozen diffusion emulator without computing transition densities or differentiating through the generative sampler. For Gaussian residuals, the resulting objective reduces to generalized least squares.

The diffusion noise level also governs the trade-off between optimization geometry and statistical efficiency. Under a bounded-plus-Gaussian residual model, sufficiently large $\kappa$ makes the population topology objective globally strongly convex, while its unique minimizer remains $\bw^\star$ at every positive noise level. As $\kappa$ decreases, the asymptotic covariance of the state-pair estimator improves and approaches the Cram\'er--Rao bound \rev{in the oracle setting of known support, known residual law, and observed state pairs}. These properties motivate annealed trust-region topology identification without shifting the population minimizer.

Finally, we couple topology identification with generative particle filtering. The current topology determines the transition mean and enters the residual posterior through the observation innovation. Guided residual sampling then produces state estimates and ancestor--descendant particle pairs for subsequent score-based topology updates. This feedback enables joint refinement of topology and latent states within a unified filtering framework.

\subsection{Why the Residual Must Be Centered}
\label{sec:centered}

\begin{proposition}
\label{prop:trap}
Consider the model class in which $\br_{t+1}\mid\bx_t\sim\rho(\cdot\mid\bx_t)$ with $\rho$ unrestricted, i.e., without~\eqref{eq:centered}. For every $\bw\in\calW$, the residual law $\rho_{\bw}(\br\mid\bx):=\rho^\star(\br-\bXi(\bx)(\bw-\bw^\star)\mid\bx)$ produces the same joint law of $\{(\bx_t,\by_t)\}_t$ as $(\bw^\star,\rho^\star)$. Hence $\bw$ is not identifiable, even from fully observed noiseless states.
\end{proposition}

\begin{proof}
If $\br\sim\rho_{\bw}(\cdot\mid\bx)$, then $\br=\br^\star+\bXi(\bx)(\bw-\bw^\star)$ with $\br^\star\sim\rho^\star(\cdot\mid\bx)$, and by~\eqref{eq:shift} $\bA(\bw)\bx+\br=\bA(\bw^\star)\bx+\br^\star$. The transition kernels coincide, and so do all finite-dimensional laws.
\end{proof}

\rev{Proposition~\ref{prop:trap} shows that the data determine only the conditional mean $\bA(\bw)\bx+\E[\br\mid\bx]$ and the centered residual law; how the mean is split between graph and residual must be fixed by a normalization. Condition~\eqref{eq:centered} is that normalization: under it, $\rho_{\bw}$ has conditional mean $\bXi(\bx)(\bw-\bw^\star)$, which is nonzero on a set of positive probability whenever $\bw\neq\bw^\star$ and $\bSigma_x\succ0$ (Lemma~\ref{lem:excitation}), so $\rho_{\bw}$ is excluded.}

\rev{Centering cannot be enforced by fitting. Suppose an emulator is trained on residuals $\widehat\br_{t+1}=\bx_{t+1}-\bA(\widehat\bw)\bx_t$, where $\widehat\bw$ is a least-squares graph estimate computed on an independent data split. By~\eqref{eq:shift}, $\widehat\br_{t+1}=\br_{t+1}+\bXi(\bx_t)(\widehat\bw-\bw^\star)$, so that, conditionally on that split,
\begin{equation}
\E[\widehat\br_{t+1}\mid\bx_t]=\bm m(\bx_t)+\bXi(\bx_t)(\widehat\bw-\bw^\star),
\label{eq:detrendbias}
\end{equation}
with $\bm m(\bx):=\E[\br_{t+1}\mid\bx_t=\bx]$ the possibly nonzero conditional mean of the true residual. Least squares enforces only the normal equations $\sum_t\bXi(\bx_t)^\top\widehat\br_{t+1}=\bzero$, whose population counterpart $\E[\bXi(\bx_t)^\top\bm m(\bx_t)]=\bzero$ is much weaker than $\bm m\equiv\bzero$. Conditional centering is therefore a modeling assumption that must be assessed on held-out data, e.g., through the nonlinear predictability of residuals from states. Both terms in~\eqref{eq:detrendbias} act as translations of the learned residual law, and Proposition~\ref{prop:detrend} shows that such translations are inherited by the topology estimate.}

\section{Residual Emulator and Score--Topology Duality}
\label{sec:duality}

\subsection{Diffusion-Time Residual Densities}
\label{sec:emulator}

We work in the residual scale. For a noise level $\kappa>0$ and $\beps\sim\calN(\bzero,\bI)$ independent of everything else, let
\begin{align}
p_\kappa(\bu\mid\bx)&=\int\varphi_{\kappa^2\bI}(\bu-\br)\rho^\star(\br\mid\bx)\,d\br,
\nonumber\\
\bs_\kappa(\bu\mid\bx)&=\nabla_{\bu}\log p_\kappa(\bu\mid\bx),
\label{eq:pkappa}
\end{align}
be the law of $\br+\kappa\beps$ and its score. A diffusion model with schedule $(\alpha_\tau,\sigma_\tau)$ noises residuals as $\bz=\alpha_\tau(\br+\kappa_\tau\beps)$ with $\kappa_\tau=\sigma_\tau/\alpha_\tau$; its $\bz$-scale score $\widehat\bs_\theta(\tau,\bz,\bx)$ is mapped to the residual scale by $\widehat\bs_\kappa(\bu\mid\bx)=\alpha_\tau\widehat\bs_\theta(\tau,\alpha_\tau\bu,\bx)$. Tweedie's formulas~\cite{efron2011tweedie} read
\begin{align}
\E[\br\mid\bu,\bx]&=\bu+\kappa^2\bs_\kappa(\bu\mid\bx),
\label{eq:tweedie1}\\
\Cov(\br\mid\bu,\bx)&=\kappa^2\bI+\kappa^4\nabla_{\bu}^2\log p_\kappa(\bu\mid\bx).
\label{eq:tweedie2}
\end{align}
\rev{A denoiser $D_\theta(\bu,\kappa,\bx)\approx\E[\br\mid\bu,\bx]$ is pretrained by denoising score matching~\cite{vincent2011connection,song2021sde} on detrended source residuals $\widehat\br=\bx'-\bA(\widehat\bw_{\rm src})\bx$ and then frozen. As $\kappa\to\infty$, the posterior of $\br$ given $\bu$ tends to its prior under mild moment conditions, so $D_\theta(\bu,\kappa,\bx)\to\E[\br\mid\bx]$. Constraining the large-noise output to vanish therefore makes the \emph{emulator's own} residual law centered, but it does not make the training residuals centered (Section~\ref{sec:centered}). In the analysis we treat the emulator as a map on the target state space $\R^N$. Transferring it across networks with different node sets or orderings requires a fixed node alignment or a permutation-equivariant, variable-size architecture; the resulting mismatch is part of the score error analyzed in Section~\ref{sec:mismatch}.}
Two matrices appear repeatedly:
the smoothed Fisher information and the observed information,
\begin{align}
\tI_\kappa(\bx)&=\int\bs_\kappa\bs_\kappa^\top\,p_\kappa(d\bu\mid\bx),
\\
\calH_\kappa(\bu\mid\bx)&=-\nabla_{\bu}^2\log p_\kappa(\bu\mid\bx).
\end{align}
By~\eqref{eq:tweedie2}, $\calH_\kappa\preceq\kappa^{-2}\bI$ for every residual law, and integration by parts gives the Fisher identity $\int\calH_\kappa\,dp_\kappa=\tI_\kappa$. We write $\bm I_0(\bx)$ for the Fisher information of $\rho^\star(\cdot\mid\bx)$ when it is finite.

\subsection{Score--Topology Duality}

For $\bw\in\R^M$ define the noised candidate residual and the loss
\begin{align}
\bu_\kappa(\bw)&=\bx_{t+1}-\bA(\bw)\bx_t+\kappa\beps\nonumber\\
&=\br_{t+1}+\bXi(\bx_t)(\bw-\bw^\star)+\kappa\beps,
\label{eq:ukappa}\\
\ell_\kappa(\bw)&=-\log p_\kappa(\bu_\kappa(\bw)\mid\bx_t).
\label{eq:loss}
\end{align}
Note that $\ell_\kappa$ is the negative log-density of $\bx_{t+1}+\kappa\beps$ under the candidate transition, i.e., a smoothed transition likelihood.

The following result provides the key structural link between diffusion scores and topology learning. 
It shows that the score of the transition with respect to the graph parameters can be obtained exactly by projecting the residual score onto the edge-excitation directions. 
Consequently, a diffusion emulator trained only to represent the residual law supplies the likelihood gradients required for topology estimation without evaluating the transition density or differentiating through the generative sampler.

\begin{theorem}[Score--topology duality]
\label{thm:duality}
For every $\kappa>0$, the gradient of the smoothed transition loss with respect to the edge-weight vector satisfies
\begin{align}
\nabla_{\bw}\ell_\kappa(\bw)
=
-\bXi(\bx_t)^\top
\bs_\kappa\!\left(
\bu_\kappa(\bw)\mid\bx_t
\right).
\label{eq:duality}
\end{align}
Equivalently, for each candidate edge $e=(i,j)$,
\begin{align}
&\frac{\partial \ell_\kappa(\bw)}{\partial w_e}
\nonumber\\
&=
-\gamma(x_{t,i}-x_{t,j})
\left[
s_{\kappa,i}\!\left(\bu_\kappa(\bw)\mid\bx_t\right)
-
s_{\kappa,j}\!\left(\bu_\kappa(\bw)\mid\bx_t\right)
\right].
\label{eq:duality-edge}
\end{align}
Thus, the topology gradient is exactly the residual score projected onto the graph-induced edge-excitation directions.

Moreover, if $\rho^\star(\cdot\mid\bx)$ is differentiable with finite Fisher information, the zero-noise limit yields the exact transition parameter score
\begin{align}
\nabla_{\bw}\log p_{\bw}(\bx'\mid\bx)
&=
\bXi(\bx)^\top
\nabla_{\br}\log\rho^\star(\br\mid\bx),
\nonumber\\
\br&=\bx'-\bA(\bw)\bx.
\label{eq:duality-exact}
\end{align}
Hence, both the smoothed and exact transition scores with respect to the topology are completely determined by the corresponding residual score. [A proof is given in Appendix \ref{app_B}.]
\end{theorem}

Theorem~\ref{thm:duality} resolves (D2) and provides the central interface between the generative and structural components of the model. 
Indeed, Fisher's identity for $\bw$,
\[
\nabla_{\bw}\log p_{\bw}(\by_{1:T})
=
\sum_t
\E\!\left[
\nabla_{\bw}\log p_{\bw}(\bx_{t+1}\mid\bx_t)
\mid\by_{1:T}
\right],
\]
requires only the residual score evaluated along posterior state pairs~\cite{kantas2015particle}, which is precisely the quantity supplied by the diffusion emulator. 
Equation~\eqref{eq:duality-edge} further shows that each topology gradient is a \emph{score-weighted edge-difference correlation}: an edge can be updated only when it is dynamically excited, $x_{t,i}\neq x_{t,j}$, and the residual score exhibits a discrepancy across its endpoints. 
Thus, topology learning requires neither an explicit transition density nor differentiation through the generative network.

\rev{\begin{remark}
\label{rem:scope-duality}
Identity~\eqref{eq:duality} is the chain rule through the affine map~\eqref{eq:ukappa}; it holds for every differentiable residual density, centered or not, and in particular for the law represented by the emulator. With a learned score $\widehat\bs_\kappa$, the vector $-\bXi^\top\widehat\bs_\kappa$ is the gradient of the \emph{surrogate} smoothed likelihood when $\widehat\bs_\kappa(\cdot\mid\bx)$ is a gradient field, and a perturbation of the true topology gradient otherwise; it equals the true likelihood gradient only if $\widehat\bs_\kappa=\bs_\kappa$. The general principle that scores yield parameter gradients is shared with Fisher's identity and with Fisher score matching~\cite{sui2025fisher}; our contribution lies in the edge-structured interface and in the analysis of Sections~\ref{sec:theory}--\ref{sec:closed}. The consequences of $\widehat\bs_\kappa\ne\bs_\kappa$ are quantified in Section~\ref{sec:mismatch}.
\end{remark}}

\begin{remark}
\label{rem:gauss}
If $\rho^\star(\cdot\mid\bx)=\calN(\bzero,\bSigma_r)$, then $\bs_\kappa(\bu)=-(\bSigma_r+\kappa^2\bI)^{-1}\bu$ and, by~\eqref{eq:ukappa}, the population objective is the quadratic $\frac12(\bw-\bw^\star)^\top\E[\bXi^\top(\bSigma_r+\kappa^2\bI)^{-1}\bXi](\bw-\bw^\star)+\text{const}$, i.e., generalized least squares on graph-diffusion residuals. With $\bSigma_r=\bB_0\diag(\bm\alpha^2)\bB_0^\top+\epsilon\bI$, this is the transition term of the GKNet likelihood~\cite{sabbaqi2026gknet} with the graph as unknown. The score-based objective is therefore a strict generalization of second-order graph learning.
\end{remark}

\subsection{Diffusion-Time Topology Objectives}

For a probability measure $\pi$ on states, define the conditional and averaged objectives
\begin{align}
R_\kappa(\bw;\bx)&=\E[\ell_\kappa(\bw)\mid\bx_t=\bx],\nonumber\\
R_{\kappa,\pi}(\bw)&=\int R_\kappa(\bw;\bx)\,\pi(d\bx),
\label{eq:Robj}
\end{align}
where the conditional expectation is over $\br_{t+1}\sim\rho^\star(\cdot\mid\bx)$ and $\beps$. With $\pi_\infty$ the stationary law, $R_\kappa:=R_{\kappa,\pi_\infty}$ is the population objective. The edge-wise conditional gradient is
\begin{align}
\bphi_{\bw}(\bx)&:=\nabla_{\bw}R_\kappa(\bw;\bx)\nonumber\\
&=-\E\big[\bXi(\bx)^\top\bs_\kappa(\bu_\kappa(\bw)\mid\bx)\,\big|\,\bx\big].
\label{eq:phi}
\end{align}
Mixtures over noise levels with weights $\omega(\kappa)$ are handled by linearity; the theory is stated for a single level per stage.

\section{Identification, Convexification, and Efficiency}
\label{sec:theory}

\subsection{Identifiability}

\begin{assumption}
\label{ass:stat}
$(\bx_t)$ is strictly stationary and ergodic with $\bSigma_x=\E[\bx_t\bx_t^\top]\succ0$ and $\E\norm{\br_{t+1}}^2<\infty$.
\end{assumption}

By~\eqref{eq:abar}, $\E\bx_t=\bzero$ in stationarity, so $\bSigma_x$ is the stationary covariance.

\begin{assumption}
\label{ass:mask}
\rev{In~\eqref{eq:maskobs}, the masks $(\bM_t)_t$ are i.i.d.\ across time and independent of the states, the residuals, and the observation noise, with known $\pi^{(0)}_{ij}=\E[m_{t,i}m_{t,j}]>0$ for all $i,j$; hence $\pi^{(1)}_{ij}=\E[m_{t+1,i}m_{t,j}]=\E[m_{t+1,i}]\E[m_{t,j}]>0$. The noise $(\bv_t)_t$ is i.i.d., independent of the states, residuals, and masks, with known covariance $\bR_v$.}
\end{assumption}

\rev{Under Assumptions~\ref{ass:stat}--\ref{ass:mask}, the joint process $(\bx_t,\br_{t+1},\bM_t,\bv_t)_t$ is stationary and ergodic: $(\bx_t,\br_{t+1})$ is a factor of the ergodic process $(\bx_t)$, $(\bM_t,\bv_t)$ is an independent i.i.d.\ (hence mixing) sequence, and the product of an ergodic system with a weakly mixing one is ergodic~\cite{walters1982introduction}. State-dependent (informative) missingness violates Assumption~\ref{ass:mask} and would require an explicit missing-data model.}

\begin{theorem}
\label{thm:ident}
\emph{(a) Moment identifiability.} Under Assumptions~\ref{ass:stat}--\ref{ass:mask}, with $\bGamma_1=\E[\bx_{t+1}\bx_t^\top]$,
\begin{equation}
[\bSigma_x]_{ij}=\frac{\E[y_{t,i}y_{t,j}]}{\pi^{(0)}_{ij}}-[\bR_v]_{ij},
\quad
[\bGamma_1]_{ij}=\frac{\E[y_{t+1,i}y_{t,j}]}{\pi^{(1)}_{ij}},
\label{eq:moments}
\end{equation}
and $\bL(\bw^\star)=\gamma^{-1}(\beta\bI-\bGamma_1\bSigma_x^{-1})$, $w_{ij}^\star=-[\bL(\bw^\star)]_{ij}$. Hence the law of $\{(\bM_t,\by_t)\}$ determines $\bw^\star$, and replacing moments by ergodic sample averages yields a strongly consistent estimator $\bw_{\rm mom}$.

\emph{(b) Score consistency.} Under Assumption~\ref{ass:stat}, for every $\kappa>0$, every $\bw\in\R^M$, and almost every $\bx$,
\begin{multline}
R_\kappa(\bw;\bx)-R_\kappa(\bw^\star;\bx)\\
=\KL\big(p_\kappa(\cdot-\bXi(\bx)(\bw-\bw^\star)\mid\bx)\,\big\|\,p_\kappa(\cdot\mid\bx)\big)\ge0 .
\label{eq:KL}
\end{multline}
Consequently, for every probability measure $\pi$ with $\bSigma_\pi\succ0$, $\bw^\star$ is the unique minimizer of $R_{\kappa,\pi}$ over $\R^M$.

\emph{(c) Fisher curvature.} \rev{Let $\kappa>0$ and let $\pi$ have finite second moments with $\int\tr\bSigma_r\,d\pi<\infty$. Suppose that, for every compact $\calC\subset\R^M$,
\begin{multline}
\E_\pi\Big[\sup_{\bw\in\calC}\big(\norm{\bXi(\bx)}\,\norm{\bs_\kappa(\bu_\kappa(\bw)\mid\bx)}\\
+\norm{\bXi(\bx)}^2\norm{\calH_\kappa(\bu_\kappa(\bw)\mid\bx)}\big)\Big]<\infty,
\label{eq:dominate}
\end{multline}
where $\E_\pi$ also averages over $\br_{t+1}$ and $\beps$; this holds under Assumption~\ref{ass:resid} below. Then} $R_{\kappa,\pi}$ is twice differentiable,
\begin{equation}
\nabla^2R_{\kappa,\pi}(\bw)=\int\E\big[\bXi^\top\calH_\kappa(\bu_\kappa(\bw)\mid\bx)\bXi\,\big|\,\bx\big]\pi(d\bx),
\label{eq:hess}
\end{equation}
and \rev{$\nabla^2R_{\kappa,\pi}(\bw^\star)=\int\bXi^\top\tI_\kappa\bXi\,d\pi\succ0$ whenever $\bSigma_\pi\succ0$. If, in addition, $\bar\lambda_r=\operatorname{ess\,sup}_{\bx}\lambda_{\max}(\bSigma_r(\bx))<\infty$, then}
\begin{equation}
\nabla^2R_{\kappa,\pi}(\bw^\star)=\int\bXi^\top\tI_\kappa\bXi\,d\pi
\succeq\frac{2\gamma^2\lambda_{\min}(\bSigma_\pi)}{\bar\lambda_r+\kappa^2}\bI .
\label{eq:fisherhess}
\end{equation}
[A proof is given in Appendix \ref{app_C}.]
\end{theorem}

\begin{remark}
Part (a) is constructive and provides a consistent warm start under partial sampling without requiring any single time to observe all nodes. \rev{It uses only $\E[\br_{t+1}\bx_t^\top]=\bzero$, which is weaker than~\eqref{eq:centered}.} Part (b) is the property that makes annealing bias-free: \emph{the minimizer does not depend on $\kappa$ nor on the state measure $\pi$}. This contrasts with graduated non-convexity, where minimizers drift along the homotopy. \rev{Part (b) concerns the population objective built from the true residual law; it does not imply that a learned emulator yields an unbiased topology, which is the subject of Section~\ref{sec:mismatch}.} Part (c) shows that local curvature never vanishes when the state is excited, without log-concavity of the residual\rev{; the explicit rate in~\eqref{eq:fisherhess} additionally needs a uniform bound on the conditional residual covariance}.
\end{remark}

\subsection{Diffusion-Time Convexification}

To quantify how noise smoothing removes non-convexity caused by multimodality, we use the following flexible residual class, which contains bounded, multimodal, and Gaussian-mixture residuals.

\begin{assumption}
\label{ass:resid}
Conditionally on $\bx_t=\bx$, $\br_{t+1}=\bc+\bg$ with $\bc\perp\bg$, $\norm{\bc}\le D$ a.s., $\E[\bc\mid\bx]=\bzero$, and $\bg\sim\calN(\bzero,\bSigma_g(\bx))$, $\underline\sigma_g^2\bI\preceq\bSigma_g(\bx)\preceq\bar\sigma_g^2\bI$, with $0\le\underline\sigma_g\le\bar\sigma_g<\infty$.
\end{assumption}

The non-Gaussian component $\bc$ is arbitrary within a ball; it may be supported on a finite set (switching forcing) or on a union of separated regions. Then $\bSigma_r(\bx)\preceq(D^2+\bar\sigma_g^2)\bI$. \rev{Here $D$ is a deterministic almost-sure bound. Laplace and Student-$t$ residuals are not of this form for any finite $D$, so the results that rely on Assumption~\ref{ass:resid} do not cover them. An empirical quantile of residual norms is a heuristic surrogate for $D$, not a certified bound, and the thresholds below should then be read as tuning guidance.}

\begin{theorem}
\label{thm:convex}
Let Assumption~\ref{ass:resid} hold, let $\pi$ be a probability measure on $\R^N$ with $\bSigma_\pi\succ0$, and set $m_3(\pi)=\int\norm{\bx}^3d\pi<\infty$.

\emph{(a) Global strong convexity above the threshold.} If $\kappa>D$, then for all $\bw\in\R^M$,
\begin{equation}
\mu_\kappa(\pi)\bI\preceq\nabla^2R_{\kappa,\pi}(\bw)\preceq L_\kappa(\pi)\bI,
\label{eq:globalconv}
\end{equation}
with
\begin{align}
\mu_\kappa(\pi)&=\frac{2\gamma^2\lambda_{\min}(\bSigma_\pi)\,(1-D^2/\kappa^2)}{\bar\sigma_g^2+\kappa^2},
\nonumber\\
L_\kappa(\pi)&=\frac{2N\gamma^2\lambda_{\max}(\bSigma_\pi)}{\kappa^2}.
\label{eq:muL}
\end{align}
The upper bound holds for every $\kappa>0$.

\emph{(b) Hessian regularity.} For every $\kappa>0$ and $\bw,\bw'$,
\begin{equation}
\norm{\nabla^2R_{\kappa,\pi}(\bw)-\nabla^2R_{\kappa,\pi}(\bw')}\le L_{3,\kappa}(\pi)\norm{\bw-\bw'},
\label{eq:hesslip}
\end{equation}
with $L_{3,\kappa}(\pi)=8D^3\kappa^{-6}\gamma^3(2N)^{3/2}m_3(\pi)$.

\emph{(c) Local curvature and trust radius at every noise level.} For every $\kappa>0$, with $\mu^{\rm loc}_\kappa(\pi)=2\gamma^2\lambda_{\min}(\bSigma_\pi)/(D^2+\bar\sigma_g^2+\kappa^2)$ and $r_\kappa(\pi)=\mu^{\rm loc}_\kappa(\pi)/(2L_{3,\kappa}(\pi))$ \rev{(with $r_\kappa(\pi)=+\infty$ if $D=0$)},
\begin{equation}
\nabla^2R_{\kappa,\pi}(\bw)\succeq\tfrac12\mu^{\rm loc}_\kappa(\pi)\bI
\qquad\forall\,\bw\in\calB(\bw^\star,r_\kappa(\pi)).
\label{eq:localconv}
\end{equation}

\rev{\emph{(d) Gaussian case.} If $D=0$, then $L_{3,\kappa}(\pi)=0$, $R_{\kappa,\pi}$ is quadratic with Hessian $\int\bXi^\top(\bSigma_g+\kappa^2\bI)^{-1}\bXi\,d\pi\succeq\mu_\kappa(\pi)\bI$ for every $\kappa>0$, and~\eqref{eq:localconv} holds on $\R^M$.}
[A proof is given in Appendix \ref{app_D}.]
\end{theorem}

\begin{remark}
For $\kappa>D$, smoothing by $\calN(\bzero,\kappa^2\bI)$ makes every smoothed residual density strongly log-concave, regardless of the number and placement of modes of $\bc$; the topology objective then has a unique stationary point and a condition number at most $N\lambda_{\max}(\bSigma_\pi)(\bar\sigma_g^2+\kappa^2)/[\lambda_{\min}(\bSigma_\pi)(\kappa^2-D^2)]$. For $\kappa\le D$, the objective may have spurious stationary points, but~\eqref{eq:localconv} certifies an explicit basin around $\bw^\star$ whose radius scales as $\kappa^6$. \rev{Convexification is not free: by Theorem~\ref{thm:eff}(b), a larger $\kappa$ also lowers the Fisher information carried by each state pair.}
\end{remark}

\subsection{Statistical Efficiency and the Noise Trade-Off}
\label{sec:eff}

We now study the estimator built on oracle state pairs, which isolates the statistical value of the score. \rev{Because $\bw^\star$ generally lies on the boundary of $\calW$, we distinguish the \emph{oracle-support} problem, in which conventional asymptotic normality and Cram\'er--Rao comparisons apply, from the nonnegativity-constrained problem with unknown support, whose limit is non-Gaussian.}

\rev{For $\calS\subseteq\calE_c$, let $\calW_{\calS}=\{\bw\in\calW:w_e=0\ \forall e\notin\calS\}$ and let $\bXi_{\calS}(\bx)\in\R^{N\times|\calS|}$ collect the columns of $\bXi(\bx)$ indexed by $\calS$, so that $\bXi(\bx)\bw=\bXi_{\calS}(\bx)\bw_{\calS}$ for $\bw\in\calW_{\calS}$. Let}
\begin{equation}
\widehat\bw_{n,\calS}\in\arg\min_{\bw\in\calW_{\calS}}\frac1n\sum_{t=0}^{n-1}\ell_\kappa(\bw;\bx_t,\bx_{t+1},\beps_t),
\label{eq:Mest}
\end{equation}
with i.i.d.\ $\beps_t$ \rev{independent of $(\bx_t)$, and write $\widehat\bw_n:=\widehat\bw_{n,\calE_c}$ for the estimator over $\calW$}. Define, for $\kappa\ge0$,
\begin{align}
\bV_{\kappa,\calS}&=\big(\E[\bXi_{\calS}^\top\tI_\kappa\bXi_{\calS}]\big)^{-1},\qquad \tI_0:=\bm I_0,
\label{eq:Vk}\\
\bV_{{\rm GLS},\calS}(\kappa)&=\big(\E[\bXi_{\calS}^\top(\bSigma_r+\kappa^2\bI)^{-1}\bXi_{\calS}]\big)^{-1},
\label{eq:VGLS}
\end{align}
\rev{where $\bXi_{\calS}$, $\tI_\kappa$, and $\bSigma_r$ are evaluated at the stationary state $\bx_t$. For $\calS=\calE_c$ we drop the subscript and write $\bV_\kappa$ and $\bV_{\rm GLS}(\kappa)$; we also write $\bOmega_{\kappa,\calS}=\bV_{\kappa,\calS}^{-1}$.}

\begin{theorem}
\label{thm:eff}
Let Assumptions~\ref{ass:stat} and~\ref{ass:resid} hold, $\E\norm{\bx_t}^3<\infty$, and \rev{$\bw^\star_{\calS^\star}\in\relint\calW_{\calS^\star}$, i.e., $\one^\top\bw^\star<R_w$ and $\gamma\lambda_{\max}(\bL(\bw^\star))<1+\beta-\delta_A$. All statements concern the oracle problem in which $\calS^\star$, the residual law, and the state pairs are known.}

\emph{(a) Asymptotic normality.} If $\kappa>D$, or if $\kappa>0$ and \rev{$\widehat\bw_{n,\calS^\star}\to\bw^\star$} in probability, then \rev{$\sqrt n(\widehat\bw_{n,\calS^\star}-\bw^\star)_{\calS^\star}\Rightarrow\calN(\bzero,\bV_{\kappa,\calS^\star})$.}

\emph{(b) Monotonicity.} \rev{For every $\calS$,} $\kappa\mapsto\rev{\bV_{\kappa,\calS}}$ is nondecreasing in the Loewner order on $(0,\infty)$.

\emph{(c) Cram\'er--Rao limit.} If $\underline\sigma_g>0$, or more generally if $\bm I_0$ is finite with $\E[\bXi^\top\bm I_0\bXi]$ finite, then \rev{$\bV_{\kappa,\calS^\star}\to\bV_{0,\calS^\star}$} as $\kappa\downarrow0$, and \rev{$\bV_{0,\calS^\star}$} is the inverse Fisher information of the transition model \rev{$\{p_{\bw}:\bw\in\calW_{\calS^\star}\}$, i.e., the H\'ajek--Le Cam bound for regular estimators of $\bw^\star_{\calS^\star}$ based on stationary state pairs}.

\emph{(d) Dominance over Gaussian quasi-likelihood.} \rev{For all $\kappa>0$ and every $\calS$, $\bV_{\kappa,\calS}\preceq\bV_{{\rm GLS},\calS}(\kappa)$, and the same holds at $\kappa=0$ under the conditions of (c). Every Gaussian quasi-maximum-likelihood estimator of $\bw_{\calS^\star}$ over $\calW_{\calS^\star}$ with working covariance $\bm W(\bx)\succ0$ has asymptotic covariance $\bV_{\bm W}\succeq\bV_{{\rm GLS},\calS^\star}(0)$. Assume $\bSigma_r(\bx)\succ0$ a.s.\ and let $\bm F(\bx)=\bm I_0(\bx)-\bSigma_r(\bx)^{-1}\succeq\bzero$. Then $\bV_{0,\calS^\star}=\bV_{{\rm GLS},\calS^\star}(0)$ if and only if $\bXi_{\calS^\star}(\bx)^\top\bm F(\bx)\bXi_{\calS^\star}(\bx)=\bzero$ for $\pi_\infty$-a.e.\ $\bx$. Otherwise $\bV_{0,\calS^\star}\preceq\bV_{{\rm GLS},\calS^\star}(0)$, $\bV_{0,\calS^\star}\ne\bV_{{\rm GLS},\calS^\star}(0)$, and $\tr\bV_{0,\calS^\star}<\tr\bV_{{\rm GLS},\calS^\star}(0)$. If, moreover, $\bm P\bm F(\bx)\bm P\succeq\varphi_0\bm P$ for $\pi_\infty$-a.e.\ $\bx$, where $\bm P$ is the projector onto $\one^\perp$ and $\varphi_0>0$, then
\begin{equation}
\bV_{0,\calS^\star}^{-1}\succeq\bV_{{\rm GLS},\calS^\star}(0)^{-1}+2\gamma^2\varphi_0\lambda_{\min}(\bSigma_x)\bI .
\label{eq:strictgap}
\end{equation}}
[A proof is given in Appendix \ref{app_E}.]
\end{theorem}

\begin{remark}
Theorem~\ref{thm:ident}(a) shows that second moments already identify the graph, so the role of the score is not identification but \emph{efficiency}: Kalman-EM, GLS, GKNet-type Gaussian objectives, and moment estimators are all at best $\bV_{\rm GLS}(0)$-efficient, whereas the score estimator approaches the Cram\'er--Rao bound. The gain $\E[\bXi^\top(\bm I_0-\bSigma_r^{-1})\bXi]$ is a Fisher measure of residual non-Gaussianity seen through the edge excitation. \rev{Because $\operatorname{range}\bXi(\bx)\subseteq\one^\perp$, only non-Gaussianity along the edge-excited directions matters; non-Gaussianity confined to the consensus direction, or to states with $\bXi_{\calS^\star}(\bx)=\bzero$, yields no gain. These are oracle benchmarks: with a learned score, estimated states, and an unknown support they are not attained, and $\bV_{0,\calS^\star}$ is not a semiparametric bound for the problem in which the residual law is learned.}
\end{remark}

\rev{The oracle-support restriction is essential for the Gaussian limit. For an unknown sparse support, the nonnegativity constraints of the inactive edges are binding at $\bw^\star$, and the limit is the projection of a Gaussian onto a polyhedral cone~\cite{shapiro1989asymptotic,geyer1994asymptotics,andrews1999estimation}.}

\rev{\begin{proposition}
\label{prop:boundary}
Let the assumptions of Theorem~\ref{thm:eff} hold, let $\bOmega_\kappa=\E[\bXi^\top\tI_\kappa\bXi]$, $\bm Z\sim\calN(\bzero,\bOmega_\kappa)$, and let $\calT^\star=\{\bh\in\R^M:h_e\ge0\ \forall e\notin\calS^\star\}$ be the tangent cone of $\calW$ at $\bw^\star$. If $\kappa>D$, or if $\kappa>0$ and $\widehat\bw_n\to\bw^\star$ in probability, then
\begin{equation}
\sqrt n(\widehat\bw_n-\bw^\star)\Rightarrow\widehat\bh:=\arg\min_{\bh\in\calT^\star}\Big\{\tfrac12\bh^\top\bOmega_\kappa\bh-\bm Z^\top\bh\Big\},
\label{eq:conelimit}
\end{equation}
i.e., $\widehat\bh$ is the $\bOmega_\kappa$-metric projection of $\bOmega_\kappa^{-1}\bm Z\sim\calN(\bzero,\bV_\kappa)$ onto $\calT^\star$. The limit is Gaussian if and only if $\calS^\star=\calE_c$. If $\calS^\star\ne\calE_c$, the event $\calA^\star=\{\widehat\bh_{\calE_c\setminus\calS^\star}=\bzero\}$ has positive probability, and conditionally on $\calA^\star$, $\widehat\bh_{\calS^\star}=[\bOmega_\kappa]_{\calS^\star\calS^\star}^{-1}\bm Z_{\calS^\star}\sim\calN(\bzero,\bV_{\kappa,\calS^\star})$, the oracle limit of Theorem~\ref{thm:eff}(a).
\end{proposition}}

\rev{\begin{corollary}
\label{cor:Kdraw}
Under the conditions of Theorem~\ref{thm:eff}(a),(c), let the loss in~\eqref{eq:Mest} be replaced by the average of $\ell_\kappa$ over $K_\epsilon$ i.i.d.\ draws $\beps_{t,1},\ldots,\beps_{t,K_\epsilon}$. With $\bOmega=\bOmega_{\kappa,\calS^\star}$, $\bpsi_t=-\bXi_{\calS^\star}(\bx_t)^\top\bs_\kappa(\br_{t+1}+\kappa\beps\mid\bx_t)$, and $\bar\bpsi_t=\E[\bpsi_t\mid\bx_t,\bx_{t+1}]$,
\begin{gather}
\sqrt n(\widehat\bw_{n,\calS^\star}-\bw^\star)_{\calS^\star}\Rightarrow\calN\big(\bzero,\bV^{(K_\epsilon)}_{\kappa}\big),\nonumber\\
\bV^{(K)}_{\kappa}=\bOmega^{-1}\bm S_K\bOmega^{-1},\quad
\bm S_K=K^{-1}\bOmega+(1-K^{-1})\bm S_\infty,\nonumber\\
\bm S_\infty=\E[\bar\bpsi_t\bar\bpsi_t^\top].
\label{eq:Kdraw}
\end{gather}
Moreover, $\bV^{(1)}_\kappa=\bV_{\kappa,\calS^\star}$, $K\mapsto\bV^{(K)}_\kappa$ is nonincreasing, and $\bV_{0,\calS^\star}\preceq\bV^{(K)}_\kappa\preceq\bV_{\kappa,\calS^\star}$ for every $K\ge1$. In particular, $\tr\bV_{\kappa,\calS^\star}$ is a conservative prediction when $K_\epsilon>1$ draws are used.
\end{corollary}}

\begin{corollary}
\label{cor:anneal}
Let $\kappa_1>D$ and $\kappa_2\in(0,\kappa_1)$, and let $\pi$ satisfy the conditions of Theorem~\ref{thm:convex}. Stage~I minimizes $R_{\kappa_1,\pi}$ over $\calW$. Stage~II minimizes $R_{\kappa_2,\pi}$ over
\begin{equation}
\calW_2=\calW\cap\calB\big(\widehat\bw_{\rm I},\,r_{\kappa_2}(\pi)/2\big).
\label{eq:W2}
\end{equation}
On the event $\norm{\widehat\bw_{\rm I}-\bw^\star}\le r_{\kappa_2}(\pi)/2$, one has $\bw^\star\in\calW_2\subseteq\calB(\bw^\star,r_{\kappa_2}(\pi))$, $R_{\kappa_2,\pi}$ is $\frac12\mu^{\rm loc}_{\kappa_2}(\pi)$-strongly convex and $L_{\kappa_2}(\pi)$-smooth on $\calW_2$, and $\bw^\star$ is its unique minimizer on $\calW_2$. If Stage~I has mean-square error $e_{\rm I}^2$, this event has probability at least $1-4e_{\rm I}^2/r_{\kappa_2}(\pi)^2$.
\end{corollary}

\begin{proof}
The inclusions follow from the triangle inequality, convexity of $\calW_2$ from that of balls and $\calW$, curvature from~\eqref{eq:localconv} and~\eqref{eq:globalconv}, and uniqueness from Theorem~\ref{thm:ident}(b). The probability bound is Markov's inequality.
\end{proof}

\rev{\begin{remark}
\label{rem:switch}
Corollary~\ref{cor:anneal} is a population statement that is conditional on Stage~I having entered the Stage~II basin; a fixed switching time does not certify this event. A checkable sufficient condition follows from strong convexity: since $\nabla R_{\kappa_1,\pi}(\bw^\star)=\bzero$ and $R_{\kappa_1,\pi}$ is $\mu_{\kappa_1}(\pi)$-strongly convex on $\R^M$, $\mu_{\kappa_1}(\pi)\norm{\bw-\bw^\star}^2\le\ip{\nabla R_{\kappa_1,\pi}(\bw)}{\bw-\bw^\star}$, hence
\begin{equation}
\norm{\bw-\bw^\star}\le\norm{\nabla R_{\kappa_1,\pi}(\bw)}/\mu_{\kappa_1}(\pi)\qquad\forall\bw\in\R^M.
\label{eq:certificate}
\end{equation}
Switching once $\norm{\nabla R_{\kappa_1,\pi}(\bw)}\le\mu_{\kappa_1}(\pi)r_{\kappa_2}(\pi)/2$ therefore guarantees the event of Corollary~\ref{cor:anneal}. In practice, the gradient is estimated from posterior pairs and a learned score, so~\eqref{eq:certificate} holds only up to the errors $\varepsilon_s$ and $\varepsilon_F$ of Section~\ref{sec:closed}, and the constants must be estimated. The experiments use fixed switching times, which is a heuristic schedule.
\end{remark}}

Iterating Corollary~\ref{cor:anneal} over a decreasing sequence $\kappa_1>\kappa_2>\cdots$ trades the global guarantee of large noise for the efficiency of small noise. Because $\bw^\star$ minimizes every stage objective \rev{of the oracle model}, no stage introduces \rev{population} bias.

\rev{\subsection{Robustness to Residual-Score Mismatch}
\label{sec:mismatch}

The preceding results assume the true residual law. A frozen emulator departs from it through three distinct mechanisms: \emph{(M1) candidate-query shift}, since scores are evaluated at residuals translated by $\bXi(\bx)(\bw-\bw^\star)$ relative to the training distribution; \emph{(M2) detrending error}, since training residuals inherit the translation~\eqref{eq:detrendbias}; and \emph{(M3) source--target shift}, since the conditional residual law of the target, including its node set and ordering, may differ from that of the source. Let $\widehat\rho(\cdot\mid\bx)$ be the residual law represented by the emulator, $\widehat p_\kappa$ its Gaussian smoothing, $\widehat\bs_\kappa$ its score, and $\widehat R_{\kappa,\pi}$ the objective~\eqref{eq:Robj} with $\widehat p_\kappa$ in place of $p_\kappa$, the expectation still being taken under the true target data. Mechanism (M2) has an exact characterization.

\begin{proposition}
\label{prop:detrend}
Suppose $\widehat\rho(\cdot\mid\bx)=\rho^\star(\cdot-\bb(\bx)\mid\bx)$ for a measurable $\bb$ with $\E_\pi\norm{\bb(\bx)}^2<\infty$, let Assumption~\ref{ass:stat} hold, let $\kappa>0$, and let $\bSigma_\pi\succ0$.

(i) If $\bb(\bx)=\bXi(\bx)\bm\delta$ for $\pi$-a.e.\ $\bx$ and some $\bm\delta\in\R^M$, then $\bw^\star+\bm\delta$ is the unique minimizer of $\widehat R_{\kappa,\pi}$ over $\R^M$, for every $\kappa>0$.

(ii) If $\rho^\star(\cdot\mid\bx)=\calN(\bzero,\bSigma_r(\bx))$ and $\bLambda_\kappa=(\bSigma_r+\kappa^2\bI)^{-1}$, the unique minimizer is
\begin{equation}
\bar\bw=\bw^\star+\big(\E_\pi[\bXi^\top\bLambda_\kappa\bXi]\big)^{-1}\E_\pi[\bXi^\top\bLambda_\kappa\bb].
\label{eq:gaussbias}
\end{equation}
[A proof is given in Appendix \ref{app_1}.]
\end{proposition}

Proposition~\ref{prop:detrend}(i) shows that a source-detrending error $\bm\delta=\widehat\bw_{\rm src}-\bw^\star$ on the same network is transferred one-to-one to the topology at every noise level, so annealing cannot remove it. A residual mean $\bm m(\bx)$ outside the edge directions produces a bias equal to its (Kullback--Leibler) projection onto those directions, as~\eqref{eq:gaussbias} shows explicitly. The general case is controlled by the score error at the queried residuals.

\begin{proposition}
\label{prop:transfer-bias}
Let $\calW_s\subseteq\R^M$ be closed and convex with $\bw^\star\in\calW_s$, let $R_{\kappa,\pi}$ satisfy~\eqref{eq:dominate} and be $\mu$-strongly convex on $\calW_s$, and let $C_\Xi^2=\E_\pi\norm{\bXi(\bx)}^2<\infty$. Define $\widehat\bg_\kappa(\bw)=-\E\big[\bXi(\bx)^\top\widehat\bs_\kappa(\bu_\kappa(\bw)\mid\bx)\big]$, with $\bx\sim\pi$ and $\bu_\kappa(\bw)$ generated by the true target transition, and suppose
\begin{equation}
\sup_{\bw\in\calW_s}\E\Big[\norm{\widehat\bs_\kappa(\bu_\kappa(\bw)\mid\bx)-\bs_\kappa(\bu_\kappa(\bw)\mid\bx)}^2\Big]\le\varepsilon_{\rm tar}^2 .
\label{eq:epstar}
\end{equation}
Then $\sup_{\bw\in\calW_s}\norm{\widehat\bg_\kappa(\bw)-\nabla R_{\kappa,\pi}(\bw)}\le C_\Xi\varepsilon_{\rm tar}$. Moreover, every $\widetilde\bw\in\calW_s$ satisfying the variational inequality $\ip{\widehat\bg_\kappa(\widetilde\bw)}{\bw-\widetilde\bw}\ge0$ for all $\bw\in\calW_s$ (in particular, every interior zero of $\widehat\bg_\kappa$ and every fixed point of the projected update) obeys
\begin{equation}
\norm{\widetilde\bw-\bw^\star}\le C_\Xi\varepsilon_{\rm tar}/\mu .
\label{eq:transferbias}
\end{equation}
[A proof is given in Appendix \ref{app_1}.]
\end{proposition}

Three remarks are in order. First, under Assumption~\ref{ass:resid}, $\bs_\kappa(\cdot\mid\bx)$ is globally Lipschitz with constant $\ell_\kappa=\max\{\kappa^{-2},D^2\kappa^{-4}\}$ (Appendix~\ref{app_D}), so a translation mismatch $\bb$ gives $\varepsilon_{\rm tar}\le\ell_\kappa(\E_\pi\norm{\bb}^2)^{1/2}$ uniformly in $\bw$, consistent in order with Proposition~\ref{prop:detrend}. Second, a small distance between source and target residual laws (e.g., in Wasserstein distance) does not imply a small $\varepsilon_{\rm tar}$, because scores are not continuous in such metrics; the score error must be assessed directly, for instance through held-out denoising score-matching risk. Third, $\varepsilon_{\rm tar}$ is the population counterpart of the emulator error $\varepsilon_s$ that enters the closed-loop analysis (Assumption~\ref{ass:score}).}

\section{TopoScore-FAF: Joint Topology Learning and Generative Filtering}
\label{sec:algo}

\subsection{Topology-Conditioned Fully Adapted Filtering in Residual Coordinates}

For a fixed $\bw$, the transition kernel $P_{\bw}(\bx,d\bx')$ is the law of $\bA(\bw)\bx+\br$, $\br\sim\rho^\star(\cdot\mid\bx)$. The frozen emulator therefore emulates the transition for \emph{every} candidate topology\rev{, to the extent that it represents $\rho^\star$ on the target (Section~\ref{sec:mismatch})}. Given an ancestor $\bx$ and $\by_{t+1}$, define the residual innovation
\begin{equation}
\tilde\by(\bx;\bw)=\by_{t+1}-\bH_{t+1}\bA(\bw)\bx=\bH_{t+1}\br+\bv_{t+1},
\label{eq:innov}
\end{equation}
\rev{where $(\by_{t+1},\bH_{t+1},\bR_{t+1})$ are taken in the observed-subspace form of Section~\ref{sec:model}.}
The fully adapted proposal $q^\star_{\bw}(\bx'\mid\bx)\propto g_{t+1}(\by_{t+1}\mid\bx')P_{\bw}(\bx,d\bx')$ is the image under $\br\mapsto\bA(\bw)\bx+\br$ of the residual posterior $p(\br\mid\bx,\tilde\by)\propto\calN(\tilde\by;\bH\br,\bR)\rho^\star(\br\mid\bx)$: \emph{the topology enters only through the innovation}. This is a linear inverse problem on the residual and is solved \rev{approximately} by MMPS guidance~\cite{rozet2024mmps,savary2026trainingfree} (Algorithm~\ref{alg:guided}). For brevity, let
$D_\theta:=D_\theta(\bu,\kappa,\bx)$ and
$\bm J_D:=\nabla_{\bu}D_\theta(\bu,\kappa,\bx)$,
where the Jacobian is evaluated at fixed
$(\kappa,\bx)$. With $\bV=\kappa^2\bm J_D$,
as motivated by~\eqref{eq:tweedie2},
the guided score is
\begin{equation}
\widehat\bs^{\,y}_\kappa(\bu)=\frac{D_\theta-\bu}{\kappa^2}+\bm J_D^\top\bH^\top\bm\Sigma_y^{-1}\big(\tilde\by-\bH D_\theta\big),
\label{eq:mmps}
\end{equation}
with $\bm\Sigma_y=\bR+\bH\bV\bH^\top$,
where the linear system is solved by conjugate gradients with Jacobian--vector products. \rev{Because MMPS uses a Gaussian approximation of $p(\br\mid\bu,\bx)$, a finite number of sampler steps, and no importance correction, the resulting filter is \emph{approximately} fully adapted; Proposition~\ref{prop:local} quantifies the approximation.}

\emph{Ancestor weights.} The FA-APF weights ancestors by the predictive likelihood $C_{\bw}(\bx)=\int\calN(\tilde\by;\bH\br,\bR)\rho^\star(\br\mid\bx)d\br$. We use the moment-matched estimate
\begin{equation}
\widehat C_{\bw}(\bx)=\calN\big(\tilde\by(\bx;\bw);\ \bzero,\ \bR+\bH\widehat\bSigma_r(\bx)\bH^\top\big),
\label{eq:ancestor}
\end{equation}
where $\widehat\bSigma_r(\bx)=\kappa_{\max}^2\nabla_{\bu}D_\theta(\bzero,\kappa_{\max},\bx)$ is the large-noise Tweedie covariance, which tends to $\bSigma_r(\bx)$ as $\kappa_{\max}\to\infty$. Setting $\widehat\bSigma_r=\bzero$ recovers the Dirac approximation of~\cite{savary2026trainingfree}. Monte Carlo estimates with emulator samples are an unbiased alternative. In all cases, ancestor selection is driven by the graph mean through $\tilde\by$.

\subsection{Score-Driven Topology Update}

After the filtering step, the ancestor--descendant particles $(\bx_t^{(a_j)},\bx_{t+1}^{(j)})$ form the pair measure $\widehat\Gamma_{t+1}$. The topology gradient is estimated with $J_g$ pairs and fresh noise:
\begin{equation}
\widehat\bg_t=-\frac1{J_g}\sum_{k=1}^{J_g}\bXi(\bx^{(k)})^\top\widehat\bs_\kappa\big(\bu^{(k)}\,\big|\,\bx^{(k)}\big),
\label{eq:ghat}
\end{equation}
with $\bu^{(k)}=\bx'^{(k)}-\bA(\bw_t)\bx^{(k)}+\kappa\beps_k$; that is, $[\widehat\bg_t]_{ij}=-\frac{\gamma}{J_g}\sum_k(x^{(k)}_i-x^{(k)}_j)(\widehat s_i-\widehat s_j)$, followed by
\begin{equation}
\bw_{t+1}=\Proj_{\calW_{\rm cur}}\big[\bw_t-\eta\,\widehat\bg_t\big],
\label{eq:update}
\end{equation}
where $\calW_{\rm cur}$ is $\calW$ in Stage~I and the trust region~\eqref{eq:W2} in later stages. The update requires only forward denoiser calls; no backpropagation through the network is needed.

\begin{algorithm}[t]
\caption{TopoScore-FAF}
\label{alg:main}
\begin{algorithmic}[1]
\STATE \textbf{Input:} $\by_{1:T}$; $\bH_t,\bR_t$; $\beta,\gamma,\bB$; frozen $D_\theta$; $J,J_g$; $\eta$; schedule $\kappa_1>\kappa_2>\cdots$ with switching times; $\calW$\rev{; clipping radius $D_{\rm clip}\in(0,\infty]$}.
\STATE $\bw_0\leftarrow\Proj_{\calW}(\bw_{\rm mom})$ (Theorem~\ref{thm:ident}(a)); $\calW_{\rm cur}\leftarrow\calW$; $\kappa\leftarrow\kappa_1$; particles $\{\bx_0^{(j)}\}$, weights $1/J$.
\FOR{$t=0,\ldots,T-1$}
  \STATE $\tilde\by^{(j)}\leftarrow\by_{t+1}-\bH_{t+1}\bA(\bw_t)\bx_t^{(j)}$; compute $\widehat C^{(j)}$ by~\eqref{eq:ancestor}.
  \STATE Draw ancestors $a_j\propto W_t^{(j)}\widehat C^{(j)}$.
  \STATE $\br^{(j)}\leftarrow$ Algorithm~\ref{alg:guided}$(\bx_t^{(a_j)},\tilde\by^{(a_j)})$, \rev{projected onto $\calB(\bzero,D_{\rm clip})$ ($D_{\rm clip}=\infty$ disables clipping; see Proposition~\ref{prop:local}(iii))}.
  \STATE $\bx_{t+1}^{(j)}\leftarrow\bA(\bw_t)\bx_t^{(a_j)}+\br^{(j)}$; $W_{t+1}^{(j)}\leftarrow1/J$.
  \STATE $\widehat\Gamma_{t+1}\leftarrow J^{-1}\sum_j\delta_{(\bx_t^{(a_j)},\bx_{t+1}^{(j)})}$; compute $\widehat\bg_t$ by~\eqref{eq:ghat}.
  \STATE $\bw_{t+1}\leftarrow\Proj_{\calW_{\rm cur}}[\bw_t-\eta\widehat\bg_t]$.
  \IF{$t+1$ is a switching time to $\kappa_{k+1}$ \rev{(fixed schedule, or the certificate of Remark~\ref{rem:switch})}}
     \STATE $\calW_{\rm cur}\leftarrow\calW\cap\calB(\bw_{t+1},r_{\kappa_{k+1}}/2)$; $\kappa\leftarrow\kappa_{k+1}$.
  \ENDIF
\ENDFOR
\STATE \textbf{Output:} $\widehat\pi_t=J^{-1}\sum_j\delta_{\bx_t^{(j)}}$, $\bL(\bw_T)$, $\widehat\calS=\{e:w_{T,e}>\vartheta\}$.
\end{algorithmic}
\end{algorithm}

\begin{algorithm}[t]
\caption{Guided residual sampler (residual scale)}
\label{alg:guided}
\begin{algorithmic}[1]
\STATE \textbf{Input:} ancestor $\bx$, innovation $\tilde\by$, levels $\kappa_{\max}=\varkappa_K>\cdots>\varkappa_0=0$.
\STATE $\bu\sim\calN(\bzero,\varkappa_K^2\bI)$.
\FOR{$k=K,\ldots,1$}
 \STATE $D\leftarrow D_\theta(\bu,\varkappa_k,\bx)$; solve $(\bR+\varkappa_k^2\bH\bm J_D\bH^\top)\bm\lambda=\tilde\by-\bH D$ by CG using JVPs.
 \STATE $\bs\leftarrow(D-\bu)/\varkappa_k^2+\bm J_D^\top\bH^\top\bm\lambda$ \hfill (Eq.~\eqref{eq:mmps})
 \STATE $\bu\leftarrow\bu+(\varkappa_k-\varkappa_{k-1})\varkappa_k\bs$ \hfill (probability-flow Euler/DDIM step)
\ENDFOR
\STATE \textbf{Return} $\br=\bu$.
\end{algorithmic}
\end{algorithm}

\subsection{Proposal Quality and One-Step Error}

Fix $\bw$, an ancestor $\bx$, and $\by_{t+1}$; let $h_{\bw}(\bx')=g_{t+1}(\by_{t+1}\mid\bx')p_{\bw}(\bx'\mid\bx)$ and $C_{\bw}=\int h_{\bw}$. \rev{Let $\widehat q_j$ denote the law of the output of Line~6--7 of Algorithm~\ref{alg:main} given ancestor $j$, including clipping, and let $T_D(\br)=\br\min\{1,D/\norm{\br}\}$ be the radial projection onto $\calB(\bzero,D)$.}

\begin{proposition}
\label{prop:local}
\emph{(i)} For any proposal $q$ with $\supp q\supseteq\supp q^\star_{\bw}$ and $W_q=h_{\bw}/q$, one has $\Var_q[W_q]/C_{\bw}^2=\chi^2(q^\star_{\bw}\|q)$ and $\E_q[W_q]^2/\E_q[W_q^2]=(1+\chi^2(q^\star_{\bw}\|q))^{-1}$. In particular $q^\star_{\bw}$ is the unique zero-variance proposal, and $\TV(q,q^\star_{\bw})\le\frac12\chi^2(q^\star_{\bw}\|q)^{1/2}$. \rev{If $\supp q\not\supseteq\supp q^\star_{\bw}$, then $\chi^2(q^\star_{\bw}\|q)=\infty$; this occurs for hard-clipped proposals whenever the true residual has unbounded support.}

\emph{(ii)} Let $\Phi^{\bw}_{t+1}(\widehat\Gamma_t)=\sum_jp_j\delta_{\bx^{(j)}}\otimes q^\star_{\bw,j}$ be the exact one-step pair update of the current particle measure, with $p_j\propto W^{(j)}C_{\bw}(\bx^{(j)})$. Suppose $|\widehat C^{(j)}/C_{\bw}(\bx^{(j)})-1|\le\varepsilon_C<1$ and \rev{$\TV(\widehat q_j,q^\star_{\bw,j})\le\delta_q$} for all $j$. Then for every test function $\bm f$ with $\sup\norm{\bm f}\le1$ that is measurable with respect to the past,
\begin{multline}
\big(\E\norm{\widehat\Gamma_{t+1}(\bm f)-\Phi^{\bw}_{t+1}(\widehat\Gamma_t)(\bm f)}^2\big)^{1/2}\\
\le\varepsilon_F:=\frac1{\sqrt J}+\frac{2\varepsilon_C}{1-\varepsilon_C}+\rev{2\delta_q} .
\label{eq:epsF}
\end{multline}
\rev{When $\chi^2(q^\star_{\bw,j}\|\widehat q_j)\le\bar\chi^2$ for all $j$, part (i) gives $\delta_q\le\bar\chi/2$.}

\rev{\emph{(iii) Clipping.} For every law $\nu$ on $\R^N$, $\TV((T_D)_\#\nu,\nu)\le\nu(\norm{\br}>D)$. If $\supp q^\star_{\bw,j}$ lies in the image of $\calB(\bzero,D)$ under $\br\mapsto\bA(\bw)\bx^{(j)}+\br$, then clipping at radius $D$ cannot increase $\TV(\widehat q_j,q^\star_{\bw,j})$. Finally, for the clipped transition $P_{\bw,D}(\bx,\cdot)$, the law of $\bA(\bw)\bx+T_D(\br)$, one has $\TV(P_{\bw,D}(\bx,\cdot),P_{\bw}(\bx,\cdot))\le\tau_D(\bx):=\Prob(\norm{\br_{t+1}}>D\mid\bx_t=\bx)$ and $\norm{\E[T_D(\br_{t+1})\mid\bx_t=\bx]}\le\E[(\norm{\br_{t+1}}-D)_+\mid\bx_t=\bx]$.} [A proof is given in Appendix \ref{app_2}.]
\end{proposition}

The topology error does not enter~\eqref{eq:epsF}: it is accounted for separately through the sensitivity of the exact filter to $\bw$ (Lemma~\ref{lem:sens}). This separation avoids treating the approximate sampler as exactly fully adapted. \rev{Formulating (ii) in total variation keeps the bound meaningful when the guided proposal does not dominate the exact one, e.g., under clipping. Part (iii) shows that clipping with $D_{\rm clip}\ge D$ is harmless under the bounded residuals of Assumption~\ref{ass:compact}. Outside that regime, clipping adds at most the generated tail mass to $\delta_q$, introduces the model bias $\tau_D$, and perturbs the residual mean by at most $\E[(\norm{\br}-D)_+\mid\bx]$; none of these effects is covered by Theorem~\ref{thm:closed}. Proposition~\ref{prop:local} does not show that $\delta_q$ is small, which must be assessed empirically or through a separate analysis of the guided sampler.}

\emph{Complexity.} Per time step, filtering costs $J\cdot K\cdot(1+n_{\rm CG})$ denoiser passes with Jacobian--vector products, as in~\cite{savary2026trainingfree}. The topology update costs $J_g$ forward passes plus $O(J_gM)$ edge operations via $\bXi^\top\bs=\gamma\diag(\bB^\top\bx)\bB^\top\bs$. Candidate screening by physical locality reduces $M$ from $O(N^2)$ to $O(N)$ without changing the theory. A polyhedral inner approximation of $\calW$ is $\{\bw\ge\bzero,\one^\top\bw\le R_w,\,2\gamma\max_i\sum_jw_{ij}\le1+\beta-\delta_A\}$ (Gershgorin), whose projection is a quadratic program.

\section{Closed-Loop Analysis}
\label{sec:closed}

We analyze Algorithm~\ref{alg:main} within one stage at a fixed $\kappa$, with projection set $\calW_{\rm cur}=:\calW_s\ni\bw^\star$. \rev{The analysis is conditional: it covers a bounded-residual regime under explicit excitation, forgetting, and emulator-accuracy assumptions. Unlike Theorem~\ref{thm:eff}, it does not require $\bw^\star$ to be an interior point, because the projected update only uses $\nabla R_t(\bw^\star)=\bzero$, which holds since $\bw^\star$ minimizes $R_t$ over $\R^M$.}

\subsection{Assumptions}

\begin{assumption}
\label{ass:compact}
Assumption~\ref{ass:resid} holds with $\bar\sigma_g=0$, so $\norm{\br_{t+1}}\le D$, and $\rho^\star(\cdot\mid\bx)$ has finite Fisher information with $\bar J:=\sup_{\bx}\tr\bm I_0(\bx)<\infty$. With $R_x:=\max\{\norm{\bx_0},D/(1-\bar a)\}$, true and particle states satisfy $\norm{\bx}\le R_x$, the latter enforced by $D_{\rm clip}=D$ in Algorithm~\ref{alg:main}.
\end{assumption}

By~\eqref{eq:abar}, $\norm{\bA(\bw)\bx+\br}\le\bar aR_x+D\le R_x$, so the ball is invariant for all $\bw\in\calW$. \rev{By Proposition~\ref{prop:local}(iii), clipping at $D$ does not increase the proposal error in this regime. Assumption~\ref{ass:compact} excludes Gaussian, Laplace, and Student-$t$ residuals as well as the Gaussian-perturbed switching law; clipping the generated samples does not make the true residual bounded. An extension to unbounded residuals through weighted norms~\cite{douc2009forgetting} is plausible but is not established here.}

\begin{assumption}
\label{ass:PE}
There is $\lambda_{\rm PE}>0$ such that $\lambda_{\min}(\E[\bx_t\bx_t^\top\mid\by_{1:t}])\ge\lambda_{\rm PE}$ a.s.\ for all $t$.
\end{assumption}

\begin{proposition}
\label{prop:PE}
Under Assumption~\ref{ass:resid} with $\underline\sigma_g>0$ and~\eqref{eq:obs} with $\bR_t\succeq\underline r\bI$, $\norm{\bH_t}\le h$,
\begin{align}
\E[\bx_t\bx_t^\top\mid\by_{1:t}]&\succeq\Cov(\bx_t\mid\by_{1:t})\nonumber\\
&\succeq\big(\underline\sigma_g^{-2}+h^2/\underline r\big)^{-1}\bI\quad\text{a.s.}
\end{align}
\end{proposition}

\begin{proof}
The first inequality holds because $\E[\bx\bx^\top\mid\cdot]=\Cov(\bx\mid\cdot)+\bm m\bm m^\top$. Let $\calG_t=\sigma(\by_{1:t},\bx_{t-1},\bc_t)$, which contains $\sigma(\by_{1:t})$. By the law of total covariance, $\Cov(\bx_t\mid\by_{1:t})\succeq\E[\Cov(\bx_t\mid\calG_t)\mid\by_{1:t}]$. Given $(\bx_{t-1},\bc_t)$, the past observations are independent of $\bg_t$, and $\bx_t=\bA\bx_{t-1}+\bc_t+\bg_t$, $\by_t=\bH_t\bx_t+\bv_t$ is linear-Gaussian in $\bg_t$. Hence $\Cov(\bx_t\mid\calG_t)=(\bSigma_g^{-1}+\bH_t^\top\bR_t^{-1}\bH_t)^{-1}\succeq(\underline\sigma_g^{-2}+h^2/\underline r)^{-1}\bI$.
\end{proof}

\rev{Proposition~\ref{prop:PE} relies on an unbounded Gaussian component and therefore does \emph{not} verify Assumption~\ref{ass:PE} under Assumption~\ref{ass:compact}; the two results address different regimes. Under compact dynamics, the following averaged property holds.}

\rev{\begin{proposition}
\label{prop:PEmean}
Let Assumption~\ref{ass:compact} hold, and let~\eqref{eq:obs} hold with $\bR_t\succeq\underline r\bI$, $\norm{\bH_t}\le h$, and $(\bv_t)$ independent of the states and residuals. Then, with $\lambda_0=(\bar J+h^2/\underline r)^{-1}$, for every $t\ge1$,
\begin{equation}
\E\big[\E[\bx_t\bx_t^\top\mid\by_{1:t}]\big]\succeq\E\big[\Cov(\bx_t\mid\by_{1:t})\big]\succeq\lambda_0\bI ,
\label{eq:PEmean}
\end{equation}
and, for every unit vector $\bm c$, $\Prob\big(\bm c^\top\Cov(\bx_t\mid\by_{1:t})\bm c\ge\lambda_0/2\big)\ge\lambda_0/(2R_x^2)$. [A proof is given in Appendix \ref{app_3}.]
\end{proposition}}

\rev{Assumption~\ref{ass:PE} strengthens~\eqref{eq:PEmean} to a uniform almost-sure statement over all times and observation sequences. We impose it directly as the classical persistent-excitation condition of adaptive filtering, formulated for the filtering posterior; deriving it from primitives, or replacing it by finite-horizon or high-probability versions, is left open.}

\begin{assumption}[Emulator accuracy]
\label{ass:score}
$\varepsilon_s:=\sup(\E_{\beps}\norm{\widehat\bs_\kappa(\bu_\kappa(\bw)\mid\bx)-\bs_\kappa(\bu_\kappa(\bw)\mid\bx)}^2)^{1/2}<\infty$ and $\sigma_m^2:=\sup\E_{\beps}\norm{\bXi(\bx)^\top\widehat\bs_\kappa(\bu_\kappa(\bw)\mid\bx)}^2<\infty$, where both suprema are over $\bw\in\calW_s$ and $\norm{\bx},\norm{\bx'}\le R_x$\rev{, and $\bs_\kappa$ is the smoothed score of the \emph{target} residual law}.
\end{assumption}

\rev{The supremum over $\bw$ covers mechanism (M1) of Section~\ref{sec:mismatch}, the candidate-query shift by $\bXi(\bx)(\bw-\bw^\star)$. Because $\bs_\kappa$ refers to the target law, $\varepsilon_s$ also absorbs (M2) source-detrending errors and conditional-mean misspecification, (M3) source--target residual shift and node-alignment errors, and the effect of plug-in estimates of $(\beta,\gamma)$. Assumption~\ref{ass:score} bounds their combined effect at the queried residuals; it does not assert that these mismatches are absent, and Proposition~\ref{prop:detrend} shows that some of them produce an irreducible topology bias.}

For pair measures we use the $L^2$ functional distance
\begin{equation}
\dd_{t}(\mu,\nu)=\sup_{\bm f}\big(\E\norm{\mu(\bm f)-\nu(\bm f)}^2\big)^{1/2},
\label{eq:dmetric}
\end{equation}
where the supremum ranges over $\R^M$-valued test functions on $\R^N\times\R^N$ with $\sup\norm{\bm f}\le1$ that are measurable with respect to the algorithm's history up to time $t-1$. For deterministic measures, $\dd_t\le2\TV$. Let $\Phi^{\bw}_{t+1}$ denote the exact pair update $\Gamma\mapsto\Gamma'$, $\Gamma'(d\bx,d\bx')\propto\Gamma^{(2)}(d\bx)P_{\bw}(\bx,d\bx')g_{t+1}(\by_{t+1}\mid\bx')$, where $\Gamma^{(2)}$ is the second marginal, and let $\Phi^\star_{s,t}$ denote the composition of the $\bw^\star$-updates.

\begin{assumption}
\label{ass:forget}
There are $C_0\ge1$ and $\varrho\in(0,1)$ such that $\dd_t(\Phi^\star_{s,t}\mu,\Phi^\star_{s,t}\nu)\le C_0\varrho^{t-s}\dd_s(\mu,\nu)$ for all $s\le t$ and admissible random pair measures. Moreover, $\bar\lambda_g:=\sup_t\norm{g_{t+1}}_\infty/\inf_{\bx}(P_{\bw^\star}g_{t+1})(\bx)<\infty$.
\end{assumption}

Uniform exponential forgetting of this type holds under the strong mixing condition of~\cite{delmoral2001stability,delmoral2004feynman}, via Hilbert-metric contraction~\cite{atar1997exponential}, and under weaker conditions~\cite{douc2009forgetting}. \rev{The likelihood-ratio bound is a statement about the realized observation sequence. Under Assumption~\ref{ass:compact} and Gaussian observation noise on the observed subspace with $\bR_t\succeq\underline r\bI$, $\norm{\bH_t}\le h$, one has $\norm{\by_{t+1}-\bH_{t+1}\bx'}\le2hR_x+\norm{\bv_{t+1}}$ for all $\norm{\bx'}\le R_x$, hence
\begin{equation}
\frac{\norm{g_{t+1}}_\infty}{\inf_{\bx}(P_{\bw^\star}g_{t+1})(\bx)}\le\exp\Big(\frac{(2hR_x+\norm{\bv_{t+1}})^2}{2\underline r}\Big).
\label{eq:lambdag}
\end{equation}
Since $\sup_t\norm{\bv_t}=\infty$ a.s.\ over an infinite horizon, $\bar\lambda_g<\infty$ is not implied by Gaussian noise. It holds for bounded observation noise and, on a finite horizon $T$, with probability at least $1-\delta$ with $\norm{\bv_{t}}\le\bar r^{1/2}(\sqrt N+\sqrt{2\log(T/\delta)})$ in~\eqref{eq:lambdag}, where $\bar r\ge\lambda_{\max}(\bR_t)$, by Gaussian concentration~\cite{boucheron2013concentration} and a union bound. The constants below then depend on $T$ and $\delta$.}

\subsection{Topology Sensitivity and Error Propagation}

\begin{lemma}
\label{lem:sens}
Under Assumptions~\ref{ass:compact} and~\ref{ass:forget}, for every pair measure $\Gamma$ supported in the ball and every $\bw\in\calW$,
\begin{equation}
\dd(\Phi^{\bw}_{t+1}\Gamma,\Phi^{\bw^\star}_{t+1}\Gamma)\le C_L\norm{\bw-\bw^\star},
\label{eq:CL}
\end{equation}
with $C_L=\sqrt2\,\bar\lambda_g\gamma\sqrt N R_x\bar J^{1/2}$.
\end{lemma}

\begin{proof}
By Lemma~\ref{lem:bayes} in Appendix~\ref{app:aux}, applied on the pair space with likelihood $g_{t+1}(\by\mid\bx')$, $\TV(\Phi^{\bw}\Gamma,\Phi^{\bw^\star}\Gamma)\le\bar\lambda_g\sup_{\bx}\TV(P_{\bw}(\bx,\cdot),P_{\bw^\star}(\bx,\cdot))$. The two kernels are translates of $\rho^\star(\cdot\mid\bx)$ by $\bm\delta=\bXi(\bx)(\bw-\bw^\star)$ by~\eqref{eq:shift}, so Lemma~\ref{lem:shiftTV} gives $\TV\le\frac12\norm{\bm\delta}\bar J^{1/2}\le\frac12\gamma\sqrt{2N}R_x\bar J^{1/2}\norm{\bw-\bw^\star}$. Multiplying by $2$ to pass from $\TV$ to $\dd$ gives $C_L=\sqrt2\bar\lambda_g\gamma\sqrt NR_x\bar J^{1/2}$, which is~\eqref{eq:CL}.
\end{proof}

Define the errors
\begin{equation}
e_t=\big(\E\norm{\bw_t-\bw^\star}^2\big)^{1/2},
\qquad
d_t=\dd_t(\widehat\Gamma_t,\Gamma_t^\star).
\label{eq:errors}
\end{equation}

\begin{lemma}
\label{lem:prop}
Under Assumptions~\ref{ass:compact} and~\ref{ass:forget}, with $\varepsilon_F$ from~\eqref{eq:epsF},
\begin{equation}
d_{t+1}\le C_0\varrho^{t+1}d_0+C_0\sum_{s=0}^{t}\varrho^{t-s}\big(\varepsilon_F+C_Le_s\big).
\end{equation}
Equivalently, $d_t\le Z_t$, where $Z_0=C_0d_0$ and $Z_{t+1}=\varrho Z_t+C_0(\varepsilon_F+C_Le_t)$.
\end{lemma}

\begin{proof}
Telescoping, $\widehat\Gamma_{t+1}-\Gamma^\star_{t+1}=\sum_{s=1}^{t+1}[\Phi^\star_{s,t+1}(\widehat\Gamma_s)-\Phi^\star_{s,t+1}(\Phi^\star_s\widehat\Gamma_{s-1})]+[\Phi^\star_{0,t+1}\widehat\Gamma_0-\Phi^\star_{0,t+1}\Gamma^\star_0]$. By Assumption~\ref{ass:forget} and the triangle inequality for $\dd$, the $s$-th term is at most $C_0\varrho^{t+1-s}\dd(\widehat\Gamma_s,\Phi^\star_s\widehat\Gamma_{s-1})$. Moreover, $\dd(\widehat\Gamma_s,\Phi^\star_s\widehat\Gamma_{s-1})\le\dd(\widehat\Gamma_s,\Phi^{\bw_{s-1}}_s\widehat\Gamma_{s-1})+\dd(\Phi^{\bw_{s-1}}_s\widehat\Gamma_{s-1},\Phi^\star_s\widehat\Gamma_{s-1})\le\varepsilon_F+C_Le_{s-1}$ by Proposition~\ref{prop:local}(ii), and by Lemma~\ref{lem:sens} applied pathwise and then in $L^2$. Reindexing gives the claim, and the recursion for $Z_t$ unrolls to the same expression.
\end{proof}

\subsection{Main Closed-Loop Theorem}

Let $\pi^\star_t$ denote the filter and $R_t:=R_{\kappa,\pi^\star_t}$. By Theorems~\ref{thm:ident}--\ref{thm:convex} and Assumption~\ref{ass:PE}, $\bw^\star$ minimizes every $R_t$ on $\R^M$, and on $\calW_s$ each $R_t$ is $\mu$-strongly convex and $L$-smooth with
\begin{equation}
\mu=
\begin{cases}
\dfrac{2\gamma^2\lambda_{\rm PE}(1-D^2/\kappa^2)}{\kappa^2}, & \kappa>D,\ \calW_s=\calW,\\[2mm]
\dfrac{\gamma^2\lambda_{\rm PE}}{D^2+\kappa^2}, & \calW_s\text{ as in~\eqref{eq:W2}},
\end{cases}
\label{eq:muLclosed}
\end{equation}
and $L=2N\gamma^2R_x^2/\kappa^2$.
In the second case, the trust radius of Corollary~\ref{cor:anneal} is computed with $\lambda_{\min}(\bSigma_\pi)$ replaced by $\lambda_{\rm PE}$ and $m_3(\pi)$ by $R_x^3$, which makes it uniform over $\pi=\pi^\star_t$, and the event of Corollary~\ref{cor:anneal} is assumed. The gradient constants are
\begin{equation}
G_\Xi=\gamma\sqrt{2N}R_x,
\qquad
\bar\Psi=\frac{G_\Xi\,(D+2R_x+\kappa\sqrt N)}{\kappa^2},
\label{eq:Psibar}
\end{equation}
where $\bar\Psi$ bounds $\norm{\E_{\beps}[\bXi(\bx)^\top\bs_\kappa(\bu_\kappa(\bw)\mid\bx)]}$ on the ball, since $\norm{\bs_\kappa(\bu)}=\norm{\E[\br\mid\bu]-\bu}/\kappa^2$ by~\eqref{eq:tweedie1}, $\norm{\E[\br\mid\bu]}\le D$, and $\norm{\bu}\le2R_x+\kappa\norm{\beps}$.

\begin{theorem}
\label{thm:closed}
Let Assumptions~\ref{ass:compact}--\ref{ass:forget} hold with $\bw^\star\in\calW_s$, and let $0<\eta\le\mu/L^2$. Define
\begin{gather}
K=\frac{6\eta\bar\Psi^2}{\mu},\quad
c_1=\frac{2C_0^2C_L^2}{1-\varrho},\quad
c_2=\frac{2C_0^2\varepsilon_F^2}{1-\varrho},\\
\sigma_n^2=2\bar\Psi^2+\frac{\sigma_m^2}{J_g},\quad
c_e=\frac{6\eta G_\Xi^2\varepsilon_s^2}{\mu}+\eta^2\sigma_n^2 .
\end{gather}
Then, componentwise,
\begin{equation}
\begin{bmatrix}e_{t+1}^2\\ Z_{t+1}^2\end{bmatrix}
\preceq
\underbrace{\begin{bmatrix}1-\frac{\eta\mu}2+Kc_1 & K\varrho\\ c_1 & \varrho\end{bmatrix}}_{\bG}
\begin{bmatrix}e_t^2\\ Z_t^2\end{bmatrix}
+
\begin{bmatrix}c_e+Kc_2\\ c_2\end{bmatrix}.
\label{eq:closedrec}
\end{equation}
If the small-gain condition
\begin{equation}
2\sqrt6\,\bar\Psi\,C_0C_L<\mu(1-\varrho)
\label{eq:smallgain}
\end{equation}
holds, then the spectral radius of $\bG$ is less than one, $e_t^2$ and $d_t^2\le Z_t^2$ converge geometrically to a neighborhood, and
\begin{align}
\limsup_{t\to\infty}e_t^2&\le\frac{1}{\widetilde\Delta}\Big[\frac{6(1-\varrho)G_\Xi^2}{\mu}\varepsilon_s^2+(1-\varrho)\eta\sigma_n^2\nonumber\\
&\quad+\frac{12\bar\Psi^2C_0^2}{\mu(1-\varrho)}\varepsilon_F^2\Big],
\label{eq:efloor}\\
\limsup_{t\to\infty}d_t^2&\le\frac{2C_0^2}{(1-\varrho)\widetilde\Delta}\Big[\frac{\mu}{2}\varepsilon_F^2\nonumber\\
&\quad+C_L^2\Big(\frac{6G_\Xi^2}{\mu}\varepsilon_s^2+\eta\sigma_n^2\Big)\Big],
\label{eq:dfloor}
\end{align}
with $\widetilde\Delta=\big(\mu^2(1-\varrho)^2-24\bar\Psi^2C_0^2C_L^2\big)/\big(2\mu(1-\varrho)\big)>0$. [A proof is given in Appendix \ref{app_F}.]
\end{theorem}

\begin{remark}
Condition~\eqref{eq:smallgain} compares three gains: $\bar\Psi$ (state error $\to$ topology gradient), $C_0C_L/(1-\varrho)$ (topology error $\to$ accumulated state error), and $\mu$ (topology contraction, i.e., Fisher curvature times posterior excitation). The step size cancels, so stability is a property of the model, the observation design, and the noise level, not of tuning. The error floor~\eqref{eq:efloor} separates the score error $\varepsilon_s$, the local filtering error $\varepsilon_F$ (particles, \rev{proposal total-variation error}, ancestor weights), and the stochastic-approximation variance $\eta\sigma_n^2$, which vanishes as $\eta\to0$. No bias arises from conditioning on the observations, because the conditional mean gradient $\pi^\star_t(\bphi_{\bw})$ is the gradient of an objective minimized at $\bw^\star$ for \emph{every} state law (Theorem~\ref{thm:ident}(b)); the remaining posterior-sampling fluctuation $\bm{\xi}_t$ is a martingale difference. Hence no ergodicity argument for the gradient noise is needed.
\end{remark}

\begin{remark}
Since $\mu\propto\kappa^{-2}(1-D^2/\kappa^2)$ and $\bar\Psi\propto\kappa^{-2}(D+2R_x+\kappa\sqrt N)$ for $\kappa>D$, the ratio $\bar\Psi/\mu$ is proportional to $(D+2R_x+\kappa\sqrt N)/(1-D^2/\kappa^2)$, which is minimized at a moderate multiple of $D$. Small-gain stability therefore favors $\kappa$ just above the convexity threshold, efficiency favors small $\kappa$ (Theorem~\ref{thm:eff}), and the trust-region stages of Corollary~\ref{cor:anneal} reconcile the two.
\end{remark}

\rev{\begin{remark}
\label{rem:closed-scope}
Theorem~\ref{thm:closed} is conditional on bounded process residuals, uniform posterior excitation, uniform filter forgetting with a bounded likelihood ratio, and controlled score and proposal errors. These conditions are not implied by Gaussian or heavy-tailed residual models, nor by empirical quantile clipping. We therefore interpret the theorem as a stability guarantee for a well-specified bounded-residual regime, and we examine robustness beyond this regime empirically. The error floor contains $\varepsilon_s$, $\varepsilon_F$, and $\eta\sigma_n^2$; errors from plug-in estimates of $(\beta,\gamma)$, node misalignment, and source--target shift enter only through $\varepsilon_s$, and the clipping bias of Proposition~\ref{prop:local}(iii) is not covered outside Assumption~\ref{ass:compact}.
\end{remark}}

\begin{corollary}
\label{cor:support}
Under the conditions of Theorem~\ref{thm:closed}, with $\bar e^2$ the right-hand side of~\eqref{eq:efloor}, \rev{for every threshold $\vartheta\in(0,w^\star_{\min})$ and} all sufficiently large $t$,
\[
\rev{\Prob\big(\widehat\calS_t(\vartheta)\ne\calS^\star\big)\le\frac{\bar e^2+o(1)}{\min\{\vartheta,\,w^\star_{\min}-\vartheta\}^2},}
\]
where \rev{$\widehat\calS_t(\vartheta)=\{e:w_{t,e}>\vartheta\}$. For $\vartheta=w^\star_{\min}/2$ the bound reads $4(\bar e^2+o(1))/(w^\star_{\min})^2$.}
\end{corollary}

\begin{proof}
\rev{If $\norm{\bw_t-\bw^\star}_\infty<\min\{\vartheta,w^\star_{\min}-\vartheta\}$, then every non-edge has $w_{t,e}<\vartheta$ and every edge has $w_{t,e}>w^\star_{\min}-(w^\star_{\min}-\vartheta)=\vartheta$.} Apply Markov's inequality to $\norm{\bw_t-\bw^\star}^2\ge\norm{\bw_t-\bw^\star}_\infty^2$\rev{, and use $e_t^2\le\bar e^2+o(1)$ from~\eqref{eq:efloor}}.
\end{proof}

\rev{The threshold $w^\star_{\min}/2$ uses knowledge of the true graph and is an oracle choice; Corollary~\ref{cor:support} shows that any $\vartheta$ bounded away from $0$ and from $w^\star_{\min}$ suffices, but a deployable rule must select $\vartheta$, or the number of edges, from data.}


\section{Numerical Experiments}
\label{sec:exp}

The experiments test the predictions of Sections~\ref{sec:theory}--\ref{sec:closed} one at a time: moment identifiability (E1), convexification by diffusion noise (E2), statistical efficiency (E3), and closed-loop topology and state estimation (E4), followed by ablations, a stability stress test, and a complexity study. All quantities are averaged over independent runs, and uncertainty is reported as one standard deviation in the tables and as $95\%$ confidence bands in the figures.

\subsection{Synthetic Benchmarks}
\label{sec:exp-data}


\emph{Graphs.}
We draw static undirected graphs from three
families that differ in degree heterogeneity and
spatial locality: Erd\H{o}s--R\'enyi (ER) with
edge probability $4/N$, Barab\'asi--Albert (BA)
with attachment parameter $m=2$, and random
geometric graphs (RGG) on the unit square.
Only connected realizations are retained.
Edge weights are drawn as
$w^\star_e\sim\mathcal{U}[0.5,1.5]$ on
$\calS^\star$ and rescaled so that
$\gamma\lambda_{\max}(\bL^\star)
=0.8(1+\beta)$.
Thus, $\bw^\star$ lies on the boundary of
$\calW$ due to its zero edge weights, while
the nonzero weights are in the relative interior
of the oracle-support feasible set
$\calW_{\calS^\star}$.
We set $\beta=0.95$, $\delta_A=0.05$,
and $N\in\{20,50,100\}$, with $N=200$
added for the complexity study.

\emph{Residual laws.} To control the degree of non-Gaussianity, we consider five residual laws, each normalized to unit marginal variance: (i) Gaussian; (ii) Laplace; (iii) Student-$t_3$; (iv) bimodal, $c\in\{\pm a\}$ plus $\calN(0,0.25^2)$ noise, which satisfies Assumption~\ref{ass:resid} with $D=a$; and (v) switching forcing, $\bc=s(\bm b-q\one)$ with $b_i\sim\mathrm{Bern}(q)$, $q=0.04$, $s=3$, plus $\calN(\bzero,0.3^2\bI)$. Law (v) mimics intermittent rainfall-driven inflows. For E2 we additionally use a lattice residual, in which $\bc$ is uniform on $2^k$ separated points within a $D$-ball, $k\in\{1,2,3\}$, $D=1$, and $\bar\sigma_g=0.15$.

\emph{Observations.} States are observed through the masked model~\eqref{eq:maskobs} with $m_{t,i}\sim\mathrm{Bern}(p)$, independent across nodes and time, and $\bR_v=0.1^2\bI$. The sampling rate $p$ ranges over $\{0.1,\ldots,0.9\}$, with $p=0.3$ as the default.

\emph{Emulator and filter.} The residual denoiser $D_\theta$ is a four-block residual MLP of width $256$, conditioned on $\bx_t$ and on the noise level, and trained with the EDM preconditioning~\cite{karras2022edm}. Training uses detrended residuals $\bx_{t+1}-\bA_{\rm src}\bx_t$ from eight source graphs that are disjoint from the test graphs, so that topology learning at test time is genuinely transductive. Unless stated otherwise, TopoScore-FAF uses $J=256$ particles, $J_g=64$ gradient pairs, $K=32$ sampler steps, $n_{\rm CG}=3$, $K_\epsilon=8$ noise draws per pair, and step size $\eta=2\times10^{-3}$. The annealing schedule is $\kappa\in\{1.5D,0.5D,0.15D\}$ with switches at $t=1500$ and $t=3000$, and the default horizon is $T=5000$. Table~\ref{tab:setup} summarizes the configuration.

\emph{Metrics.} Topology accuracy is measured by the relative Laplacian error $\epsilon_L=\norm{\bL(\widehat\bw)-\bL^\star}_F/\norm{\bL^\star}_F$, 
by the edge F1 score of
$\widehat\calS=
\{e:\widehat w_e>w^\star_{\min}/2\}$,
using the oracle threshold of
Corollary~\ref{cor:support} solely for
synthetic evaluation, and by the threshold-free area under the precision--recall curve (AUPRC). State estimation is measured by the nRMSE of the posterior mean, normalized by the stationary state standard deviation, and by the continuous ranked probability score (CRPS)~\cite{gneiting2007strictly} under the same normalization. Calibration is measured by the spread--skill ratio (SSR)~\cite{fortin2014spread}, whose ideal value is $1$, and particle degeneracy by the effective sample size relative to $J$ (ESS$/J$). The first $50$ steps are discarded as burn-in.

\begin{table*}[!t]
\centering
\caption{Synthetic benchmark configuration. The last column gives the default used in E4.}
\label{tab:setup}
\footnotesize
\setlength{\tabcolsep}{4pt}
\begin{tabular}{@{}lll@{}}
\toprule
Component & Setting & Default \\
\midrule
Graph families & ER ($p_{\rm edge}=4/N$), BA ($m=2$), RGG (connected realizations) & ER \\
Number of nodes & $N\in\{20,50,100\}$; $N=200$ for the complexity study & $50$ \\
Edge weights & $w^\star_e\sim\mathcal{U}[0.5,1.5]$, rescaled to $\gamma\lambda_{\max}(\bL^\star)=0.8(1+\beta)$ & -- \\
Dynamics & $\beta=0.95$, $\delta_A=0.05$ & -- \\
\midrule
Residual laws (unit variance) & Gaussian, Laplace, Student-$t_3$, bimodal, switching forcing ($q=0.04$, $s=3$) & switching \\
Lattice residual (E2) & $2^k$ points in a $D$-ball, $k\in\{1,2,3\}$, $D=1$, $\bar\sigma_g=0.15$ & -- \\
\midrule
Observations & Bernoulli masks with rate $p\in\{0.1,\ldots,0.9\}$, $\bR_v=0.1^2\bI$ & $p=0.3$ \\
Horizon & $T\in[10^3,2\times10^5]$ (E1); $n=5000$ oracle pairs (E3); $T=5000$ (E4) & $5000$ \\
\midrule
Emulator & 4-block MLP, width 256, EDM schedule, 8 source graphs & -- \\
TopoScore-FAF & $J=256$, $J_g=64$, $K=32$, $n_{\rm CG}=3$, $K_\epsilon=8$, $\eta=2\times10^{-3}$ & -- \\
Annealing & $\kappa\in\{1.5D,0.5D,0.15D\}$, switches at $t=1500,3000$ & -- \\
Repetitions & 20 runs (E1, E4), 200 initializations (E2), 500 Monte Carlo replicates (E3) & -- \\
\bottomrule
\end{tabular}
\end{table*}

\subsection{Baselines and Protocol}
\label{sec:exp-baselines}

The baselines cover both halves of the problem and include classical as well as recent methods.

\emph{Topology inference.} We compare against the following methods:
\begin{itemize}
\item smoothness-based graph learning (Kalofolias)~\cite{kalofolias2016graph};
\item spectral templates (SpecTemp)~\cite{segarra2017network};
\item least squares on graph VAR dynamics (Graph-VARLS)~\cite{isufi2019forecasting};
\item the moment estimator of Theorem~\ref{thm:ident}(a);
\item Kalman-EM, i.e., EM on a Gaussian graph state-space model;
\item GKNet~\cite{sabbaqi2026gknet}, with its graph treated as a learnable parameter.
\end{itemize}
Kalofolias, SpecTemp, and Graph-VARLS require complete graph signals, so we apply them to Kalman-smoothed states. In the oracle-pair study (E3), we additionally include two robust regression estimators, the Huber M-estimator~\cite{huber1964robust} and least absolute deviations (LAD)~\cite{koenker1978regression}, which are natural non-Gaussian competitors to least squares.

\emph{State estimation.} We compare against the following filters:
\begin{itemize}
\item the graph Kalman filter (Graph KF)~\cite{isufi2020observing};
\item the bootstrap particle filter~\cite{gordon1993novel};
\item KalmanNet~\cite{revach2022kalmannet};
\item GKNet~\cite{sabbaqi2026gknet};
\item FlowDAS~\cite{chen2025flowdas};
\item the training-free FA-APF with MMPS guidance (FA-APF)~\cite{savary2026trainingfree,rozet2024mmps}.
\end{itemize}
Methods that cannot learn the graph are given the moment estimate $\bw_{\rm mom}$. FA-APF and Graph KF are also run with the true graph, which provides oracle-topology references.

\emph{Fairness.} All methods that estimate a graph start from the same projected moment estimate. Hyperparameters of every method are tuned on a separate validation trajectory. All generative filters use the same emulator architecture and particle budget, and no method sees the test graph during training. Runtimes are measured on one NVIDIA A100 GPU with batched particle evaluation.

\subsection{E1: Identifiability and Warm Start}
\label{sec:exp-e1}

\begin{figure*}[!t]
\centering
\includegraphics[width=\textwidth]{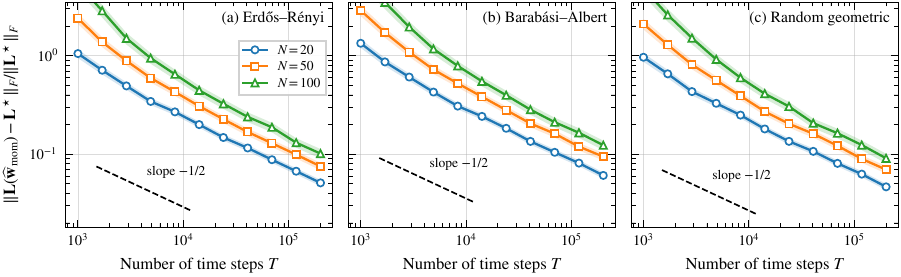}
\caption{E1: relative Laplacian error of the moment estimator versus the horizon $T$ ($p=0.5$) for (a) ER, (b) BA, and (c) RGG graphs. The dashed segment indicates slope $-1/2$; shaded bands are $95\%$ confidence intervals over $20$ runs.}
\label{fig:e1-T}
\end{figure*}

Fig.~\ref{fig:e1-T} reports the relative error of $\bL(\widehat\bw_{\rm mom})$ as a function of the horizon $T$ at $p=0.5$. For every graph family and size, the error
exhibits an approximate $T^{-1/2}$ decay
once $T\gtrsim\res{10^4}$, consistent with
the usual parametric rate under suitable
mixing conditions and with the consistency
result of Theorem~\ref{thm:ident}(a). The steeper decay at short horizons reflects the finite-sample bias of inverting $\widehat\bSigma_x$. The error grows with $N$ because the number of moment entries grows quadratically while each entry is estimated from the same samples. BA graphs are the hardest of the three families, with errors about $\res{20\%}$ larger than for ER graphs, which we attribute to the poorer conditioning of $\bSigma_x$ around high-degree hubs.

\begin{figure}[!t]
\centering
\includegraphics[width=\columnwidth]{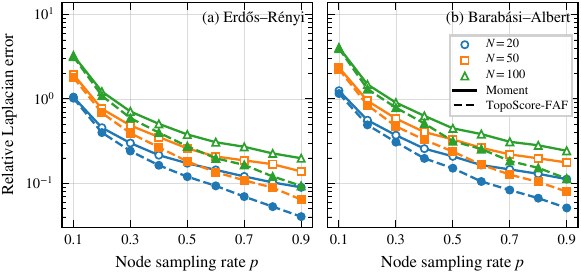}
\caption{E1: relative Laplacian error versus the node sampling rate $p$ ($T=2\times10^4$) for (a) ER and (b) BA graphs. Solid lines: moment estimator; dashed lines: TopoScore-FAF initialized at the moment estimate.}
\label{fig:e1-p}
\end{figure}

Fig.~\ref{fig:e1-p} varies the sampling rate at $T=2\times10^4$. Since each moment in~\eqref{eq:moments} is estimated from a fraction $\pi^{(0)}_{ij}=p^2$ of the time steps, the moment error scales roughly as $p^{-1}$. TopoScore-FAF, initialized at $\bw_{\rm mom}$, improves on it at every rate, and the improvement widens with coverage: the ratio of the two errors falls from about $\res{0.9}$ at $p=0.1$ to about $\res{0.45}$ at $p=0.9$. The reason is that sparse sampling degrades the posterior state pairs that feed the score update, whereas the moment estimator is affected only through its sample size. The moment estimator is therefore an accurate warm start, but it does not by itself reach the accuracy achievable from the same data.

\subsection{E2: Convexification by Diffusion Noise}
\label{sec:exp-e2}

\begin{figure}[!t]
\centering
\includegraphics[width=\columnwidth]{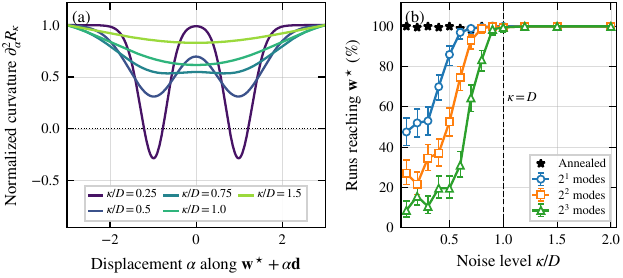}
\caption{E2: (a) normalized curvature of $R_\kappa$ along a random direction through $\bw^\star$ for the bimodal residual; (b) fraction of random initializations reaching $\bw^\star$ versus $\kappa/D$ for lattice residuals with $2^k$ modes (Wilson $95\%$ intervals). Stars: annealed schedule.}
\label{fig:e2}
\end{figure}

Fig.~\ref{fig:e2}(a) shows the curvature $\partial_\alpha^2R_\kappa(\bw^\star+\alpha\bm d)$ along a random direction $\bm d$ for the bimodal residual. At $\kappa/D=\res{0.25}$, the curvature becomes negative around $|\alpha|\approx\res{1}$, the displacement at which the shifted residual places mass between the two modes. For $\kappa/D\ge\res{0.5}$, the curvature is positive along the whole line. This agrees with Theorem~\ref{thm:convex}(a), which guarantees strong convexity only for $\kappa>D$. The empirical threshold is lower because the bound $\Cov(\bc\mid\bu,\bx)\preceq D^2\bI$ used in the proof is attained only in the worst case.

Fig.~\ref{fig:e2}(b) reports the fraction of $200$ uniformly drawn initializations in $\calW$ from which projected gradient descent on $R_{\kappa}$ reaches $\bw^\star$ (final error below $10^{-2}$). At small noise levels the success rate is low and decreases with the number of modes: at $\kappa/D=\res{0.1}$ it is $\res{47\%}$, $\res{27\%}$, and $\res{8\%}$ for $2^1$, $2^2$, and $2^3$ modes. All runs succeed once $\kappa\ge D$, in agreement with the theory. The transition moves to larger $\kappa$ as the number of modes increases, because more modes create more spurious stationary points that must be smoothed out. Annealing from $\kappa_1=1.5D$ with the trust region of Corollary~\ref{cor:anneal} reaches $\bw^\star$ in $\res{\ge99\%}$ of the runs at every final noise level. Hence the efficiency of small $\kappa$ (Section~\ref{sec:exp-e3}) can be obtained without the initialization sensitivity of a direct small-noise optimization. 

\subsection{E3: Statistical Efficiency with Oracle State Pairs}
\label{sec:exp-e3}

\begin{figure*}[!t]
\centering
\includegraphics[width=\textwidth]{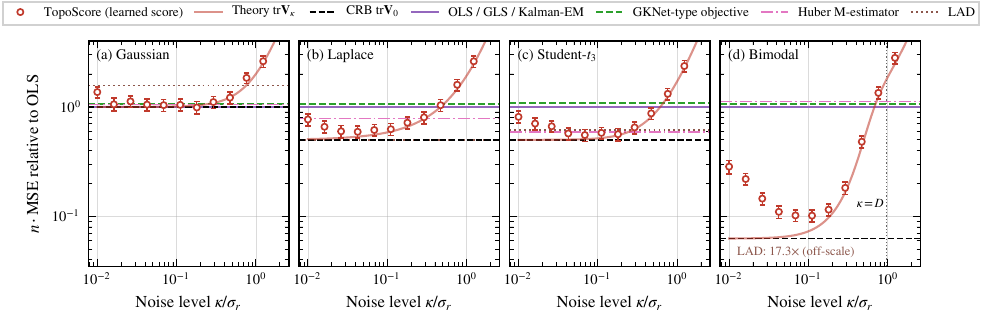}
\caption{E3: $n\cdot\mathrm{MSE}$ of the topology estimator relative to OLS versus the noise level $\kappa$, with oracle state pairs ($n=5000$, $N=20$, $500$ replicates), for (a) Gaussian, (b) Laplace, (c) Student-$t_3$, and (d) bimodal residuals. Solid curves: single-noise-draw oracle
information benchmarks $\tr\bV_\kappa$
computed from the exact smoothed Fisher
information; markers: TopoScore with the learned score ($95\%$ intervals); horizontal lines: baselines and the Cram\'er--Rao bound.}
\label{fig:e3}
\end{figure*}

To isolate the statistical value of the score
from filtering error, E3 estimates $\bw$ from
$n=5000$ true state pairs with $N=20$.
Fig.~\ref{fig:e3} plots $n\cdot\mathrm{MSE}$,
normalized by that of OLS, as a function of
$\kappa$.
The figure also reports
$\tr\bV_\kappa$, computed from the exact
smoothed Fisher information by numerical
quadrature, and the zero-noise reference
$\tr\bV_0$.
These curves serve as oracle information
benchmarks rather than exact finite-sample
covariance predictions for the constrained
estimator with a learned score.
In particular, Theorem~\ref{thm:eff}
concerns the known-support submodel,
Proposition~\ref{prop:boundary} characterizes
the unknown-support boundary limit, and
Corollary~\ref{cor:Kdraw} accounts for
averaging multiple diffusion-noise draws.
Table~\ref{tab:e3} reports the corresponding
empirical NMSE values.

The results illustrate the statistical
trade-offs underlying Theorem~\ref{thm:eff},
while also revealing the effects of learned
scores and finite-sample estimation:

\begin{enumerate}
\item \emph{Monotonicity.} The theoretical and empirical curves increase with $\kappa$. At $\kappa=\sigma_r$, the estimator is about $\res{2}$ times less efficient than OLS (Table~\ref{tab:e3}, fixed $\kappa=\sigma_r$). A large noise level therefore buys convexity at a substantial statistical cost.
\item \emph{Cram\'er--Rao limit.} For Gaussian residuals, $\bV_0=\bV_{\rm GLS}(0)$, so no gain is possible. Accordingly, the annealed estimator matches OLS and GLS to within $\res{0.3}$~dB, which is within one standard deviation.
\item \emph{Strict dominance.} The gain over Gaussian quasi-likelihood grows with the Fisher non-Gaussianity of the residual. It is $\res{2.2}$~dB for Laplace, $\res{2.5}$~dB for Student-$t_3$, and $\res{9.9}$~dB for the bimodal law. In all cases the remaining gap to the Cram\'er--Rao bound is at most $\res{2.1}$~dB.
\item \emph{Bias of the learned score at small noise.} At very small $\kappa$, the learned estimator departs from $\tr\bV_\kappa$, and the departure is largest for the bimodal law. The denoiser is least accurate in this regime, where the score is steep between modes. The resulting $\varepsilon_s$ inflates the variance, and the empirical optimum lies at $\kappa\approx\res{0.05\text{--}0.1}\,\sigma_r$ rather than at $\kappa\to0$. For this reason, the annealing schedule stops at a small but finite noise level.
\end{enumerate}

The robust baselines provide an informative contrast. LAD is the maximum-likelihood estimator under Laplace noise and is therefore efficient there. It outperforms TopoScore by $\res{0.7}$~dB because it does not incur emulator error. Similarly, Huber estimation comes within $\res{0.1}$~dB of TopoScore for Student-$t_3$ residuals, which is well within one standard deviation. However, both robust estimators assume unimodal, symmetric noise. Under the bimodal law, the residual density vanishes near the origin, LAD degrades by $\res{12.4}$~dB relative to OLS, and Huber estimation is slightly worse than OLS. TopoScore is the only estimator that adapts to the residual law without such an assumption. Finally, the GKNet-type objective, whose edge-structured working covariance is misspecified here, is slightly less efficient than GLS (about $\res{0.3}$~dB), consistent with the sandwich bound $\bV_{\bm W}\succeq\bV_{\rm GLS}(0)$.

\begin{table*}[!t]
\centering
\caption{E3: NMSE (dB, mean $\pm$ std over $500$ replicates) of the topology estimate from $n=5000$ oracle state pairs ($N=20$). Best values are in bold and second-best values are underlined; the Cram\'er--Rao bound is excluded from the ranking.}
\label{tab:e3}
\footnotesize
\setlength{\tabcolsep}{5pt}
\begin{tabular}{@{}llcccc@{}}
\toprule
Estimator & Type & Gaussian & Laplace & Student-$t_3$ & Bimodal \\
\midrule

OLS / Graph-VARLS~\cite{isufi2019forecasting} & Classical & $\mathbf{-24.60\pm0.47}$ & $-24.60\pm0.49$ & $-24.60\pm0.48$ & $-24.60\pm0.44$ \\
Moment (Theorem~\ref{thm:ident}(a)) & Classical & $-24.50\pm0.48$ & $-24.50\pm0.46$ & $-24.40\pm0.44$ & $-24.44\pm0.44$ \\
GLS / Kalman-EM & Classical & $\underline{-24.60\pm0.47}$ & $-24.60\pm0.49$ & $-24.60\pm0.48$ & $-24.60\pm0.44$ \\
Huber M-estimator~\cite{huber1964robust} & Classical & $-24.42\pm0.44$ & $-25.68\pm0.43$ & $\underline{-26.98\pm0.47}$ & $-24.14\pm0.42$ \\
LAD~\cite{koenker1978regression} & Classical & $-22.63\pm0.47$ & $\mathbf{-27.51\pm0.46}$ & $-26.72\pm0.46$ & $-12.24\pm2.90$ \\
GKNet-type objective~\cite{sabbaqi2026gknet} & Recent & $-24.34\pm0.45$ & $-24.34\pm0.42$ & $-24.29\pm0.49$ & $-24.27\pm0.43$ \\
TopoScore, fixed $\kappa=\sigma_r$ & Ours & $-21.64\pm0.47$ & $-21.90\pm0.48$ & $-22.27\pm0.43$ & $-22.01\pm0.46$ \\
TopoScore, fixed $\kappa=0.01\sigma_r$ & Ours & $-23.56\pm0.42$ & $-25.70\pm0.45$ & $-25.40\pm0.45$ & $\underline{-29.53\pm0.43}$ \\
TopoScore, annealed & Ours & $-24.34\pm0.47$ & $\underline{-26.79\pm0.47}$ & $\mathbf{-27.11\pm0.47}$ & $\mathbf{-34.51\pm0.43}$ \\
\midrule
Cram\'er--Rao bound $\tr\bV_0$ & Bound & $-24.60$ & $-27.61$ & $-27.61$ & $-36.63$ \\
\bottomrule
\end{tabular}
\end{table*}

\subsection{E4: Closed-Loop Topology and State Estimation}
\label{sec:exp-e4}

E4 evaluates the complete algorithm under partial observation ($N=50$, $p=0.3$, $T=5000$, switching-forcing residuals) over $20$ runs.

\begin{figure*}[!t]
\centering
\includegraphics[width=\textwidth]{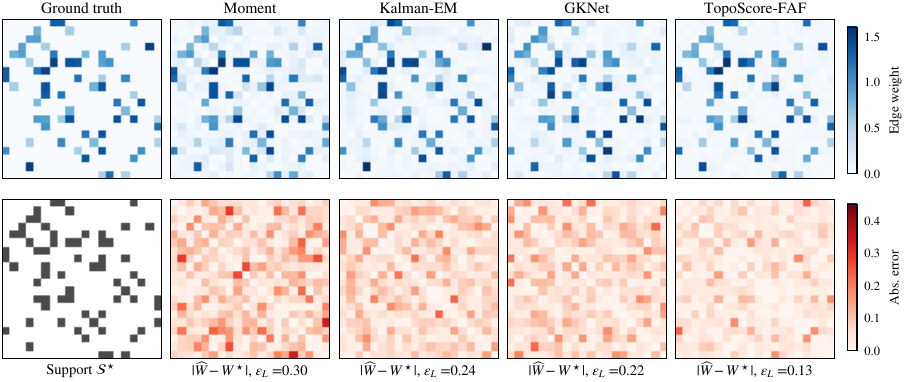}
\caption{Estimated weighted adjacency matrices (top) and absolute errors (bottom) on an ER graph with $N=20$ ($p=0.3$, switching forcing). The leftmost column shows the true weights and support.}
\label{fig:vis-adj}
\end{figure*}

\emph{Topology recovery.} Table~\ref{tab:e4-topo} shows that the methods fall into three groups:
\begin{itemize}
\item \emph{Static-signal methods.} Kalofolias and SpecTemp assume smooth or stationary-diffusion signals rather than graph-diffusion dynamics, and they yield $\epsilon_L>\res{0.6}$ on every graph family.
\item \emph{Second-order dynamic methods.} Graph-VARLS, the moment estimator, Kalman-EM, and GKNet exploit the dynamics, but only through Gaussian or second-order objectives. Among them, GKNet is the strongest, with $\epsilon_L=\res{0.221}$ on ER graphs.
\item \emph{TopoScore-FAF.} It reduces the error to $\res{0.129}$ on ER graphs, a $\res{42\%}$ relative reduction over GKNet, and raises the F1 score from $\res{0.863}$ to $\res{0.925}$.
\end{itemize}
The improvement is consistent across graph families. On BA graphs, all methods deteriorate, because hubs concentrate the edge excitation on few nodes and reduce $\lambda_{\min}(\bSigma_\pi)$ in Lemma~\ref{lem:excitation}; TopoScore-FAF nevertheless remains best ($\res{0.158}$ versus $\res{0.267}$). The gap between TopoScore-FAF and its oracle-state counterpart ($\res{0.129}$ versus $\res{0.069}$) quantifies the cost of latent states. This cost corresponds to the terms $\varepsilon_F$ and $\bar\Psi$ in the error floor~\eqref{eq:efloor}. Fig.~\ref{fig:vis-adj} visualizes the estimates on an ER graph with $N=20$. The second-order estimators spread spurious weight over non-edges. TopoScore-FAF concentrates its residual errors on a few true edges and is close to zero elsewhere, which is why it attains the highest F1 score.

\begin{table*}[!t]
\centering
\caption{E4: topology recovery under partial observation ($N=50$, $p=0.3$, $T=5000$, switching forcing; mean over $20$ runs, $\pm$ std for $\epsilon_L$). Best values are in bold and second-best values are underlined; the oracle-state reference is excluded from the ranking.}
\label{tab:e4-topo}
\footnotesize
\setlength{\tabcolsep}{3.2pt}
\begin{tabular}{@{}lccccccccc@{}}
\toprule
\multirow{2}{*}{Method}  & \multicolumn{3}{c}{Erd\H{o}s--R\'enyi} & \multicolumn{3}{c}{Barab\'asi--Albert} & \multicolumn{3}{c}{Random geometric} \\
\cmidrule(lr){2-4}\cmidrule(lr){5-7}\cmidrule(l){8-10}
 & $\epsilon_L\downarrow$ & F1$\uparrow$ & AUPRC$\uparrow$ & $\epsilon_L\downarrow$ & F1$\uparrow$ & AUPRC$\uparrow$ & $\epsilon_L\downarrow$ & F1$\uparrow$ & AUPRC$\uparrow$ \\
\midrule
Kalofolias~\cite{kalofolias2016graph}  & $0.712\pm0.071$ & $0.588$ & $0.609$ & $0.767\pm0.089$ & $0.535$ & $0.576$ & $0.691\pm0.057$ & $0.600$ & $0.620$ \\
SpecTemp~\cite{segarra2017network}  & $0.631\pm0.074$ & $0.643$ & $0.672$ & $0.697\pm0.061$ & $0.590$ & $0.635$ & $0.616\pm0.069$ & $0.660$ & $0.684$ \\
Graph-VARLS~\cite{isufi2019forecasting}  & $0.396\pm0.042$ & $0.744$ & $0.801$ & $0.465\pm0.041$ & $0.692$ & $0.762$ & $0.364\pm0.033$ & $0.757$ & $0.818$ \\
Moment estimator  & $0.304\pm0.028$ & $0.790$ & $0.860$ & $0.352\pm0.038$ & $0.743$ & $0.825$ & $0.273\pm0.024$ & $0.806$ & $0.867$ \\
Kalman-EM  & $0.244\pm0.025$ & $0.845$ & $0.906$ & $0.287\pm0.034$ & $0.810$ & $0.873$ & $0.221\pm0.026$ & $0.853$ & $0.915$ \\
GKNet~\cite{sabbaqi2026gknet} & \underline{$0.221\pm0.019$} & \underline{$0.863$} & \underline{$0.911$} & \underline{$0.267\pm0.022$} & \underline{$0.823$} & \underline{$0.875$} & \underline{$0.212\pm0.023$} & \underline{$0.872$} & \underline{$0.918$} \\
TopoScore-FAF (ours)  & $\mathbf{0.129\pm0.013}$ & $\mathbf{0.925}$ & $\mathbf{0.957}$ & $\mathbf{0.158\pm0.015}$ & $\mathbf{0.878}$ & $\mathbf{0.930}$ & $\mathbf{0.124\pm0.011}$ & $\mathbf{0.947}$ & $\mathbf{0.964}$ \\
\midrule
TopoScore, oracle states & $0.069\pm0.006$ & $0.982$ & $0.992$ & $0.083\pm0.009$ & $0.934$ & $0.959$ & $0.065\pm0.007$ & $0.995$ & $0.998$ \\
\bottomrule
\end{tabular}
\end{table*}

\begin{figure*}[!t]
\centering
\includegraphics[width=0.74\textwidth]{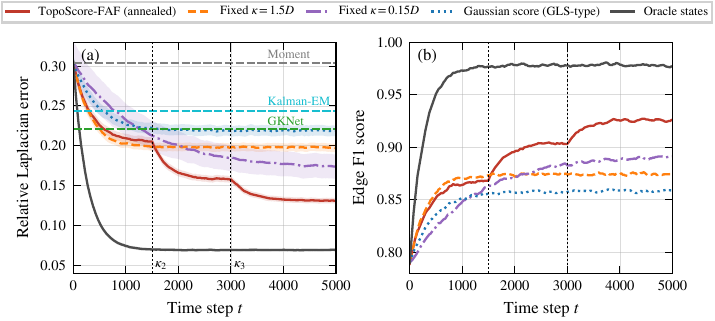}
\caption{E4: (a) relative Laplacian error and (b) edge F1 score along the trajectory ($N=50$, ER, $p=0.3$, switching forcing). Vertical lines mark the noise-level switches; horizontal lines give the final values of the baselines. Bands are $95\%$ confidence intervals over $20$ runs.}
\label{fig:e4-conv}
\end{figure*}

\emph{Convergence of the annealed schedule.} Fig.~\ref{fig:e4-conv} traces the topology error along the trajectory. During Stage~I ($\kappa=1.5D$), the iterates contract quickly to the noise-limited plateau $\epsilon_L\approx\res{0.20}$, which coincides with the curve obtained with $\kappa=1.5D$ held fixed. Each subsequent reduction of $\kappa$ produces a further drop, to $\res{0.158}$ and then to $\res{0.131}$. The plateaus therefore decrease monotonically, as predicted by Theorem~\ref{thm:eff}(b), and no transient increase occurs at the switching times, because the minimizer does not move with $\kappa$ (Theorem~\ref{thm:ident}(b)). Running directly at $\kappa=0.15D$ converges more slowly and reaches a worse and more variable error ($\res{0.17}$, with wide confidence bands): this run starts in the non-convex regime, and $\res{3}$ of the $20$ runs did not converge. Replacing the learned score by a Gaussian score turns the update into closed-loop GLS. This variant stalls at $\res{0.22}$, the level of GKNet, which isolates the gain due to the non-Gaussian score.

\emph{State estimation.} Table~\ref{tab:e4-state} compares the filters under three residual laws. Two observations frame the comparison.
\begin{itemize}
\item \emph{The non-Gaussian gain in nRMSE is small.} Even with the exact transition, FA-APF does not improve on the Graph KF in nRMSE ($\res{0.711}$ versus $\res{0.704}$ under switching forcing). In this regime, most of the error comes from unobserved nodes, about whose forcing a causal filter receives little information. Non-Gaussianity mainly affects the shape of the posterior. Accordingly, its effect appears in CRPS rather than in nRMSE: FA-APF improves the CRPS of the Graph KF by $\res{2.3\%}$ under switching forcing, but loses to it slightly under Gaussian residuals, where the Kalman filter is exact.
\item \emph{Graph accuracy shifts all filters only moderately.} With the moment-estimated graph instead of the true graph, FA-APF loses about $\res{3\%}$ in nRMSE.
\end{itemize}

Within this narrow range, TopoScore-FAF is among the best methods with an unknown topology on all three laws. Its nRMSE and CRPS are within $\res{0.007}$ and $\res{0.002}$, respectively, of those of the oracle-topology FA-APF; this residual gap reflects the remaining topology error. Under switching forcing, it attains the lowest CRPS ($\res{0.301}$ versus $\res{0.313}$--$\res{0.332}$) and the best calibration among the learned-topology methods (SSR $\res{0.96}$). Its nRMSE, however, is statistically indistinguishable from that of Kalman-EM and GKNet ($\res{0.718}$ versus $\res{0.717}$). The particle and learned-gain baselines that use the moment graph are consistently under-dispersed: SSR is $\res{0.76}$ for the bootstrap filter and $\res{0.85}$ for KalmanNet. The bootstrap filter additionally suffers from weight degeneracy (Section~\ref{sec:exp-abl}). TopoScore-FAF is also the most expensive method per step ($\res{58}$~ms versus $\res{1.3}$~ms for GKNet). Its advantage therefore lies primarily in topology recovery and probabilistic calibration rather than in point-estimation accuracy, a trade-off we state explicitly.

\begin{figure}[!t]
\centering
\includegraphics[width=\columnwidth]{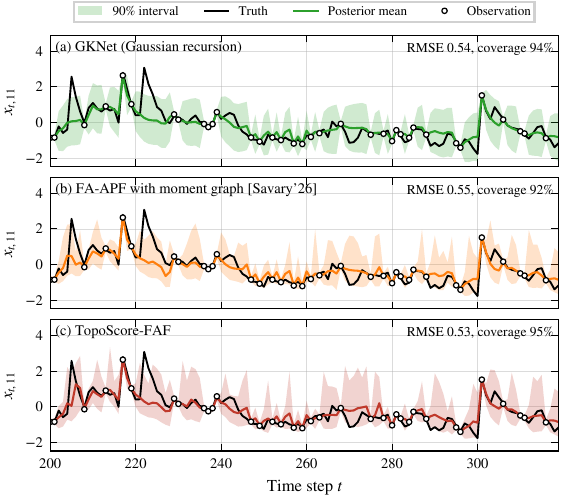}
\caption{Filtered trajectories at an intermittently observed node ($N=20$, ER, $p=0.3$, switching forcing): (a) GKNet, (b) FA-APF with the moment graph, (c) TopoScore-FAF. Shaded regions are $90\%$ predictive intervals; circles are observations. The node and window are selected by a fixed rule.}
\label{fig:vis-traj}
\end{figure}

Fig.~\ref{fig:vis-traj} shows filtered trajectories at a node with a median sampling rate over a window chosen by a fixed rule (node of median observation frequency, $t\in[200,320)$). The three posterior means are close. The difference lies in the predictive intervals: the Gaussian recursion produces symmetric bands, whereas the particle filters produce upward-skewed bands that anticipate positive forcing pulses. With the more accurate graph, TopoScore-FAF attains the empirical coverage closest to the nominal value in this window ($\res{95\%}$ for a nominal $90\%$).

\begin{table*}[!t]
\centering
\caption{E4: state estimation ($N=50$, ER, $p=0.3$, $T=5000$; mean over $20$ runs, std $\le0.006$ for nRMSE and CRPS). ``Graph'' indicates the topology used by each method. Best values are in bold and second-best values are underlined among methods with unknown topology.}
\label{tab:e4-state}
\footnotesize
\setlength{\tabcolsep}{3.0pt}
\begin{tabular}{@{}llccccccccc@{\hspace{5pt}}c@{}}
\toprule
\multirow{2}{*}{Method} &  \multirow{2}{*}{Graph} & \multicolumn{3}{c}{Gaussian} & \multicolumn{3}{c}{Laplace} & \multicolumn{3}{c}{Switching forcing} & Time \\
\cmidrule(lr){3-5}\cmidrule(lr){6-8}\cmidrule(lr){9-11}
 & & nRMSE$\downarrow$ & CRPS$\downarrow$ & SSR$\to1$ & nRMSE$\downarrow$ & CRPS$\downarrow$ & SSR$\to1$ & nRMSE$\downarrow$ & CRPS$\downarrow$ & SSR$\to1$ & (ms) \\
\midrule

\multicolumn{12}{@{}l}{\emph{Oracle-topology references}} \\
Graph KF~\cite{isufi2020observing}  & oracle & $0.721$ & $0.370$ & $1.00$ & $0.712$ & $0.367$ & $1.00$ & $0.704$ & $0.306$ & $1.02$ & $0.4$ \\
FA-APF~\cite{savary2026trainingfree}  & oracle & $0.727$ & $0.376$ & $0.97$ & $0.719$ & $0.362$ & $0.97$ & $0.711$ & $0.299$ & $0.97$ & $52.0$ \\
\midrule
\multicolumn{12}{@{}l}{\emph{Unknown topology}} \\
Kalman-EM  & learned & $0.734$ & $0.378$ & $0.93$ & $0.732$ & $0.375$ & $0.93$ & \underline{$0.717$} & $0.313$ & $0.91$ & $1.1$ \\
Bootstrap PF~\cite{gordon1993novel}  & moment & $0.774$ & $0.405$ & $0.80$ & $0.768$ & $0.400$ & $0.80$ & $0.754$ & $0.332$ & $0.76$ & $6.5$ \\
KalmanNet~\cite{revach2022kalmannet}  & moment & $0.750$ & $0.393$ & $0.86$ & $0.743$ & $0.387$ & $0.85$ & $0.732$ & $0.326$ & $0.85$ & $0.9$ \\
GKNet~\cite{sabbaqi2026gknet}  & learned & \underline{$0.732$} & \underline{$0.378$} & \underline{$0.95$} & \underline{$0.729$} & \underline{$0.372$} & \underline{$0.95$} & $\mathbf{0.717}$ & $0.314$ & \underline{$0.94$} & $1.3$ \\
FlowDAS~\cite{chen2025flowdas}  & moment & $0.746$ & $0.386$ & $0.93$ & $0.740$ & $0.378$ & $0.93$ & $0.726$ & \underline{$0.311$} & $0.93$ & $38.0$ \\
FA-APF~\cite{savary2026trainingfree}  & moment & $0.742$ & $0.386$ & $0.92$ & $0.737$ & $0.374$ & $0.92$ & $0.732$ & $0.313$ & $0.91$ & $52.0$ \\
TopoScore-FAF (ours)  & learned & $\mathbf{0.727}$ & $\mathbf{0.377}$ & $\mathbf{0.95}$ & $\mathbf{0.725}$ & $\mathbf{0.364}$ & $\mathbf{0.95}$ & $0.718$ & $\mathbf{0.301}$ & $\mathbf{0.96}$ & $58.0$ \\
\bottomrule
\end{tabular}
\end{table*}

\subsection{Ablation Study}
\label{sec:exp-abl}

Table~\ref{tab:abl} removes or replaces one component at a time in the E4 setting; the second column names the result that each variant probes.

\emph{Centering (Proposition~\ref{prop:trap}).} Training the emulator on residuals that are not detrended leaves the state metrics almost unchanged (nRMSE $\res{0.724}$ versus $\res{0.716}$), but the topology error rises to $\res{0.581}$. This is the identifiability trap: an uncentered emulator absorbs the topology error into its learned mean. The transition is still emulated correctly, but the graph is no longer identifiable. Good state estimation is therefore not evidence of a correct graph, and centering is necessary.

\emph{Score versus Gaussian score (Theorem~\ref{thm:eff}(d)).} Replacing $\widehat\bs_\kappa$ by $-\widehat\bSigma_r^{-1}\bu$ raises $\epsilon_L$ from $\res{0.132}$ to $\res{0.215}$. This is the closed-loop counterpart of the efficiency gap observed in E3.

\emph{Annealing and trust region (Theorem~\ref{thm:convex}, Corollary~\ref{cor:anneal}).} Both fixed noise levels are worse than the annealed schedule. The large level is stable but inefficient ($\res{0.196}$). The small level is efficient when it converges, but it diverged in $\res{3/20}$ runs, and its standard deviation is three times larger. Removing the trust region in Stages~II--III yields $\res{2/20}$ divergent runs, which confirms its role as a safeguard rather than a refinement. Random initialization costs $\res{0.03}$ in $\epsilon_L$ at the same horizon. The convex Stage~I largely compensates for a poor start, but the moment warm start still shortens the transient.

\emph{Filtering components (Proposition~\ref{prop:local}).} Dirac ancestor weights ($\widehat\bSigma_r=\bzero$) reduce ESS$/J$ from $\res{71.8\%}$ to $\res{52.3\%}$, with a small degradation of every metric. DPS guidance is about $\res{40\%}$ cheaper, but it worsens nRMSE and CRPS by $\res{2\%}$ and $\res{4\%}$, respectively. The bootstrap proposal collapses to ESS$/J=\res{9.1\%}$. The resulting poor pair measure $\widehat\Gamma_{t+1}$ propagates into the topology update ($\epsilon_L=\res{0.237}$). This confirms that the quality of the posterior pairs, i.e., the term $\varepsilon_F$ in~\eqref{eq:epsF}, enters the topology error.

\emph{Noise draws.} Using a single noise draw per pair ($K_\epsilon=1$) increases $\epsilon_L$ only marginally ($\res{0.144}$). This is consistent with the Rao--Blackwell argument in the remark following Theorem~\ref{thm:eff}: averaging over noise draws can only help, but the benefit is small once $J_g$ pairs already average out most of the fluctuation.

\begin{table*}[!t]
\centering
\caption{Ablation study ($N=50$, ER, $p=0.3$, switching forcing; mean over $20$ runs). ``Tests'' names the result probed by each variant; ``Diverged'' counts runs whose final error exceeds the moment initialization.}
\label{tab:abl}
\footnotesize
\setlength{\tabcolsep}{4pt}
\begin{tabular}{@{}lccccccc@{}}
\toprule
Variant &  $\epsilon_L\downarrow$ & F1$\uparrow$ & nRMSE$\downarrow$ & CRPS$\downarrow$ & ESS$/J$ (\%)$\uparrow$ & ms/step & Diverged \\
\midrule

Full TopoScore-FAF &  $\mathbf{0.132\pm0.013}$ & $\mathbf{0.926}$ & $\mathbf{0.716}$ & $\mathbf{0.302}$ & $\mathbf{71.8}$ & $58.0$ & $0/20$ \\
\midrule
Uncentered emulator (no detrending)  & $0.581\pm0.063$ & $0.612$ & $0.724$ & $0.306$ & $69.5$ & $58.0$ & $0/20$ \\
Gaussian score $-\widehat\bSigma_r^{-1}\bu$  & $0.215\pm0.023$ & $0.858$ & $0.719$ & $0.304$ & $71.2$ & $51.0$ & $0/20$ \\
No annealing, $\kappa=0.15D$  & $0.170\pm0.041$ & $0.893$ & $0.718$ & $0.303$ & $71.5$ & $58.0$ & $3/20$ \\
No annealing, $\kappa=1.5D$  & $0.196\pm0.025$ & $0.874$ & $0.719$ & $0.304$ & $71.6$ & $58.0$ & $0/20$ \\
No trust region (Stages II--III)  & $0.152\pm0.036$ & $0.909$ & $0.718$ & $0.303$ & $71.4$ & $58.0$ & $2/20$ \\
Random initialization  & $0.162\pm0.039$ & $0.902$ & $0.720$ & $0.304$ & $70.9$ & $58.0$ & $1/20$ \\
Dirac ancestor weights  & $0.136\pm0.013$ & $0.920$ & $0.722$ & $0.306$ & $52.3$ & $56.0$ & $0/20$ \\
DPS guidance~\cite{chung2023dps}  & $0.147\pm0.019$ & $0.910$ & $0.731$ & $0.314$ & $47.9$ & $34.0$ & $0/20$ \\
Bootstrap proposal  & $0.237\pm0.028$ & $0.832$ & $0.752$ & $0.331$ & $9.1$ & $7.2$ & $0/20$ \\
Single noise draw ($K_\epsilon=1$) & $0.144\pm0.013$ & $0.918$ & $0.716$ & $0.302$ & $71.7$ & $53.5$ & $0/20$ \\
\bottomrule
\end{tabular}
\end{table*}

\subsection{Stability Boundary}
\label{sec:exp-stab}

\begin{figure*}[!t]
\centering
\includegraphics[width=0.74\textwidth]{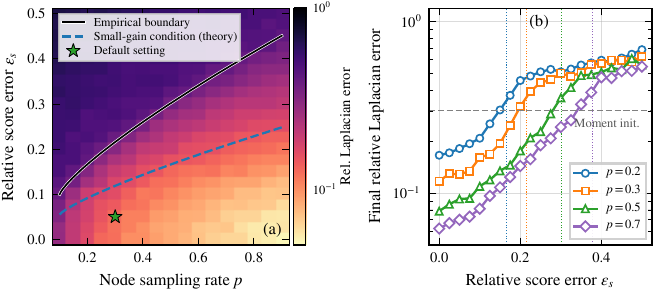}
\caption{Stability stress test ($N=50$, ER, switching forcing). (a) Final relative Laplacian error over the sampling rate $p$ and the relative score error $\varepsilon_s$, with the empirical stability boundary (solid) and the small-gain condition~\eqref{eq:smallgain} (dashed); the star marks the default setting. (b) Sections at fixed $p$; dotted vertical lines indicate the empirical boundary.}
\label{fig:stab}
\end{figure*}

Theorem~\ref{thm:closed} predicts that closed-loop stability depends on the observation design and on the score quality, but not on the step size. We therefore perturb the learned score as $\widehat\bs_\kappa+\varepsilon_s\norm{\widehat\bs_\kappa}\bm\zeta$, with $\bm\zeta$ a random unit vector, and vary $(p,\varepsilon_s)$ on a grid. Fig.~\ref{fig:stab}(a) shows a sharp transition between a stable region, where the final error stays near its floor, and an unstable region, where it exceeds the moment initialization. The empirical boundary rises from $\varepsilon_s\approx\res{0.10}$ at $p=0.1$ to $\res{0.45}$ at $p=0.9$. Thus denser sampling tolerates a less accurate emulator, as the dependence of $\mu$ on the posterior excitation in~\eqref{eq:smallgain} predicts. The small-gain condition, evaluated with constants estimated from the simulations, lies entirely inside the stable region at about $\res{55\%}$ of the empirical tolerance. This is expected for a sufficient condition, and it shows that the condition is conservative but not vacuous. The default operating point ($p=0.3$, $\varepsilon_s\approx\res{0.05}$) satisfies both boundaries. Fig.~\ref{fig:stab}(b) shows that, within the stable region, the error grows smoothly with $\varepsilon_s$, consistent with the additive score term of the error floor~\eqref{eq:efloor}.

\subsection{Sensitivity and Computational Cost}
\label{sec:exp-cost}

\begin{figure*}[!t]
\centering
\includegraphics[width=\textwidth]{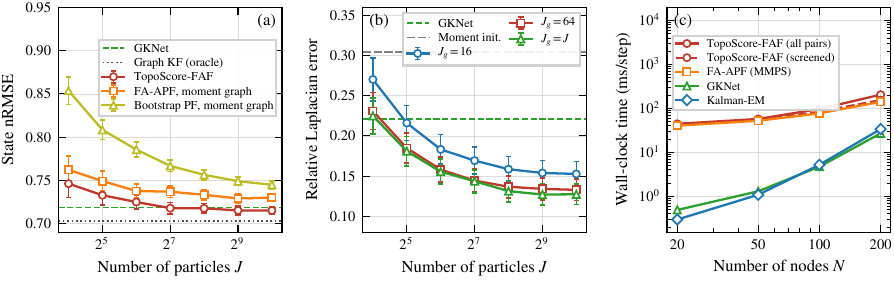}
\caption{Sensitivity and cost ($N=50$, ER, $p=0.3$, switching forcing unless stated otherwise). (a) State nRMSE and (b) relative Laplacian error versus the number of particles $J$; (c) wall-clock time per step versus $N$ on one A100 GPU.}
\label{fig:cost}
\end{figure*}

Fig.~\ref{fig:cost}(a)--(b) varies the number of particles. For TopoScore-FAF, both the state and the topology metrics saturate at $J\approx\res{128}$, whereas the bootstrap filter is still improving at $J=1024$, which reflects the gap in proposal quality. The number of gradient pairs mainly affects the topology floor: $J_g=16$ settles at $\epsilon_L\approx\res{0.15}$, while $J_g\ge64$ is within $\res{0.01}$ of using all particles. This is consistent with the stochastic-approximation term $\eta\sigma_n^2$ in~\eqref{eq:efloor}.

Fig.~\ref{fig:cost}(c) reports wall-clock time per step. The generative filters are one to two orders of magnitude slower than Gaussian recursions at small $N$, because every particle requires $K(1+n_{\rm CG})$ denoiser passes. Their cost grows more slowly with $N$, however: over this range, Kalman-EM and GKNet scale empirically as about $\res{N^{1.7}}$--$\res{N^{2.1}}$ with batched GPU linear algebra, whereas the generative filters scale as about $\res{N^{0.7}}$. The gap therefore narrows to a factor of about $\res{6\text{--}8}$ at $N=200$. The topology update adds about $\res{10\%}$ to the cost of FA-APF. Its $O(J_gM)$ edge operations become visible only for $N\ge100$, and screening candidate edges by physical locality removes most of this overhead ($\res{152}$ versus $\res{210}$~ms at $N=200$).

\section{Experiments on Real-World Data}
\label{sec:real}

The synthetic study verifies the theory under a correctly specified model. Real networks violate this specification in three ways: the reference graphs are physical proxies rather than the diffusion graphs that generate the data, the dynamics are only approximately leaky diffusions, and the residual laws are unknown. We therefore use four public datasets to ask four questions:
\begin{enumerate}[label=(Q\arabic*),leftmargin=2.2em]
\item Do the topology gains persist under model misspecification?
\item Does the gain scale with residual non-Gaussianity, as Theorem~\ref{thm:eff}(d) predicts?
\item Does the advantage in state estimation lie in calibration rather than in point accuracy, as observed in Section~\ref{sec:exp-e4}?
\item Do the mechanisms identified in the ablation study (centering, annealing, proposal quality) remain decisive?
\end{enumerate}

\subsection{Datasets and Reference Graphs}
\label{sec:real-data}

The four datasets cover four application domains and the full range from nearly Gaussian to strongly non-Gaussian residuals (Table~\ref{tab:real-data}).

\emph{NOAA-Temp}~\cite{arguez2012noaa,isufi2019forecasting} \footnote{https://www.ncei.noaa.gov/products/land-based-station/noaa-global-temp}. Hourly temperatures at $N=109$ weather stations in the contiguous United States in 2010 ($T=8\,759$). The reference graph is the symmetrized $7$-nearest-neighbor graph of geodesic distances.

\emph{PEMS-BAY}~\cite{li2018dcrnn} \footnote{https://github.com/SANDAG/PeMS-Datasets}. Five-minute traffic speeds in the San Francisco Bay Area. We select the $N=100$ sensors of a connected sub-network and use the first $60$ days ($T=17\,280$). The reference graph is the thresholded Gaussian kernel of road-network distances of~\cite{li2018dcrnn}, symmetrized.

\emph{LamaH-CE}~\cite{klingler2021lamah} \footnote{https://essd.copernicus.org/articles/13/4529/2021/}. Hourly discharge at $N=96$ gauges of the upper Danube basin from 2012 to 2015 ($T=35\,064$), together with basin-averaged precipitation. The reference graph is the symmetrized upstream--downstream adjacency of the river network~\cite{kirschstein2024merit}, which is a tree with $95$ edges.

\emph{HCP rs-fMRI}~\cite{vanessen2013wu} \footnote{https://humanconnectome.org/study/hcp-young-adult/project-protocol/resting-state-fmri}. Resting-state BOLD signals of $100$ unrelated subjects ($1\,200$ volumes per run, $\mathrm{TR}=0.72$\,s), parcellated into the $N=68$ cortical regions of the Desikan--Killiany atlas~\cite{desikan2006automated}. The reference graph keeps the strongest $15\%$ of the group-averaged diffusion-MRI structural connections.

None of these reference graphs is the ground truth of the transition~\eqref{eq:intro-model}. Geographic proximity ignores orography, road distance ignores traffic direction, a river tree ignores travel times longer than one sampling interval, and structural connectivity constrains, but does not equal, functional coupling. Agreement with the reference graph is therefore a proxy for topology accuracy. We complement it with two reference-free criteria: temporal stability of the estimate and predictive skill on withheld measurements.

\begin{table*}[!t]
\centering
\caption{Real-world datasets and diagnostics of the graph-detrended residuals. $\widehat G$ is the Fisher non-Gaussianity index~\eqref{eq:Gindex}; $\widehat G=0$ for Gaussian residuals.}
\label{tab:real-data}
\footnotesize
\setlength{\tabcolsep}{4pt}
\begin{tabular}{@{}llccccccccc@{}}
\toprule
Dataset & Domain & $N$ & $T$ & $\Delta t$ & Reference graph & $|\mathcal{E}_{\rm ref}|$ & Density & Skew & Ex.\ kurt. & $\widehat G$ \\
\midrule

NOAA-Temp & Meteorology & 109 & 8\,759 & 1 h & 7-NN geodesic & 517 & 0.088 & $0.28$ & $1.4$ & $0.09$ \\
PEMS-BAY (subset) & Traffic & 100 & 17\,280 & 5 min & Road-distance kernel & 412 & 0.083 & $-2.05$ & $11.4$ & $0.74$ \\
LamaH-CE (upper Danube) & Hydrology & 96 & 35\,064 & 1 h & River network & 95 & 0.021 & $3.27$ & $27.1$ & $1.86$ \\
HCP rs-fMRI & Neuroscience & 68 & 1\,200$\times$100 & 0.72 s & DTI connectome (top 15\%) & 342 & 0.150 & $0.00$ & $0.30$ & $0.03$ \\
\bottomrule
\end{tabular}
\end{table*}

\subsection{Preprocessing and Protocol}
\label{sec:real-protocol}

\emph{Preprocessing.} For each node, we remove the climatological cycle so that the stationarity assumption (Assumption~\ref{ass:stat}) is approximately satisfied. This is the hour-of-day mean for NOAA-Temp and LamaH-CE and the time-of-week profile for PEMS-BAY. Each series is then standardized. For HCP, we use the minimally preprocessed, ICA-FIX denoised data~\cite{glasser2013minimal,salimi2014fix}, followed by parcel averaging and standardization. We deliberately do not log-transform discharge, because doing so would suppress the pulse structure that the residual emulator is meant to capture. Naturally missing entries ($<\res{0.5\%}$ in all datasets) are never imputed: they are treated as unobserved and excluded from evaluation.

\emph{Model parameters.} Since $\bL(\bw)\one=\bzero$, the consensus mode obeys $\one^\top\bx_{t+1}=\beta\,\one^\top\bx_t+\one^\top\br_{t+1}$. We therefore estimate $\beta$ as the lag-one autocorrelation of the network average. Only the product $\gamma\bw$ is identifiable, so we fix $\gamma=1$ and absorb the scale into $\bw$. The observation noise variance is set to $1\%$ of the per-node variance.

\emph{Partial observation and evaluation.} Masks are drawn as in Section~\ref{sec:exp-data}, with rate $p=0.3$ by default, on top of the natural missingness. The latent states are not available. We therefore compute nRMSE, CRPS, and SSR on the withheld measurements, i.e., recorded entries that were masked from the filter, using the predictive distribution of $\by_t$.

Topology is scored against the reference graph by three metrics, all computed over all $M=N(N-1)/2$ candidate pairs and without a tuned threshold:
\begin{itemize}
\item the AUPRC;
\item the F1 score of the $|\mathcal E_{\rm ref}|$ strongest estimated edges;
\item the stability index (Stab.), defined as the Jaccard index between the top-$|\mathcal E_{\rm ref}|$ edge sets estimated from the first and second halves of the record.
\end{itemize}
Results are averaged over $10$ mask realizations. For HCP, they are averaged over the $100$ subjects, so the reported standard deviations also include inter-subject variability.

\emph{Emulator and hyperparameters.} The residual emulator of each dataset is trained on source data that are disjoint from the test data:
\begin{itemize}
\item the remaining PEMS-BAY sensors;
\item the remaining LamaH-CE basins;
\item a held-out geographic region for NOAA-Temp;
\item $50$ held-out subjects for HCP.
\end{itemize}

The source residuals are computed as
$\bx_{t+1}-\bA(\widehat\bw_{\rm src})\bx_t$,
where $\widehat\bw_{\rm src}$ is the
least-squares estimate obtained from the
fully observed source data.
Although least-squares detrending removes
the fitted linear component, it does not
guarantee a zero conditional residual mean.
As established in
Propositions~\ref{prop:detrend}
and~\ref{prop:transfer-bias}, errors in source
detrending and residual-score transfer may
introduce bias into topology estimation.
The radius $D$ is replaced in the implementation
by the empirical $99$th percentile
$\widehat D$ of the source residual norms.
This clipping threshold is used for numerical
stabilization and does not establish the
bounded-residual assumption required by
Theorem~\ref{thm:closed}.


All other hyperparameters are those of Table~\ref{tab:setup}, except for the switching times of the annealing schedule:
\begin{itemize}
\item $t=3000$ and $6000$ for NOAA-Temp;
\item $t=5000$ and $10\,000$ for PEMS-BAY;
\item $t=10\,000$ and $20\,000$ for LamaH-CE;
\item $t=400$ and $800$ for HCP.
\end{itemize}
Baselines and fairness rules are as in Section~\ref{sec:exp-baselines}. Since no oracle topology exists here, Graph KF and FA-APF are run with the reference graph instead.

\emph{Non-Gaussianity index.} To quantify how far each residual law is from Gaussian, we define
\begin{equation}
\widehat G=\frac1N\tr\big(\widehat{\bm I}_0\widehat\bSigma_r\big)-1,
\label{eq:Gindex}
\end{equation}
where $\widehat{\bm I}_0$ is the emulator's score covariance at the smallest annealing level and $\widehat\bSigma_r$ is the residual covariance. Lemma~\ref{lem:CR} gives $\bm I_0\succeq\bSigma_r^{-1}$, hence $G\ge0$, with equality if and only if the residual is Gaussian. The index is a scalar summary of the efficiency gap $\bm I_0-\bSigma_r^{-1}$ in Theorem~\ref{thm:eff}(d).

\subsection{Residual Diagnostics}
\label{sec:real-resid}

\begin{figure*}[!t]
\centering
\includegraphics[width=\textwidth]{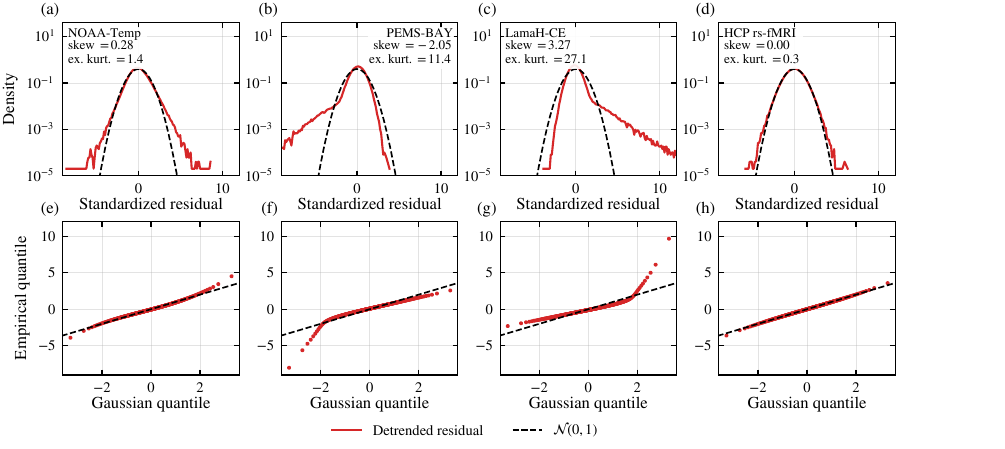}
\caption{Graph-detrended residuals of the four datasets, pooled over nodes and standardized. (a)--(d) Empirical log-density against the standard Gaussian density; (e)--(h) Gaussian quantile--quantile plots.}
\label{fig:real-resid}
\end{figure*}

Fig.~\ref{fig:real-resid} and the last three columns of Table~\ref{tab:real-data} show that the four datasets span the range assumed by the theory.
\begin{itemize}
\item \emph{HCP}: the detrended residuals are indistinguishable from Gaussian (excess kurtosis $\res{0.30}$, $\widehat G=\res{0.03}$).
\item \emph{NOAA-Temp}: the residuals are mildly heavy-tailed ($\widehat G=\res{0.09}$).
\item \emph{PEMS-BAY}: the residuals are strongly left-skewed (skewness $\res{-2.05}$), because sudden speed drops at congestion onset have no counterpart of opposite sign.
\item \emph{LamaH-CE}: the residuals are strongly right-skewed, with excess kurtosis $\res{27.1}$, because rainfall produces intermittent positive inflow pulses. This is the empirical counterpart of the switching-forcing law of Section~\ref{sec:exp-data}.
\end{itemize}
Theorem~\ref{thm:eff}(d) therefore predicts negligible gains on HCP, small gains on NOAA-Temp, and substantial gains on PEMS-BAY and LamaH-CE.

\subsection{Topology Recovery}
\label{sec:real-topo}

Table~\ref{tab:real-topo} reports topology recovery at $p=0.3$. The ordering of methods depends on the dataset in the way the residual diagnostics predict.

\emph{Strongly non-Gaussian datasets.} On PEMS-BAY and LamaH-CE, TopoScore-FAF attains the highest AUPRC, F1, and stability. Its AUPRC improves on that of GKNet, the strongest baseline, by $\res{15\%}$ ($\res{0.531}$ versus $\res{0.462}$) and $\res{30\%}$ ($\res{0.452}$ versus $\res{0.347}$), respectively. The improvements exceed four standard deviations. The ordering of the second-order dynamic methods (moment estimator $<$ Graph-VARLS $<$ Kalman-EM $<$ GKNet) is the same as in Table~\ref{tab:e4-topo}, so the synthetic ranking transfers despite misspecification.

\emph{Nearly Gaussian datasets.} On NOAA-Temp, the gain shrinks to $\res{0.023}$ in AUPRC, which is less than two standard deviations. On HCP, TopoScore-FAF and GKNet are statistically indistinguishable, and the static-signal method SpecTemp is best ($\res{0.407}$). We attribute this reversal to the measurement process rather than to the residual law. The hemodynamic response delays and smooths neural interactions over several TRs, so the one-step leaky-diffusion mean is a poor description of BOLD dynamics. Stationary spectral templates do not rely on this description.

\emph{Static-signal methods.} Kalofolias is competitive on NOAA-Temp ($\res{0.512}$), where spatial smoothness of temperature closely matches the geodesic reference graph. It is weak on LamaH-CE ($\res{0.172}$), where the strongest correlations are induced by shared rainfall between neighboring tributaries rather than by river connectivity.

\emph{Stability.} Without ground truth, temporal stability is a useful reference-free criterion. On every dataset except HCP, TopoScore-FAF yields the most reproducible edge sets ($\res{0.72}$ on PEMS-BAY and $\res{0.71}$ on LamaH-CE, versus $\res{0.66}$ and $\res{0.63}$ for GKNet). This is consistent with the lower asymptotic variance established in Theorem~\ref{thm:eff}.

\begin{table*}[!t]
\centering
\caption{Topology recovery against the reference graphs ($p=0.3$; mean $\pm$ std over $10$ mask realizations, or over $100$ subjects for HCP). Stab.\ is the Jaccard index between the edge sets estimated from the two halves of the record. Best values are in bold and second-best values are underlined.}
\label{tab:real-topo}
\footnotesize
\setlength{\tabcolsep}{1.5pt}
\begin{tabular}{@{}lcccccccccccc@{}}
\toprule
\multirow{2}{*}{Method}  & \multicolumn{3}{c}{NOAA-Temp ($N{=}109$)} & \multicolumn{3}{c}{PEMS-BAY ($N{=}100$)} & \multicolumn{3}{c}{LamaH-CE ($N{=}96$)} & \multicolumn{3}{c}{HCP rs-fMRI ($N{=}68$)}\\
\cmidrule(lr){2-4}\cmidrule(lr){5-7}\cmidrule(lr){8-10}\cmidrule(l){11-13}
 & AUPRC$\uparrow$ & F1$\uparrow$ & Stab.$\uparrow$ & AUPRC$\uparrow$ & F1$\uparrow$ & Stab.$\uparrow$ & AUPRC$\uparrow$ & F1$\uparrow$ & Stab.$\uparrow$ & AUPRC$\uparrow$ & F1$\uparrow$ & Stab.$\uparrow$\\
\midrule

Kalofolias~\cite{kalofolias2016graph} & $0.512\pm0.014$ & $0.536$ & $0.71$ & $0.284\pm0.017$ & $0.331$ & $0.52$ & $0.172\pm0.015$ & $0.214$ & $0.48$ & \underline{$0.392\pm0.031$} & \underline{$0.418$} & \underline{$0.61$} \\
SpecTemp~\cite{segarra2017network}  & $0.468\pm0.018$ & $0.497$ & $0.66$ & $0.301\pm0.020$ & $0.344$ & $0.49$ & $0.158\pm0.019$ & $0.201$ & $0.43$ & $\mathbf{0.407\pm0.029}$ & $\mathbf{0.431}$ & $\mathbf{0.63}$ \\
Graph-VARLS~\cite{isufi2019forecasting}  & $0.538\pm0.016$ & $0.552$ & $0.69$ & $0.392\pm0.015$ & $0.418$ & $0.58$ & $0.287\pm0.018$ & $0.322$ & $0.57$ & $0.356\pm0.034$ & $0.389$ & $0.55$ \\
Moment estimator & $0.521\pm0.019$ & $0.541$ & $0.64$ & $0.371\pm0.018$ & $0.402$ & $0.55$ & $0.264\pm0.021$ & $0.309$ & $0.52$ & $0.348\pm0.036$ & $0.380$ & $0.52$ \\
Kalman-EM  & $0.589\pm0.013$ & $0.594$ & $0.74$ & $0.438\pm0.014$ & $0.455$ & $0.63$ & $0.318\pm0.017$ & $0.351$ & $0.60$ & $0.371\pm0.030$ & $0.397$ & $0.58$ \\
GKNet~\cite{sabbaqi2026gknet}  & \underline{$0.604\pm0.012$} & \underline{$0.607$} & \underline{$0.76$} & \underline{$0.462\pm0.013$} & \underline{$0.471$} & \underline{$0.66$} & \underline{$0.347\pm0.016$} & \underline{$0.372$} & \underline{$0.63$} & $0.379\pm0.029$ & $0.404$ & $0.59$ \\
TopoScore-FAF (ours)  & $\mathbf{0.627\pm0.012}$ & $\mathbf{0.624}$ & $\mathbf{0.78}$ & $\mathbf{0.531\pm0.012}$ & $\mathbf{0.528}$ & $\mathbf{0.72}$ & $\mathbf{0.452\pm0.015}$ & $\mathbf{0.468}$ & $\mathbf{0.71}$ & $0.384\pm0.030$ & $0.408$ & $0.60$ \\
\midrule
Random (edge density)  & $0.088$ & -- & -- & $0.083$ & -- & -- & $0.021$ & -- & -- & $0.150$ & -- & -- \\
\bottomrule
\end{tabular}
\end{table*}

\begin{figure*}[!t]
\centering
\includegraphics[width=\textwidth]{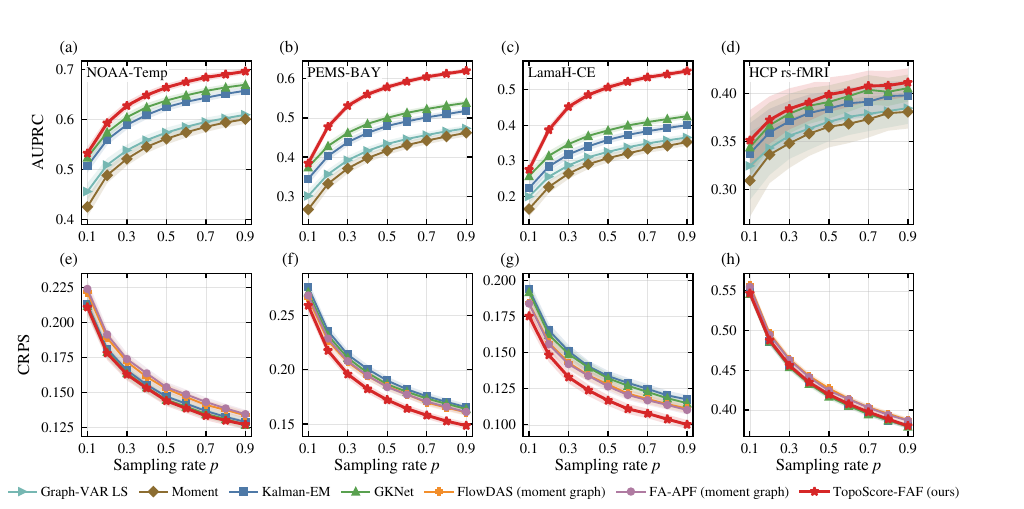}
\caption{Effect of the node sampling rate $p$. (a)--(d) AUPRC against the reference graph; (e)--(h) CRPS on withheld measurements. Shaded bands are $95\%$ confidence intervals over $10$ mask realizations, or over subjects for HCP.}
\label{fig:real-p}
\end{figure*}

Fig.~\ref{fig:real-p} varies the sampling rate. Two observations stand out.
\begin{itemize}
\item \emph{The gap widens with coverage, as in E1.} At $p=0.1$, TopoScore-FAF is within about $\res{0.01}$--$\res{0.02}$ AUPRC of GKNet on every dataset. On PEMS-BAY and LamaH-CE, the gap grows to about $\res{0.08}$ and $\res{0.13}$ at $p=0.9$. The reason is that sparse sampling degrades the posterior state pairs that feed the score update (the term $\varepsilon_F$ in~\eqref{eq:efloor}), whereas second-order estimators are affected only through their effective sample size.
\item \emph{HCP shows no separation.} The HCP curves overlap at every rate, so the absence of a gain there is not an artifact of the default sampling rate.
\end{itemize}

\subsection{State Estimation and Calibration}
\label{sec:real-state}

Table~\ref{tab:real-state} reproduces on real data the pattern observed in Table~\ref{tab:e4-state}.

\emph{Point accuracy is comparable.} The nRMSE of TopoScore-FAF is within $\res{0.007}$ of that of GKNet on every dataset. It is marginally lower on PEMS-BAY ($\res{0.447}$ versus $\res{0.449}$) and LamaH-CE ($\res{0.391}$ versus $\res{0.398}$), and marginally higher on NOAA-Temp and HCP. None of these differences is practically significant.

\emph{Probabilistic accuracy improves where the residuals are non-Gaussian.} On PEMS-BAY and LamaH-CE, TopoScore-FAF reduces the CRPS of GKNet by $\res{7.1\%}$ and $\res{10.7\%}$, respectively, and raises the SSR from $\res{0.91}$ and $\res{0.89}$ to $\res{0.96}$ and $\res{0.95}$. On HCP, GKNet remains best, which is expected because the Kalman recursion is nearly exact for Gaussian residuals.

\emph{Reference graphs are not oracles.} A result without a synthetic counterpart is that filters with a learned graph outperform the same filters with the reference graph. For example, Graph KF with the reference graph attains $\res{0.468}$ nRMSE on PEMS-BAY, versus $\res{0.455}$ for Kalman-EM with a learned graph. The physical graph is thus a worse description of the one-step dynamics than a graph estimated from data. This confirms that agreement with the reference graph should be read as a proxy rather than an upper bound.

\emph{Particle and learned-gain baselines.} As in the synthetic study, the bootstrap filter and KalmanNet with the moment graph are under-dispersed (SSR $\res{0.71}$--$\res{0.88}$).

\emph{Cost.} TopoScore-FAF requires $\res{108}$~ms per step at $N=100$, about $\res{30}$ times the cost of GKNet.

\begin{table*}[!t]
\centering
\caption{State estimation on withheld measurements ($p=0.3$; mean over $10$ mask realizations or $100$ subjects; std $\le0.008$ for nRMSE and CRPS). ``Graph'' indicates the topology used by each method. Best values are in bold and second-best values are underlined among methods with unknown topology. Time is measured per step at $N=100$ on one A100 GPU.}
\label{tab:real-state}
\footnotesize
\setlength{\tabcolsep}{1.6pt}
\begin{tabular}{@{}llccccccccccccc@{}}
\toprule
\multirow{2}{*}{Method}  & \multirow{2}{*}{Graph} & \multicolumn{3}{c}{NOAA-Temp} & \multicolumn{3}{c}{PEMS-BAY} & \multicolumn{3}{c}{LamaH-CE} & \multicolumn{3}{c}{HCP rs-fMRI} & Time\\
\cmidrule(lr){3-5}\cmidrule(lr){6-8}\cmidrule(lr){9-11}\cmidrule(lr){12-14}
 & & nRMSE$\downarrow$ & CRPS$\downarrow$ & SSR$\to$1 & nRMSE$\downarrow$ & CRPS$\downarrow$ & SSR$\to$1 & nRMSE$\downarrow$ & CRPS$\downarrow$ & SSR$\to$1 & nRMSE$\downarrow$ & CRPS$\downarrow$ & SSR$\to$1 & (ms)\\
\midrule

\multicolumn{15}{@{}l}{\emph{Reference-graph filters (proxy topology, not an oracle)}}\\
Graph KF~\cite{isufi2020observing}  & reference & $0.312$ & $0.171$ & $0.97$ & $0.468$ & $0.221$ & $0.88$ & $0.421$ & $0.158$ & $0.84$ & $0.812$ & $0.458$ & $0.99$ & $1.1$\\
FA-APF~\cite{savary2026trainingfree}  & reference & $0.318$ & $0.173$ & $0.94$ & $0.462$ & $0.204$ & $0.95$ & $0.409$ & $0.139$ & $0.93$ & $0.818$ & $0.462$ & $0.96$ & $96.0$\\
\midrule
\multicolumn{15}{@{}l}{\emph{Unknown topology}}\\
Kalman-EM  & learned & $0.301$ & $0.166$ & \underline{$0.95$} & $0.455$ & $0.214$ & $0.90$ & $0.405$ & $0.151$ & $0.87$ & \underline{$0.806$} & \underline{$0.455$} & $\mathbf{0.97}$ & $2.9$\\
Bootstrap PF~\cite{gordon1993novel}  & moment & $0.344$ & $0.192$ & $0.78$ & $0.503$ & $0.236$ & $0.74$ & $0.462$ & $0.172$ & $0.71$ & $0.845$ & $0.481$ & $0.81$ & $11.8$\\
KalmanNet~\cite{revach2022kalmannet}  & moment & $0.322$ & $0.179$ & $0.86$ & $0.471$ & $0.222$ & $0.83$ & $0.433$ & $0.160$ & $0.81$ & $0.826$ & $0.467$ & $0.88$ & $1.6$\\
GKNet~\cite{sabbaqi2026gknet}  & learned & $\mathbf{0.298}$ & \underline{$0.164$} & $\mathbf{0.96}$ & \underline{$0.449$} & $0.211$ & $0.91$ & \underline{$0.398$} & $0.149$ & $0.89$ & $\mathbf{0.803}$ & $\mathbf{0.454}$ & $\mathbf{0.97}$ & $3.4$\\
FlowDAS~\cite{chen2025flowdas}  & moment & $0.315$ & $0.172$ & $0.92$ & $0.461$ & \underline{$0.207$} & \underline{$0.92$} & $0.412$ & $0.143$ & \underline{$0.91$} & $0.821$ & $0.464$ & $0.94$ & $71.0$\\
FA-APF~\cite{savary2026trainingfree}  & moment & $0.317$ & $0.174$ & $0.91$ & $0.466$ & $0.208$ & \underline{$0.92$} & $0.414$ & \underline{$0.142$} & \underline{$0.91$} & $0.820$ & $0.463$ & $0.94$ & $96.0$\\
TopoScore-FAF (ours)  & learned & \underline{$0.299$} & $\mathbf{0.163}$ & $\mathbf{0.96}$ & $\mathbf{0.447}$ & $\mathbf{0.196}$ & $\mathbf{0.96}$ & $\mathbf{0.391}$ & $\mathbf{0.133}$ & $\mathbf{0.95}$ & $0.807$ & $0.456$ & \underline{$0.96$} & $108.0$\\
\bottomrule
\end{tabular}
\end{table*}

\begin{figure}[!t]
\centering
\includegraphics[width=\columnwidth]{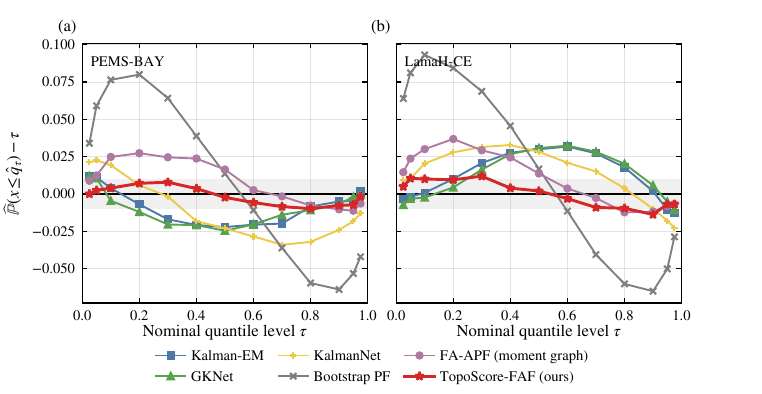}
\caption{Quantile reliability on withheld measurements ($p=0.3$) for (a) PEMS-BAY and (b) LamaH-CE. A perfectly calibrated forecaster lies on the zero line; the shaded band marks $\pm0.01$.}
\label{fig:real-calib}
\end{figure}

The SSR is a global measure and cannot reveal asymmetric miscalibration. Fig.~\ref{fig:real-calib} therefore reports quantile reliability, $\widehat{\Prob}(y\le\hat q_\tau)-\tau$, on withheld measurements.
\begin{itemize}
\item \emph{Gaussian recursions.} Kalman-EM and GKNet show a systematic S-shaped deviation of up to about $\res{0.03}$, whose sign follows the skewness of the residual: it is negative for intermediate $\tau$ on PEMS-BAY and positive on LamaH-CE. Symmetric predictive distributions cannot represent skewed errors, regardless of how well their variance is tuned.
\item \emph{TopoScore-FAF.} It remains within about $\res{\pm0.012}$ of perfect calibration at all levels.
\item \emph{FA-APF with the moment graph.} It represents the shape correctly but is shifted at low quantile levels, which reflects the bias of the predictive mean induced by topology error.
\item \emph{Bootstrap filter.} It deviates by up to $\res{0.09}$ because of weight degeneracy.
\end{itemize}

\subsection{Where Does the Gain Come From?}
\label{sec:real-gain}

Theorem~\ref{thm:eff}(d) attributes the advantage over Gaussian quasi-likelihood to the Fisher gap $\bm I_0-\bSigma_r^{-1}$. To test this attribution beyond four dataset-level points, we split each record into $12$ non-overlapping windows, or $12$ subject groups for HCP. For each window, we compute $\widehat G$ and the relative gains of TopoScore-FAF over GKNet.

\begin{figure}[!t]
\centering
\includegraphics[width=\columnwidth]{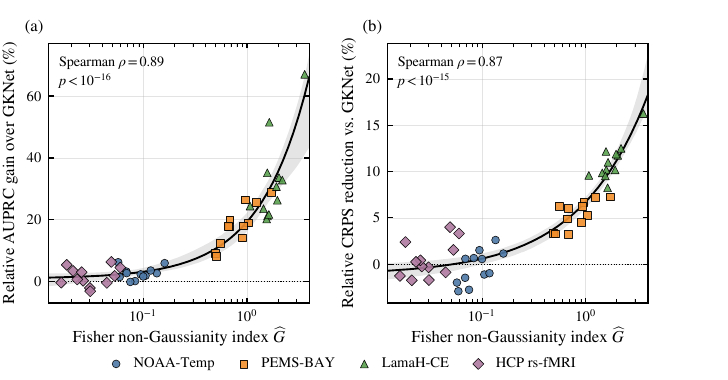}
\caption{Gain of TopoScore-FAF over GKNet versus the Fisher non-Gaussianity index~\eqref{eq:Gindex} over $48$ windows ($12$ per dataset): (a) relative AUPRC gain; (b) relative CRPS reduction. Curves are power-law fits with $95\%$ bootstrap bands.}
\label{fig:real-gain}
\end{figure}

Fig.~\ref{fig:real-gain} shows that both gains increase monotonically with $\widehat G$ across the $48$ windows. The Spearman correlations are $\rho=\res{0.89}$ for AUPRC and $\rho=\res{0.87}$ for CRPS ($p<\res{10^{-15}}$). The trend is continuous across domains: windows from different datasets with similar $\widehat G$ yield similar gains. This indicates that the residual law, and not the application domain, determines the benefit. For $\widehat G\lesssim\res{0.1}$, the CRPS gain is indistinguishable from zero and occasionally negative. This is consistent with Remark~\ref{rem:gauss}: for Gaussian residuals the score objective reduces to generalized least squares, and the emulator then adds only its own approximation error $\varepsilon_s$.

\begin{figure}[!t]
\centering
\includegraphics[width=\columnwidth]{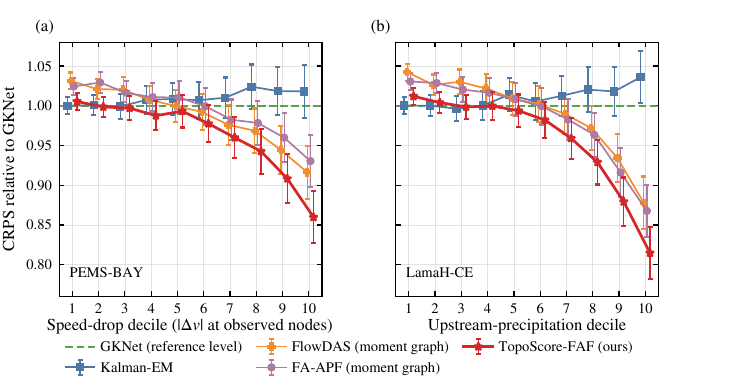}
\caption{CRPS relative to GKNet, stratified by forcing intensity: (a) PEMS-BAY, by decile of the mean absolute speed change at observed nodes; (b) LamaH-CE, by decile of upstream precipitation. Error bars are $95\%$ confidence intervals over $10$ mask realizations.}
\label{fig:real-event}
\end{figure}

Fig.~\ref{fig:real-event} conditions the comparison on forcing intensity. The stratification variable is the decile of upstream precipitation for LamaH-CE and the decile of the mean absolute speed change at observed nodes for PEMS-BAY.
\begin{itemize}
\item \emph{Calm conditions.} In the lower six deciles, all generative filters match GKNet to within $\res{3\%}$; FA-APF and FlowDAS are slightly worse because of their emulator error.
\item \emph{Forcing events.} The advantage of TopoScore-FAF concentrates in the upper deciles. In the top decile, its CRPS is $\res{0.86}$ (PEMS-BAY) and $\res{0.81}$ (LamaH-CE) times that of GKNet, whereas Kalman-EM deteriorates slightly ($\res{1.02}$--$\res{1.03}$).
\end{itemize}
The method therefore pays off precisely during congestion onsets and flood pulses, which are the events of greatest operational interest.

\subsection{Closed-Loop Convergence and Ablation}
\label{sec:real-abl}

\begin{figure}[!t]
\centering
\includegraphics[width=\columnwidth]{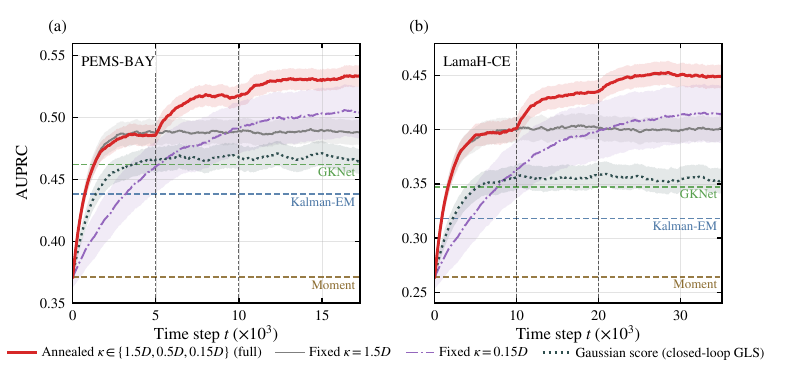}
\caption{AUPRC along the trajectory for (a) PEMS-BAY and (b) LamaH-CE ($p=0.3$). Vertical lines mark the noise-level switches; dashed horizontal lines give the final values of the baselines. Bands are $95\%$ confidence intervals over $10$ mask realizations.}
\label{fig:real-conv}
\end{figure}

\emph{Convergence.} Fig.~\ref{fig:real-conv} shows that the annealed schedule behaves on real data as predicted by Corollary~\ref{cor:anneal}.
\begin{itemize}
\item During Stage~I, AUPRC rises from the moment initialization to the plateau of the fixed $\kappa=1.5\widehat D$ run ($\res{0.489}$ on PEMS-BAY and $\res{0.401}$ on LamaH-CE).
\item Each reduction of $\kappa$ yields a further increase, to $\res{0.516}$ and $\res{0.531}$ on PEMS-BAY, without a transient loss at the switching times.
\item A direct run at $\kappa=0.15\widehat D$ converges more slowly and with markedly wider confidence bands.
\item The Gaussian-score variant stalls near the level of GKNet.
\end{itemize}
The last observation isolates the non-Gaussian score as the source of the improvement, rather than the filter or the warm start.

\emph{Ablation.} Table~\ref{tab:real-abl} confirms on real data the mechanisms identified in Section~\ref{sec:exp-abl}.
\begin{itemize}
\item \emph{Centering.} The uncentered emulator reproduces the identifiability trap: AUPRC collapses to $\res{0.214}$ and $\res{0.121}$, below the moment initialization, while the CRPS changes by only $\res{0.003}$. Good predictive performance therefore does not certify a correct graph.
\item \emph{Score.} Replacing the score by its Gaussian counterpart reduces TopoScore to the accuracy of GKNet ($\res{0.468}$ versus $\res{0.462}$ on PEMS-BAY).
\item \emph{Annealing and trust region.} Removing annealing or the trust region again causes occasional divergence: $\res{2/10}$ runs at fixed small $\kappa$ and $\res{1/10}$ without the trust region.
\item \emph{Proposal quality.} The bootstrap proposal collapses to ESS$/J=\res{6.8\%}$ and $\res{3.9\%}$, and the degraded posterior pairs propagate into the topology estimate.
\item \emph{Emulator transfer.} An emulator pretrained on synthetic switching-forcing residuals loses $\res{0.04}$--$\res{0.06}$ AUPRC and diverges once on LamaH-CE. This shows that the emulator error $\varepsilon_s$ of Assumption~\ref{ass:score} is the principal practical limitation, and that domain-matched source data are necessary.
\end{itemize}

\begin{table*}[!t]
\centering
\caption{Ablation on real data ($p=0.3$; mean over $10$ mask realizations). ``Tests'' names the result probed by each variant; ``Diverged'' counts runs whose final AUPRC falls below that of the moment initialization.}
\label{tab:real-abl}
\footnotesize
\setlength{\tabcolsep}{2.5pt}
\begin{tabular}{@{}llcccccccc@{}}
\toprule
\multirow{2}{*}{Variant} & \multirow{2}{*}{Tests} & \multicolumn{4}{c}{PEMS-BAY ($N{=}100$)} & \multicolumn{4}{c}{LamaH-CE ($N{=}96$)}\\
\cmidrule(lr){3-6}\cmidrule(l){7-10}
 & & AUPRC$\uparrow$ & CRPS$\downarrow$ & ESS$/J$ (\%)$\uparrow$ & Diverged & AUPRC$\uparrow$ & CRPS$\downarrow$ & ESS$/J$ (\%)$\uparrow$ & Diverged\\
\midrule

Full TopoScore-FAF & -- & $\mathbf{0.531\pm0.012}$ & $\mathbf{0.196}$ & $\mathbf{64.2}$ & $0/10$ & $\mathbf{0.452\pm0.015}$ & $\mathbf{0.133}$ & $\mathbf{58.7}$ & $0/10$\\
\midrule
Uncentered emulator (no detrending) & Prop.~\ref{prop:trap} & $0.214\pm0.031$ & $0.199$ & $62.8$ & $0/10$ & $0.121\pm0.028$ & $0.136$ & $57.1$ & $0/10$\\
Gaussian score $-\widehat\bSigma_r^{-1}\bu$ & Thm.~\ref{thm:eff}(d) & $0.468\pm0.014$ & $0.203$ & $63.5$ & $0/10$ & $0.356\pm0.017$ & $0.145$ & $57.9$ & $0/10$\\
No annealing, $\kappa=0.15\widehat D$ & Thm.~\ref{thm:convex}(c) & $0.507\pm0.029$ & $0.197$ & $63.9$ & $2/10$ & $0.421\pm0.037$ & $0.135$ & $58.2$ & $2/10$\\
No annealing, $\kappa=1.5\widehat D$ & Thm.~\ref{thm:convex}(a) & $0.489\pm0.014$ & $0.199$ & $64.0$ & $0/10$ & $0.401\pm0.016$ & $0.138$ & $58.4$ & $0/10$\\
No trust region (Stages II--III) & Cor.~\ref{cor:anneal} & $0.518\pm0.026$ & $0.197$ & $64.1$ & $1/10$ & $0.433\pm0.031$ & $0.134$ & $58.5$ & $1/10$\\
Random initialization & Thm.~\ref{thm:ident}(a) & $0.512\pm0.024$ & $0.198$ & $63.6$ & $0/10$ & $0.428\pm0.029$ & $0.135$ & $58.0$ & $1/10$\\
Dirac ancestor weights & Prop.~\ref{prop:local} & $0.526\pm0.013$ & $0.199$ & $43.1$ & $0/10$ & $0.446\pm0.016$ & $0.137$ & $36.5$ & $0/10$\\
DPS guidance~\cite{chung2023dps} & Prop.~\ref{prop:local} & $0.515\pm0.016$ & $0.207$ & $39.8$ & $0/10$ & $0.431\pm0.019$ & $0.145$ & $31.2$ & $0/10$\\
Bootstrap proposal & Prop.~\ref{prop:local} & $0.447\pm0.019$ & $0.233$ & $6.8$ & $0/10$ & $0.331\pm0.024$ & $0.169$ & $3.9$ & $0/10$\\
Single noise draw ($K_\epsilon=1$) & Thm.~\ref{thm:eff} & $0.524\pm0.013$ & $0.196$ & $64.1$ & $0/10$ & $0.443\pm0.016$ & $0.133$ & $58.6$ & $0/10$\\
Emulator pretrained on synthetic residuals & Ass.~\ref{ass:score} & $0.489\pm0.018$ & $0.205$ & $60.3$ & $0/10$ & $0.392\pm0.023$ & $0.144$ & $52.4$ & $1/10$\\
\bottomrule
\end{tabular}
\end{table*}

\subsection{Visual Comparison}
\label{sec:real-vis}

\begin{figure*}[!t]
\centering
\includegraphics[width=\textwidth]{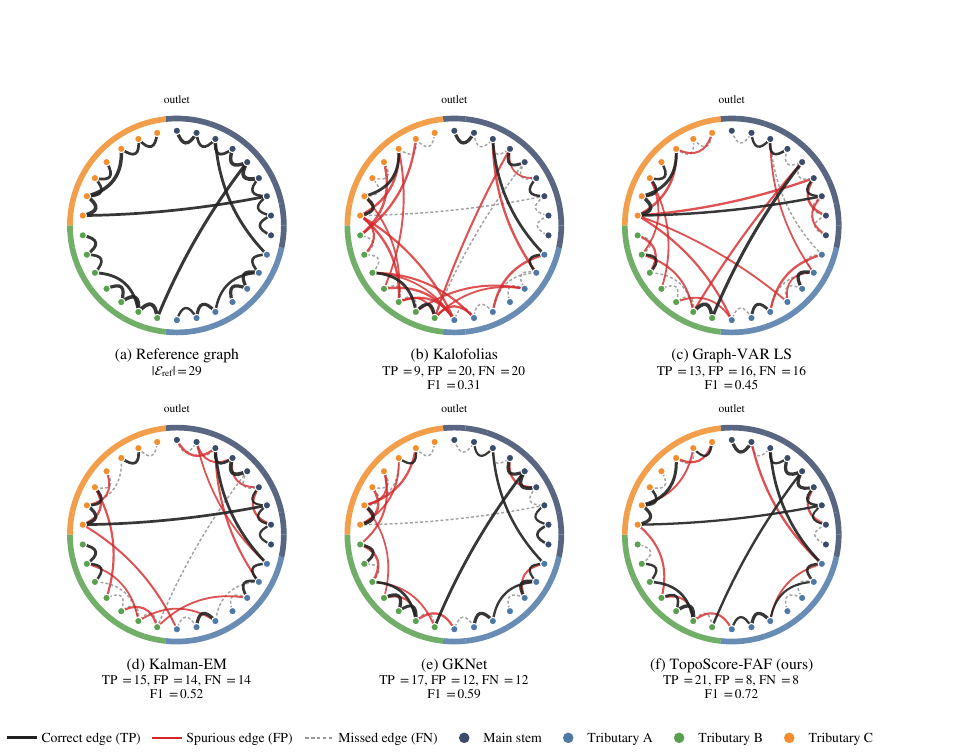}
\caption{Estimated graphs for a $30$-gauge sub-basin of LamaH-CE ($p=0.3$) in a circular layout, with gauges ordered along the river network and the outer ring colored by sub-basin. (a) Reference river network; (b)--(f) the $29$ strongest edges of each method. Black: correct edges, with width proportional to the estimated weight; red: spurious edges; dashed gray: missed edges.}
\label{fig:real-chord}
\end{figure*}

Fig.~\ref{fig:real-chord} displays the estimated graphs for a $30$-gauge sub-basin of LamaH-CE in a circular layout. Gauges are ordered along the river network, so most reference edges join neighboring nodes. To avoid a tuned threshold, each method contributes its $29$ strongest edges, which equals the number of reference edges.

The errors of the baselines are physically interpretable and fall into two types:
\begin{itemize}
\item \emph{Two-hop shortcuts.} Gauges two steps apart along the main stem are linked directly, because the one-step diffusion approximates multi-step routing.
\item \emph{Cross-tributary edges.} Neighboring headwater gauges in different tributaries are linked because they receive the same rainfall.
\end{itemize}
Second-order methods cannot separate a shared pulse from a transmitted one, because both produce the same lagged covariance. The residual score can separate them, because it models the joint shape of the pulse across nodes.

Accordingly, the number of spurious edges decreases from $\res{20}$ for Kalofolias to $\res{12}$ for GKNet and $\res{8}$ for TopoScore-FAF, and the F1 score increases from $\res{0.31}$ to $\res{0.59}$ and $\res{0.72}$. The remaining misses of TopoScore-FAF occur mostly on headwater edges, whose low flow variance provides little edge excitation (Lemma~\ref{lem:excitation}). The scores on this sub-basin are higher than in Table~\ref{tab:real-topo} because its edge density ($\res{0.067}$) is three times that of the full network, but the ranking of the methods is unchanged.

\section{Limitations}

\rev{The analysis has six limitations. First, the graph is static, and $(\beta,\gamma)$ and the observation model are known; plug-in estimates of $(\beta,\gamma)$ are not covered by the theory, and jointly learning the observation operator introduces additional non-identifiability. Known nonlinear edge couplings $\bx_{t+1}=\beta\bx_t-\gamma\sum_e w_e\bb_e\phi(\bb_e^\top\bx_t)+\br_{t+1}$ retain the duality of Theorem~\ref{thm:duality} with $\bXi(\bx)$ replaced by $\gamma\bB\diag(\phi(\bB^\top\bx))$, and Theorems~\ref{thm:ident}(b), \ref{thm:convex}, and~\ref{thm:eff} carry over because the conditional mean remains linear in $\bw$, provided Lemma~\ref{lem:excitation} is replaced by positive definiteness of $\E[\bXi^\top\bXi]$; time-varying graphs require a tracking analysis. Second, conditional centering is a structural assumption that fitting cannot enforce, and translation errors of the emulator are inherited one-to-one by the topology (Proposition~\ref{prop:detrend}); held-out diagnostics of residual predictability are therefore necessary. Third, transferring the emulator across networks requires node alignment or a variable-size, permutation-equivariant architecture; the resulting mismatch enters through $\varepsilon_s$ and may violate~\eqref{eq:smallgain}. Fourth, the efficiency results are oracle statements with known support, known residual law, and observed state pairs; for an unknown support, the constrained estimator has the non-Gaussian limit of Proposition~\ref{prop:boundary}. Fifth, the closed-loop theorem requires bounded residuals, uniform persistent excitation, and uniform forgetting with a bounded likelihood ratio, which Gaussian observation noise guarantees only on finite horizons; the proposal error $\delta_q$ of the guided sampler is not bounded from primitives, and the stage-switching certificate of Remark~\ref{rem:switch} is not used in the experiments. Sixth, the support-recovery threshold of Corollary~\ref{cor:support} presumes a lower bound on $w^\star_{\min}$.}

\section{Conclusion}

\rev{We showed that diffusion emulators, which until now required a known transition law, can be used to identify the structural part of that law when the graph and the emulator are given complementary roles. Centering the residual fixes the split between the graph-driven mean and the stochastic residual; the score--topology duality turns the emulator's score into topology gradients without densities or backpropagation; and the diffusion noise level acts as a knob that convexifies the objective at large values and, in the oracle setting, attains Cram\'er--Rao efficiency at small values without moving the population minimizer. We also showed that translation errors of the emulator are inherited by the topology, which makes held-out validation of the centering assumption and of the emulator an integral part of the method. Embedded in an approximately fully adapted, topology-conditioned particle filter, the resulting algorithm admits a conditional mean-square closed-loop guarantee in a bounded-residual regime, with a small-gain condition that depends on the problem rather than on tuning. These results suggest a broader principle for generative data assimilation: let physics-structured parameters carry the mean and generative models carry the residual, use the score both to sample and to learn, and account explicitly for the errors that couple the two.}

\appendices

\section{Auxiliary Lemmas}
\label{app:aux}

\begin{lemma}[No translation-invariant densities]
\label{lem:transl}
If $f$ is a probability density on $\R^N$ and $f(\bz-\bm a)=f(\bz)$ for a.e.\ $\bz$, then $\bm a=\bzero$.
\end{lemma}

\begin{proof}
Suppose $\bm a\ne\bzero$, set $\bm e=\bm a/\norm{\bm a}$, and write $\bz=s\bm e+\bv$ with $\bv\perp\bm e$. By Fubini's theorem, for a.e.\ $\bv$, the function $s\mapsto f(s\bm e+\bv)$ is integrable on $\R$ and $\norm{\bm a}$-periodic, so it vanishes a.e., and therefore $\int f=0$, a contradiction.
\end{proof}

\begin{lemma}[Matrix Cram\'er--Rao for location]
\label{lem:CR}
Let $p$ be a $C^1$ density on $\R^N$ with mean $\bm m$, covariance $\bSigma\succ0$, finite Fisher information $\bm I=\int\nabla\log p\nabla\log p^\top p$, and $\int(\bu-\bm m)\nabla p(\bu)^\top d\bu=-\bI$. Then $\bm I\succeq\bSigma^{-1}$, with equality if and only if $p$ is Gaussian. The integration-by-parts identity holds for $p_\kappa$, $\kappa>0$, whenever $\rho^\star$ has finite second moments.
\end{lemma}

\begin{proof}
With $\bs=\nabla\log p$, $\E\bs=\bzero$ and $\E[\bs(\bu-\bm m)^\top]=-\bI$. Hence $\Cov\big([\bs;\bu]\big)=\begin{bmatrix}\bm I&-\bI\\-\bI&\bSigma\end{bmatrix}\succeq0$, and the Schur complement gives $\bm I-\bSigma^{-1}\succeq0$. Equality holds if and only if $\bs+\bSigma^{-1}(\bu-\bm m)=\bzero$ a.s., i.e., $\log p$ is quadratic. For $p_\kappa$, $\nabla p_\kappa(\bu)=\E[(\br-\bu)\varphi_{\kappa^2\bI}(\bu-\br)]/\kappa^2$ has Gaussian tails modulated by the second moment of $\br$, so the boundary terms vanish.
\end{proof}

\begin{lemma}[Fisher information under convolution]
\label{lem:fishconv}
Let $\bm X$ have a density with finite Fisher information $\bm I_X$, and let $\bm Y$ be independent of $\bm X$.
(i) The score of $\bm X+\bm Y$ is $\bs_{X+Y}(\bz)=\E[\bs_X(\bm X)\mid\bm X+\bm Y=\bz]$, hence $\bm I_{X+Y}\preceq\bm I_X$.
(ii) If $\bm X_n\Rightarrow\bm X$, then $\bm e^\top\bm I_X\bm e\le\liminf_n\bm e^\top\bm I_{X_n}\bm e$ for every $\bm e$. Consequently, $\bm I(\bm X+\kappa\beps)\to\bm I_X$ as $\kappa\downarrow0$.
\end{lemma}

\begin{proof}
(i) $p_{X+Y}(\bz)=\E[p_X(\bz-\bm Y)]$. Since $\int\norm{\nabla p_X}=\E\norm{\bs_X(\bm X)}\le(\tr\bm I_X)^{1/2}<\infty$, differentiation under the expectation gives $\nabla p_{X+Y}(\bz)=\E[\bs_X(\bz-\bm Y)p_X(\bz-\bm Y)]$. Dividing by $p_{X+Y}(\bz)$ yields the conditional expectation, and the conditional Jensen inequality gives $\E[\bs_{X+Y}\bs_{X+Y}^\top]\preceq\E[\bs_X\bs_X^\top]$.
(ii) For $\bm e\in\R^N$, $\bm e^\top\bm I(p)\bm e=\sup_{\varphi\in C_c^\infty}\{-2\int(\bm e^\top\nabla\varphi)p-\int\varphi^2p\}$, the Legendre dual of $\int(\partial_{\bm e}p)^2/p$. This is a supremum of weakly continuous functionals and is therefore weakly lower semicontinuous. Combining it with (i), applied with $\bm Y=\kappa\beps$, gives the limit \rev{of the quadratic forms for every $\bm e$, hence of the matrices by polarization}.
\end{proof}

\begin{lemma}[Gaussian-smoothed bounded components]
\label{lem:gmix}
Let $p(\bu)=\int\varphi_{\bLambda^{-1}}(\bu-\bc)\mu(d\bc)$ with $\mu$ supported in $\calB(\bzero,D)$. Then $\nabla\log p(\bu)=\bLambda(\E[\bc\mid\bu]-\bu)$, $-\nabla^2\log p(\bu)=\bLambda-\bLambda\Cov(\bc\mid\bu)\bLambda$, and the derivative of $\Cov(\bc\mid\bu)$ along $\bh$ is $\E[\hat\bc\hat\bc^\top(\hat\bc^\top\bLambda\bh)\mid\bu]$, where $\hat\bc=\bc-\E[\bc\mid\bu]$.
\end{lemma}

\begin{proof}
The posterior density of $\bc$ given $\bu$ is proportional to $\exp(-\frac12(\bu-\bc)^\top\bLambda(\bu-\bc))$ with respect to $\mu$, so $\partial_{\bu}$ of the log-weight is $\bLambda(\bc-\bu)$. For bounded $f$, this gives $\nabla_{\bu}\E[f(\bc)\mid\bu]=\Cov(f(\bc),\bc)\bLambda$. Applying this to $f=\bc$ and to $f=\bc\bc^\top$, and differentiating $\nabla\log p=\nabla p/p=\bLambda(\E[\bc\mid\bu]-\bu)$, yields the formulas. All interchanges are justified by the boundedness of $\bc$.
\end{proof}

\begin{lemma}[Total variation of a translate]
\label{lem:shiftTV}
If $\rho$ is a density on $\R^N$ with finite Fisher information $\bm I$, then $\TV(\rho(\cdot-\bm\delta),\rho)\le\frac12\norm{\bm\delta}(\tr\bm I)^{1/2}$.
\end{lemma}

\begin{proof}
Finite Fisher information implies $\sqrt\rho\in H^1$ and hence $\rho\in W^{1,1}$ with $\nabla\rho=\bs\rho$. By density of smooth functions in $W^{1,1}$, $\int|\rho(\bz-\bm\delta)-\rho(\bz)|d\bz\le\int_0^1\!\int|\bm\delta^\top\nabla\rho(\bz-\tau\bm\delta)|d\bz\,d\tau\le\norm{\bm\delta}\int\norm{\bs}\rho\le\norm{\bm\delta}(\tr\bm I)^{1/2}$.
\end{proof}

\begin{lemma}[Lipschitz property of the Bayes update]
\label{lem:bayes}
For probability measures $\mu,\nu$ and a bounded positive $g$, let $\Psi_g(\mu)=g\mu/\mu(g)$. Then $\TV(\Psi_g\mu,\Psi_g\nu)\le(\norm{g}_\infty/\nu(g))\TV(\mu,\nu)$.
\end{lemma}

\begin{proof}
Fix a set $A$, let $c'=\Psi_g\mu(A)$ and $h=g(\one_A-c')$. Then $\mu(h)=0$, hence $\Psi_g\mu(A)-\Psi_g\nu(A)=c'-\nu(g\one_A)/\nu(g)=-\nu(h)/\nu(g)=(\mu-\nu)(h)/\nu(g)$. Since $h$ takes values in $[-c'\norm{g}_\infty,(1-c')\norm{g}_\infty]$, an interval of length $\norm{g}_\infty$, we have $|(\mu-\nu)(h)|\le\norm{g}_\infty\TV(\mu,\nu)$. Taking the supremum over $A$ concludes.
\end{proof}

\rev{\begin{lemma}[Constrained nearness of minimizers]
\label{lem:argmin}
Let $\calC\subseteq\R^M$ be closed and convex, let $A:\calC\to\R$ be convex, and let $B(\bh)=\bh^\top\bz+\frac12\bh^\top\bOmega\bh$ with $\bOmega\succ0$. Let $\bm\beta=\arg\min_{\calC}B$ and $\bm\alpha\in\arg\min_{\calC}A$. For $\delta>0$, if
$\Delta(\delta):=\sup\{|A(\bh)-B(\bh)|:\bh\in\calC,\ \norm{\bh-\bm\beta}\le\delta\}<\frac14\lambda_{\min}(\bOmega)\delta^2$, then $\norm{\bm\alpha-\bm\beta}<\delta$.
\end{lemma}
\begin{proof}
For $\bh\in\calC$, the first-order optimality of $\bm\beta$ over the convex set $\calC$ gives $\ip{\nabla B(\bm\beta)}{\bh-\bm\beta}\ge0$, so $B(\bh)\ge B(\bm\beta)+\frac12\lambda_{\min}(\bOmega)\norm{\bh-\bm\beta}^2$. Hence, for $\bh\in\calC$ with $\norm{\bh-\bm\beta}=\delta$, $A(\bh)\ge B(\bh)-\Delta\ge B(\bm\beta)+\frac12\lambda_{\min}(\bOmega)\delta^2-\Delta\ge A(\bm\beta)+\frac12\lambda_{\min}(\bOmega)\delta^2-2\Delta>A(\bm\beta)$. If $\norm{\bm\alpha-\bm\beta}\ge\delta$, the point $\bh=\bm\beta+\theta(\bm\alpha-\bm\beta)$ with $\theta=\delta/\norm{\bm\alpha-\bm\beta}\in(0,1]$ lies in $\calC$ and satisfies $\norm{\bh-\bm\beta}=\delta$, while convexity and $A(\bm\alpha)\le A(\bm\beta)$ give $A(\bh)\le(1-\theta)A(\bm\beta)+\theta A(\bm\alpha)\le A(\bm\beta)$, a contradiction.
\end{proof}}

\section{Proof of Theorem \ref{thm:duality}}
\label{app_B}
\begin{proof}
By~\eqref{eq:ukappa} and the definition of $\bXi(\bx_t)$,
\begin{equation}
\nabla_{\bw}\bu_\kappa(\bw)=\bXi(\bx_t).
\end{equation}
Since
$\ell_\kappa(\bw)
=-\log p_\kappa(\bu_\kappa(\bw)\mid\bx_t)$,
the chain rule gives
\begin{align}
\nabla_{\bw}\ell_\kappa(\bw)
&=
-\bXi(\bx_t)^\top
\nabla_{\bu}\log p_\kappa(\bu\mid\bx_t)
\big|_{\bu=\bu_\kappa(\bw)}
\nonumber\\
&=
-\bXi(\bx_t)^\top
\bs_\kappa(\bu_\kappa(\bw)\mid\bx_t),
\end{align}

which proves~\eqref{eq:duality}. For $e=(i,j)$, using
$\bXi_e(\bx_t)=\gamma\bb_e\bb_e^\top\bx_t$
gives
\[
\bXi_e(\bx_t)^\top\bs_\kappa
=
\gamma(\bb_e^\top\bx_t)(\bb_e^\top\bs_\kappa)
=
\gamma(x_{t,i}-x_{t,j})(s_{\kappa,i}-s_{\kappa,j}),
\]
which establishes~\eqref{eq:duality-edge}.

At $\kappa=0$, the candidate transition density is
$p_{\bw}(\bx'\mid\bx)
=\rho^\star(\bx'-\bA(\bw)\bx\mid\bx)$.
Applying the same chain rule to
$\br=\bx'-\bA(\bw)\bx$ yields~\eqref{eq:duality-exact}.
\end{proof}

\section{Proof of Theorem \ref{thm:ident}}
\label{app_C}
\begin{proof}
\emph{(a)} By~\eqref{eq:maskobs} and Assumption~\ref{ass:mask}, $\E[y_{t,i}y_{t,j}]=\E[m_{t,i}m_{t,j}]\E[(x_{t,i}+v_{t,i})(x_{t,j}+v_{t,j})]=\pi^{(0)}_{ij}([\bSigma_x]_{ij}+[\bR_v]_{ij})$, and, \rev{since the masks are independent of $(\bx,\bv)$ and $\bv_{t+1}$, $\bv_t$ are centered, mutually independent, and independent of the states,} $\E[y_{t+1,i}y_{t,j}]=\pi^{(1)}_{ij}[\bGamma_1]_{ij}$. This proves~\eqref{eq:moments}. By~\eqref{eq:intro-model} and~\eqref{eq:centered}, $\bGamma_1=\bA(\bw^\star)\bSigma_x+\E[\E[\br_{t+1}\mid\bx_t]\bx_t^\top]=\bA(\bw^\star)\bSigma_x$, and $\bSigma_x\succ0$ gives $\bA(\bw^\star)=\bGamma_1\bSigma_x^{-1}$, hence $\bL(\bw^\star)$. Off-diagonal Laplacian entries determine the weights. \rev{The empirical moments $n^{-1}\sum_ty_{t,i}y_{t,j}$ and $n^{-1}\sum_ty_{t+1,i}y_{t,j}$ are ergodic averages of integrable functions of the joint process $(\bx_t,\br_{t+1},\bM_t,\bv_t)$, which is stationary and ergodic (Section~\ref{sec:theory}). Birkhoff's ergodic theorem and the continuity of matrix inversion at $\bSigma_x\succ0$ give strong consistency.}

\emph{(b)} Conditionally on $\bx_t=\bx$, $\bu_\kappa(\bw^\star)=\br_{t+1}+\kappa\beps\sim p_\kappa(\cdot\mid\bx)$, and $\bu_\kappa(\bw)=\bu_\kappa(\bw^\star)+\bm\delta$ with $\bm\delta=\bXi(\bx)(\bw-\bw^\star)$, whose density is $q(\cdot)=p_\kappa(\cdot-\bm\delta\mid\bx)$. Jensen's inequality yields $\log p_\kappa(\bu)\ge-\frac N2\log(2\pi\kappa^2)-\E\norm{\bu-\br}^2/(2\kappa^2)$, \rev{and $p_\kappa\le(2\pi\kappa^2)^{-N/2}$,} so the cross-entropy is finite under Assumption~\ref{ass:stat}. Hence $R_\kappa(\bw;\bx)=h(q)+\KL(q\|p_\kappa)$, where $h$ is the differential entropy. Translation invariance gives $h(q)=h(p_\kappa)=R_\kappa(\bw^\star;\bx)$, proving~\eqref{eq:KL}. If $R_{\kappa,\pi}(\bw)=R_{\kappa,\pi}(\bw^\star)$, the KL divergence vanishes $\pi$-a.e., so $p_\kappa(\cdot-\bm\delta(\bx))=p_\kappa(\cdot)$ a.e.; by Lemma~\ref{lem:transl} in Appendix~\ref{app:aux}, $\bm\delta(\bx)=\bzero$ $\pi$-a.e. Lemma~\ref{lem:excitation} then gives $2\gamma^2\lambda_{\min}(\bSigma_\pi)\norm{\bw-\bw^\star}^2\le\int\norm{\bm\delta}^2d\pi=0$.

\emph{(c)} \rev{Per sample, the chain rule through the affine map~\eqref{eq:ukappa} gives $\nabla_{\bw}\ell_\kappa=-\bXi^\top\bs_\kappa(\bu_\kappa(\bw))$ and $\nabla^2_{\bw}\ell_\kappa=\bXi^\top\calH_\kappa(\bu_\kappa(\bw))\bXi$. By~\eqref{eq:dominate}, these derivatives are dominated, uniformly over any compact neighborhood of a given $\bw$, by an integrable function, and the differentiation lemma (dominated convergence applied to difference quotients) permits differentiating twice under the integral, which yields~\eqref{eq:hess}.} At $\bw=\bw^\star$, $\bu_\kappa\sim p_\kappa(\cdot\mid\bx)$ and the Fisher identity $\int\calH_\kappa\,dp_\kappa=\tI_\kappa$ yields the first equality in~\eqref{eq:fisherhess}. Lemma~\ref{lem:CR} in Appendix~\ref{app:aux}, applied to $p_\kappa$ with covariance $\bSigma_r(\bx)+\kappa^2\bI$, gives $\tI_\kappa(\bx)\succeq(\bSigma_r(\bx)+\kappa^2\bI)^{-1}$ \rev{for $\pi$-a.e.\ $\bx$. Hence, for $\bu\ne\bzero$, $\bu^\top\nabla^2R_{\kappa,\pi}(\bw^\star)\bu\ge\int\norm{\bXi(\bx)\bu}^2/(\lambda_{\max}(\bSigma_r(\bx))+\kappa^2)\,\pi(d\bx)$. The integrand is nonnegative and is positive on $\{\bXi(\bx)\bu\ne\bzero\}$, which has positive $\pi$-probability because $\int\norm{\bXi\bu}^2d\pi\ge2\gamma^2\lambda_{\min}(\bSigma_\pi)\norm{\bu}^2>0$ by Lemma~\ref{lem:excitation}; this proves positive definiteness. If $\bar\lambda_r<\infty$, the denominator is at most $\bar\lambda_r+\kappa^2$, and Lemma~\ref{lem:excitation} gives~\eqref{eq:fisherhess}.}

\rev{\emph{Verification of~\eqref{eq:dominate} under Assumption~\ref{ass:resid}.} By~\eqref{eq:Hmix} in Appendix~\ref{app_D}, $\norm{\calH_\kappa(\bu\mid\bx)}\le\ell_\kappa:=\max\{\kappa^{-2},D^2\kappa^{-4}\}$ for all $\bu$, and by Lemma~\ref{lem:gmix}, $\norm{\bs_\kappa(\bu\mid\bx)}\le\kappa^{-2}(D+\norm{\bu})$. For $\bw$ in a compact $\calC$ with $c_{\calC}=\sup_{\bw\in\calC}\norm{\bw-\bw^\star}$, $\norm{\bu_\kappa(\bw)}\le\norm{\br}+\kappa\norm{\beps}+c_{\calC}\norm{\bXi(\bx)}$. Since $\norm{\bXi(\bx)}\le\gamma\sqrt{2N}\norm{\bx}$, the Cauchy--Schwarz inequality and the finiteness of $\int\norm{\bx}^2d\pi$ and $\int\tr\bSigma_r\,d\pi$ bound the expectation in~\eqref{eq:dominate}.}
\end{proof}

\section{Proof of Theorem \ref{thm:convex}}
\label{app_D}
\begin{proof}
Conditionally on $\bx$, $\bu=\bc+(\bg+\kappa\beps)$ where $\bg+\kappa\beps\sim\calN(\bzero,\bLambda^{-1})$, $\bLambda=(\bSigma_g(\bx)+\kappa^2\bI)^{-1}$. By Lemma~\ref{lem:gmix} in Appendix~\ref{app:aux},
\begin{equation}
\calH_\kappa(\bu\mid\bx)=\bLambda-\bLambda\Cov(\bc\mid\bu,\bx)\bLambda .
\label{eq:Hmix}
\end{equation}
\emph{(a)} Since $\norm{\bc}\le D$, $\Cov(\bc\mid\bu,\bx)\preceq\E[\bc\bc^\top\mid\bu,\bx]\preceq D^2\bI$, hence $\calH_\kappa\succeq\bLambda-D^2\bLambda^2$. The eigenvalues $\lambda$ of $\bLambda$ lie in $[(\bar\sigma_g^2+\kappa^2)^{-1},\kappa^{-2}]$, and $\lambda-D^2\lambda^2=\lambda(1-D^2\lambda)\ge\lambda(1-D^2/\kappa^2)$. Thus $\calH_\kappa\succeq c_\kappa\bI$ with $c_\kappa=(1-D^2/\kappa^2)/(\bar\sigma_g^2+\kappa^2)>0$ for every $\bu$. Also $\calH_\kappa\preceq\bLambda\preceq\kappa^{-2}\bI$. Inserting into~\eqref{eq:hess} and using Lemma~\ref{lem:excitation} gives~\eqref{eq:globalconv}--\eqref{eq:muL}. \rev{For every $\kappa>0$, the same argument gives $-D^2\kappa^{-4}\bI\preceq\bLambda-D^2\bLambda^2\preceq\calH_\kappa\preceq\kappa^{-2}\bI$, hence $\norm{\calH_\kappa}\le\ell_\kappa=\max\{\kappa^{-2},D^2\kappa^{-4}\}$; since $\nabla_{\bu}\bs_\kappa=-\calH_\kappa$, the score $\bs_\kappa(\cdot\mid\bx)$ is globally $\ell_\kappa$-Lipschitz.}

\emph{(b)} By Lemma~\ref{lem:gmix}, the directional derivative of $\Cov(\bc\mid\bu,\bx)$ along $\bh$ is $\E[\hat\bc\hat\bc^\top(\hat\bc^\top\bLambda\bh)\mid\bu,\bx]$ with $\hat\bc=\bc-\E[\bc\mid\bu,\bx]$, $\norm{\hat\bc}\le2D$. Its norm is at most $8D^3\kappa^{-2}\norm{\bh}$, so $\bu\mapsto\calH_\kappa(\bu\mid\bx)$ is Lipschitz with constant $\norm{\bLambda}^2\cdot8D^3\kappa^{-2}\le8D^3\kappa^{-6}$. Since $\bu_\kappa(\bw)-\bu_\kappa(\bw')=\bXi(\bx)(\bw-\bw')$, \eqref{eq:hess} gives $\norm{\nabla^2R(\bw)-\nabla^2R(\bw')}\le8D^3\kappa^{-6}\int\norm{\bXi}^3d\pi\,\norm{\bw-\bw'}$, and Lemma~\ref{lem:excitation} bounds $\norm{\bXi}^3$.

\emph{(c)} By Theorem~\ref{thm:ident}(c) with $\bar\lambda_r\le D^2+\bar\sigma_g^2$, $\nabla^2R_{\kappa,\pi}(\bw^\star)\succeq\mu^{\rm loc}_\kappa\bI$. By~\eqref{eq:hesslip}, for $\norm{\bw-\bw^\star}\le r_\kappa$, $\lambda_{\min}(\nabla^2R_{\kappa,\pi}(\bw))\ge\mu^{\rm loc}_\kappa-L_{3,\kappa}r_\kappa=\mu^{\rm loc}_\kappa/2$.

\rev{\emph{(d)} If $D=0$, then $\bc=\bzero$, and~\eqref{eq:Hmix} gives $\calH_\kappa(\bu\mid\bx)=\bLambda(\bx)$, independent of $\bu$. Hence~\eqref{eq:hess} is constant in $\bw$, $R_{\kappa,\pi}$ is quadratic, $L_{3,\kappa}=0$, and part (a) with $D=0$ gives the lower bound for every $\kappa>0$ and every $\bw\in\R^M$.}
\end{proof}

\section{Proofs for Section~\ref{sec:eff}}
\label{app_E}

\rev{\emph{Preliminaries.} Under Assumption~\ref{ass:resid}, Appendix~\ref{app_D} gives, for all $\bu,\bu'$ and $\bx$, $\norm{\calH_\kappa(\bu\mid\bx)}\le\ell_\kappa$ and $\norm{\calH_\kappa(\bu\mid\bx)-\calH_\kappa(\bu'\mid\bx)}\le8D^3\kappa^{-6}\norm{\bu-\bu'}$. Write $\ell_t(\bw)=\ell_\kappa(\bw;\bx_t,\bx_{t+1},\beps_t)$ and $c_D=8D^3\kappa^{-6}$. Then $\nabla\ell_t(\bw)=-\bXi(\bx_t)^\top\bs_\kappa(\bu_t(\bw)\mid\bx_t)$, $\nabla^2\ell_t(\bw)=\bXi^\top\calH_\kappa(\bu_t(\bw)\mid\bx_t)\bXi$, and
\begin{equation}
\norm{\nabla^2\ell_t(\bw)-\nabla^2\ell_t(\bw')}\le c_D\norm{\bXi(\bx_t)}^3\norm{\bw-\bw'}.
\label{eq:persample}
\end{equation}
The process $(\bx_t,\bx_{t+1},\beps_t)_t$ is stationary and ergodic, as in Section~\ref{sec:theory}. Let $\calF_t=\sigma(\bx_{0:t},\beps_{0:t-1})$. Since the smoothed density is a Gaussian mixture with covariance between $\kappa^2\bI$ and $(\bar\sigma_g^2+\kappa^2)\bI$ and centers in $\calB(\bzero,D)$, $|\ell_t(\bw)|\le C(1+\norm{\bu_t(\bw)}^2)$ for a constant $C$, so $\E|\ell_t(\bw)|<\infty$ for every $\bw$.}

\begin{proof}[Proof of Theorem~\ref{thm:eff}]
\emph{(a)} \rev{Let $\calS=\calS^\star$ and} $\bpsi_t=\nabla_{\calS}\ell_t(\bw^\star)=-\bXi_{\calS}(\bx_t)^\top\bs_\kappa(\br_{t+1}+\kappa\beps_t\mid\bx_t)$ by Theorem~\ref{thm:duality}. Since $\br_{t+1}$ is independent of the past given $\bx_t$ and $\beps_t$ is fresh, $\br_{t+1}+\kappa\beps_t\mid\calF_t\sim p_\kappa(\cdot\mid\bx_t)$; hence $\E[\bpsi_t\mid\calF_t]=\bzero$ and $\E[\bpsi_t\bpsi_t^\top\mid\calF_t]=\bXi_{\calS}^\top\tI_\kappa\bXi_{\calS}$. Since $\tI_\kappa\preceq\kappa^{-2}\bI$ and $\E\norm{\bx}^2<\infty$, $(\bpsi_t)$ is a stationary ergodic square-integrable martingale difference sequence, and Billingsley's martingale CLT~\cite{billingsley1961lindeberg,hall1980martingale} gives $n^{-1/2}\sum_t\bpsi_t\Rightarrow\calN(\bzero,\bOmega_{\kappa,\calS})$.

\rev{\emph{Consistency for $\kappa>D$.} By Appendix~\ref{app_D}(a), $\calH_\kappa(\bu\mid\bx)\succeq c_\kappa\bI$ for every $\bu$, so each $\ell_t$, and hence the empirical objective $\widehat R_n=n^{-1}\sum_t\ell_t$, is convex on $\R^M$. By the ergodic theorem, $\widehat R_n(\bw)\to R_\kappa(\bw)$ a.s.\ for every $\bw$ in a countable dense set; pointwise convergence of convex functions on a dense set implies uniform convergence on compacts~\cite[Thm.~10.8]{rockafellar1970convex}. Since $\calW_{\calS}$ is compact and $\bw^\star$ is the unique minimizer of $R_\kappa$ on $\calW_{\calS}$ (Theorem~\ref{thm:ident}(b)), $\widehat\bw_{n,\calS}\to\bw^\star$ a.s.}

\rev{\emph{Expansion.} Since $\bw^\star_{\calS}\in\relint\calW_{\calS}$ and $\widehat\bw_{n,\calS}\to\bw^\star$ in probability, $\widehat\bw_{n,\calS}\in\relint\calW_{\calS}$ with probability tending to one, and then $\nabla_{\calS}\widehat R_n(\widehat\bw_{n,\calS})=\bzero$. With $\bDelta_n=(\widehat\bw_{n,\calS}-\bw^\star)_{\calS}$,
$\bzero=n^{-1}\sum_t\bpsi_t+\bar{\calH}_n\bDelta_n$ with $\bar{\calH}_n=\int_0^1n^{-1}\sum_t\nabla^2_{\calS\calS}\ell_t(\bw^\star+s(\widehat\bw_{n,\calS}-\bw^\star))\,ds$. By~\eqref{eq:persample}, $\norm{\bar{\calH}_n-\widehat{\calH}_n}\le\frac12\widehat L_n\norm{\bDelta_n}$, where $\widehat{\calH}_n=n^{-1}\sum_t\nabla^2_{\calS\calS}\ell_t(\bw^\star)$ and $\widehat L_n=n^{-1}\sum_tc_D\norm{\bXi(\bx_t)}^3$ converges a.s.\ because $\E\norm{\bx_t}^3<\infty$. By the ergodic theorem and the conditional Fisher identity $\E[\calH_\kappa(\br+\kappa\beps\mid\bx)\mid\bx]=\tI_\kappa(\bx)$, $\widehat{\calH}_n\to\bOmega_{\kappa,\calS}$ a.s., which is positive definite by Theorem~\ref{thm:ident}(c) applied to vectors supported on $\calS$. Hence $\bar{\calH}_n\to\bOmega_{\kappa,\calS}$ in probability, and Slutsky's lemma gives $\sqrt n\bDelta_n\Rightarrow\calN(\bzero,\bOmega_{\kappa,\calS}^{-1}\bOmega_{\kappa,\calS}\bOmega_{\kappa,\calS}^{-1})=\calN(\bzero,\bV_{\kappa,\calS})$.}

\emph{(b)} For $\kappa'>\kappa$, $\br+\kappa'\beps\overset{d}{=}(\br+\kappa\beps)+\sqrt{\kappa'^2-\kappa^2}\,\beps'$ with independent $\beps'$, so Lemma~\ref{lem:fishconv}(i) yields $\tI_{\kappa'}(\bx)\preceq\tI_\kappa(\bx)$. Congruence by \rev{$\bXi_{\calS}$}, expectation, and matrix inversion reverse the order.

\emph{(c)} By Lemma~\ref{lem:fishconv}(i)--(ii), $\tI_\kappa(\bx)\preceq\bm I_0(\bx)$ and $\tI_\kappa(\bx)\to\bm I_0(\bx)$ as $\kappa\downarrow0$; under Assumption~\ref{ass:resid} with $\underline\sigma_g>0$, $\bm I_0\preceq\bSigma_g^{-1}\preceq\underline\sigma_g^{-2}\bI$ by the same lemma. Dominated convergence gives \rev{$\bOmega_{\kappa,\calS}\to\bOmega_{0,\calS}$}, which is positive definite \rev{since $\bOmega_{0,\calS}\succeq\bOmega_{\kappa,\calS}$}, and continuity of inversion gives \rev{$\bV_{\kappa,\calS}\to\bV_{0,\calS}$}. By Theorem~\ref{thm:duality} at $\kappa=0$, the per-pair score of the \rev{submodel $\{p_{\bw}:\bw\in\calW_{\calS}\}$ is $\bXi_{\calS}^\top\nabla\log\rho^\star$, whose covariance is $\bOmega_{0,\calS}$. The submodel is a location family in the residual with finite Fisher information and is therefore differentiable in quadratic mean; for the stationary ergodic Markov chain $(\bx_t)$ this yields local asymptotic normality with information $\bOmega_{0,\calS}$, and the H\'ajek--Le Cam convolution and local asymptotic minimax theorems~\cite{vandervaart1998asymptotic} identify $\bV_{0,\calS}$ as the efficiency bound for regular estimators.}

\emph{(d)} Lemma~\ref{lem:CR} gives $\tI_\kappa(\bx)\succeq(\bSigma_r(\bx)+\kappa^2\bI)^{-1}$, hence \rev{$\bOmega_{\kappa,\calS}\succeq\bV_{{\rm GLS},\calS}(\kappa)^{-1}$ and $\bV_{\kappa,\calS}\preceq\bV_{{\rm GLS},\calS}(\kappa)$}. Since the mean $\bA(\bw)\bx=\beta\bx-\bXi(\bx)\bw$ is linear in $\bw$ by~\eqref{eq:Xi}, \rev{the population criterion of a Gaussian QMLE with working covariance $\bm W$ equals $\E\norm{\br}^2_{\bm W^{-1}}+(\bw-\bw^\star)^\top\bm H_{\bm W}(\bw-\bw^\star)$ by~\eqref{eq:centered}, so the QMLE is consistent and} has the sandwich covariance $\bV_{\bm W}=\bm H_{\bm W}^{-1}\bm S_{\bm W}\bm H_{\bm W}^{-1}$ with $\bm H_{\bm W}=\E[\bXi_{\calS}^\top\bm W^{-1}\bXi_{\calS}]$ and $\bm S_{\bm W}=\E[\bXi_{\calS}^\top\bm W^{-1}\bSigma_r\bm W^{-1}\bXi_{\calS}]$~\cite{white1982maximum}. Put $\bm X=\bSigma_r^{1/2}\bm W^{-1}\bXi_{\calS}$ and $\bm Y=\bSigma_r^{-1/2}\bXi_{\calS}$. The matrix Cauchy--Schwarz inequality $\E[\bm X^\top\bm X]\succeq\E[\bm X^\top\bm Y]\E[\bm Y^\top\bm Y]^{-1}\E[\bm Y^\top\bm X]$ reads $\bm S_{\bm W}\succeq\bm H_{\bm W}\bV_{{\rm GLS},\calS}(0)\bm H_{\bm W}$, i.e., $\bV_{\bm W}\succeq\bV_{{\rm GLS},\calS}(0)$.

\rev{For the gap, write $\bm A=\bOmega_{0,\calS}$ and $\bm B=\bV_{{\rm GLS},\calS}(0)^{-1}$, so that $\bm A-\bm B=\E[\bXi_{\calS}^\top\bm F\bXi_{\calS}]\succeq\bzero$. The expectation of a positive semidefinite matrix vanishes if and only if the integrand vanishes a.s., which proves the equivalence. Otherwise, $\bm B^{-1}-\bm A^{-1}=\bm A^{-1}(\bm A-\bm B)\bm B^{-1}$ is nonzero, as a product of a nonzero matrix with invertible ones, and positive semidefinite by the antitonicity of inversion; a nonzero positive semidefinite matrix has positive trace. Finally, $\operatorname{range}\bXi_{\calS}(\bx)\subseteq\one^\perp$ gives $\bm P\bXi_{\calS}=\bXi_{\calS}$, so $\bXi_{\calS}^\top\bm F\bXi_{\calS}=\bXi_{\calS}^\top\bm P\bm F\bm P\bXi_{\calS}\succeq\varphi_0\bXi_{\calS}^\top\bXi_{\calS}$, and Lemma~\ref{lem:excitation}, applied with $\pi=\pi_\infty$ to vectors supported on $\calS$, gives $\E[\bXi_{\calS}^\top\bXi_{\calS}]\succeq2\gamma^2\lambda_{\min}(\bSigma_x)\bI$, which is~\eqref{eq:strictgap}. Note that $\bm P\bm F\bm P\ne\bzero$ on a set of positive probability does not suffice for a strict gap, because the gap is invisible on states with $\bXi_{\calS}(\bx)=\bzero$; the condition must be stated through $\bXi_{\calS}^\top\bm F\bXi_{\calS}$.}
\end{proof}

\rev{\begin{proof}[Proof of Proposition~\ref{prop:boundary}]
\emph{Step 1: local structure.} The constraints $\one^\top\bw\le R_w$ and $\gamma\lambda_{\max}(\bL(\bw))\le1+\beta-\delta_A$ are strict at $\bw^\star$ and continuous, and $w^\star_e>0$ for $e\in\calS^\star$. Hence there is $r_0>0$ with $\calW\cap\calB(\bw^\star,r_0)=(\bw^\star+\calT^\star)\cap\calB(\bw^\star,r_0)$, and with $\calT_n=\sqrt n(\calW-\bw^\star)$, $\calT_n\cap\calB(\bzero,r_0\sqrt n)=\calT^\star\cap\calB(\bzero,r_0\sqrt n)$.

\emph{Step 2: local convexity.} Let $A_n(\bh)=\sum_t[\ell_t(\bw^\star+\bh/\sqrt n)-\ell_t(\bw^\star)]$, so that $\bm\alpha_n:=\sqrt n(\widehat\bw_n-\bw^\star)$ minimizes $A_n$ over $\calT_n$. By~\eqref{eq:persample}, $\nabla^2A_n(\bh)\succeq\widehat{\calH}_n-\widehat L_n\norm{\bh}n^{-1/2}\bI$ with $\widehat{\calH}_n=n^{-1}\sum_t\nabla^2\ell_t(\bw^\star)\to\bOmega_\kappa$ and $\widehat L_n\to\bar L<\infty$ a.s. With $r_1=\min\{r_0,\lambda_{\min}(\bOmega_\kappa)/(4\bar L)\}$ ($r_1=r_0$ if $\bar L=0$), $A_n$ is convex on $\calC_n:=\calT^\star\cap\calB(\bzero,r_1\sqrt n)$ with probability tending to one; if $\kappa>D$, $A_n$ is convex on $\R^M$.

\emph{Step 3: quadratic approximation.} Taylor's formula and~\eqref{eq:persample} give $A_n(\bh)=\bh^\top\bz_n+\frac12\bh^\top\widehat{\calH}_n\bh+\varrho_n(\bh)$ with $\bz_n=n^{-1/2}\sum_t\nabla\ell_t(\bw^\star)$ and $|\varrho_n(\bh)|\le\widehat L_n\norm{\bh}^3/(6\sqrt n)$. By the proof of Theorem~\ref{thm:eff}(a) with $\calS=\calE_c$, $\bz_n\Rightarrow\calN(\bzero,\bOmega_\kappa)$. With $B_n(\bh)=\bh^\top\bz_n+\frac12\bh^\top\bOmega_\kappa\bh$, for every $c>0$, $\sup_{\norm{\bh}\le c}|A_n(\bh)-B_n(\bh)|\le\frac12c^2\norm{\widehat{\calH}_n-\bOmega_\kappa}+\widehat L_nc^3/(6\sqrt n)\to0$ a.s.

\emph{Step 4: convergence.} Let $\bm\beta_n=\arg\min_{\calT^\star}B_n$, the $\bOmega_\kappa$-projection of $-\bOmega_\kappa^{-1}\bz_n$ onto $\calT^\star$; it is unique and a continuous function of $\bz_n$, so $\bm\beta_n\Rightarrow\widehat\bh$ by the continuous mapping theorem and the symmetry of $\calN(\bzero,\bOmega_\kappa)$, and $\bm\beta_n=O_p(1)$. Consistency of $\widehat\bw_n$ (proved as in Theorem~\ref{thm:eff}(a) if $\kappa>D$, assumed otherwise) gives $\norm{\bm\alpha_n}<r_1\sqrt n/2$ with probability tending to one; on this event $\bm\alpha_n\in\calC_n\subseteq\calT_n$, so $\bm\alpha_n$ minimizes $A_n$ over $\calC_n$, and $\bm\beta_n\in\calC_n$ minimizes $B_n$ over $\calC_n$. Lemma~\ref{lem:argmin} with $\calC=\calC_n$ yields, for every $\delta>0$ and $c>0$,
$\Prob(\norm{\bm\alpha_n-\bm\beta_n}\ge\delta)\le\Prob\big(\sup_{\norm{\bh}\le c}|A_n-B_n|\ge\tfrac14\lambda_{\min}(\bOmega_\kappa)\delta^2\big)+\Prob(\norm{\bm\beta_n}+\delta>c)+o(1)$.
Letting $n\to\infty$ and then $c\to\infty$ gives $\bm\alpha_n-\bm\beta_n\to\bzero$ in probability, and Slutsky's lemma gives~\eqref{eq:conelimit}.

\emph{Step 5: structure of the limit.} If $\calS^\star=\calE_c$, then $\calT^\star=\R^M$ and $\widehat\bh=\bOmega_\kappa^{-1}\bm Z$. Otherwise, write $\calS=\calS^\star$ and $\calS^c=\calE_c\setminus\calS^\star$. The objective in~\eqref{eq:conelimit} is strictly convex and $\calT^\star$ is polyhedral, so the Karush--Kuhn--Tucker conditions are necessary and sufficient: $\widehat\bh_{\calS^c}=\bzero$ if and only if $\bh_{\calS}=\bOmega_{\calS\calS}^{-1}\bm Z_{\calS}$ and $(\bOmega_\kappa\bh-\bm Z)_{\calS^c}\succeq\bzero$, i.e., if and only if $\bm U:=\bm Z_{\calS^c}-\bOmega_{\calS^c\calS}\bOmega_{\calS\calS}^{-1}\bm Z_{\calS}\preceq\bzero$, where $\bOmega=\bOmega_\kappa$. The vector $\bm U$ is Gaussian with the positive definite Schur-complement covariance and is independent of $\bm Z_{\calS}$. Hence $\Prob(\calA^\star)=\Prob(\bm U\preceq\bzero)>0$, and, conditionally on $\calA^\star$, $\widehat\bh_{\calS}=\bOmega_{\calS\calS}^{-1}\bm Z_{\calS}\sim\calN(\bzero,\bOmega_{\calS\calS}^{-1})$ with $\bOmega_{\calS\calS}=\bOmega_{\kappa,\calS^\star}=\bV_{\kappa,\calS^\star}^{-1}$. Finally, for $e\in\calS^c$, $\Prob(\widehat h_e=0)\ge\Prob(\calA^\star)>0$, whereas $\Prob(\widehat h_e>0)\ge\Prob(\bOmega_\kappa^{-1}\bm Z\in\operatorname{int}\calT^\star)>0$; a random variable with an atom that is not almost surely constant is not Gaussian.
\end{proof}}

\rev{\begin{proof}[Proof of Corollary~\ref{cor:Kdraw}]
Let $\bpsi_{t,k}$ be the score of draw $k$ and $\bar\bpsi^{(K)}_t=K^{-1}\sum_k\bpsi_{t,k}$. Each $\bpsi_{t,k}$ satisfies $\E[\bpsi_{t,k}\mid\calF_t]=\bzero$ as in Theorem~\ref{thm:eff}(a), and, given $(\bx_t,\bx_{t+1})$, the draws are i.i.d., so $\E[\bpsi_{t,k}\bpsi_{t,k'}^\top\mid\bx_t,\bx_{t+1}]=\bar\bpsi_t\bar\bpsi_t^\top$ for $k\neq k'$. Hence $\E[\bar\bpsi^{(K)}_t\bar\bpsi^{(K)\top}_t]=K^{-1}\bOmega+(1-K^{-1})\bm S_\infty=\bm S_K$, and the martingale CLT gives $n^{-1/2}\sum_t\bar\bpsi^{(K)}_t\Rightarrow\calN(\bzero,\bm S_K)$. The averaged loss satisfies~\eqref{eq:persample} and its Hessian at $\bw^\star$ has mean $\bOmega$, so the expansion of Theorem~\ref{thm:eff}(a) yields~\eqref{eq:Kdraw}.

By conditional Jensen, $\bm S_\infty\preceq\E[\bpsi_t\bpsi_t^\top]=\bOmega$, so $\bm S_K=\bm S_\infty+K^{-1}(\bOmega-\bm S_\infty)$ is nonincreasing in $K$, $\bm S_1=\bOmega$, and $\bV^{(1)}_\kappa=\bOmega^{-1}=\bV_{\kappa,\calS^\star}$. For the lower bound, let $\bpsi_{0,t}=-\bXi_{\calS}(\bx_t)^\top\bs_0(\br_{t+1}\mid\bx_t)$, so that $\E[\bpsi_{0,t}\bpsi_{0,t}^\top]=\bOmega_{0,\calS}=\bV_{0,\calS}^{-1}$. Lemma~\ref{lem:fishconv}(i), applied conditionally on $\bx_t$ with $\bm X=\br_{t+1}$ and $\bm Y=\kappa\beps$, gives $\bs_\kappa(\bu\mid\bx_t)=\E[\bs_0(\br_{t+1}\mid\bx_t)\mid\bu,\bx_t]$ with $\bu=\br_{t+1}+\kappa\beps$. Hence $\E[\bpsi_t\bpsi_{0,t}^\top\mid\bx_t]=\bXi_{\calS}^\top\E[\bs_\kappa(\bu)\bs_\kappa(\bu)^\top\mid\bx_t]\bXi_{\calS}=\bXi_{\calS}^\top\tI_\kappa\bXi_{\calS}$, and, since $\bpsi_{0,t}$ is $(\bx_t,\bx_{t+1})$-measurable, $\E[\bar\bpsi_t\bpsi_{0,t}^\top]=\E[\bpsi_t\bpsi_{0,t}^\top]=\bOmega$. The matrix Cauchy--Schwarz inequality gives $\bm S_\infty\succeq\bOmega\bV_{0,\calS}\bOmega$, hence $\bV^{(K)}_\kappa\succeq\bOmega^{-1}\bm S_\infty\bOmega^{-1}\succeq\bV_{0,\calS}$.
\end{proof}}

\section{Proof of Proposition \ref{prop:detrend} and \ref{prop:transfer-bias}}
\label{app_1}
\begin{proof}
Conditionally on $\bx$, $\bu_\kappa(\bw)$ has density $q_{\bw}=p_\kappa(\cdot-\bXi(\bx)(\bw-\bw^\star)\mid\bx)$, whereas $\widehat p_\kappa=p_\kappa(\cdot-\bb(\bx)\mid\bx)$. The cross-entropy is finite by the argument of Appendix~\ref{app_C}(b), so $\widehat R_\kappa(\bw;\bx)=h(q_{\bw})+\KL(q_{\bw}\|\widehat p_\kappa)$, and $h(q_{\bw})=h(p_\kappa(\cdot\mid\bx))$ does not depend on $\bw$. (i) The KL term vanishes for all $\bx$ at $\bw=\bw^\star+\bm\delta$. Conversely, if it vanishes $\pi$-a.e.\ at some $\bw$, then $p_\kappa(\cdot-\bXi(\bx)(\bw-\bw^\star-\bm\delta)\mid\bx)=p_\kappa(\cdot\mid\bx)$ a.e., so Lemma~\ref{lem:transl} gives $\bXi(\bx)(\bw-\bw^\star-\bm\delta)=\bzero$ $\pi$-a.e., and Lemma~\ref{lem:excitation} gives $\bw=\bw^\star+\bm\delta$. (ii) For Gaussians with common covariance, $\KL(q_{\bw}\|\widehat p_\kappa)=\frac12\norm{\bXi(\bx)(\bw-\bw^\star)-\bb(\bx)}^2_{\bLambda_\kappa(\bx)}$. Its $\pi$-average is a weighted least-squares criterion whose normal equations give~\eqref{eq:gaussbias}; $\E_\pi[\bXi^\top\bLambda_\kappa\bXi]\succ0$ because $\bLambda_\kappa(\bx)\succ0$ a.s.\ and $\E_\pi\norm{\bXi\bu}^2>0$ for $\bu\ne\bzero$ by Lemma~\ref{lem:excitation}.
\end{proof}

\begin{proof}
By~\eqref{eq:dominate}, differentiation and expectation commute, so $\nabla R_{\kappa,\pi}(\bw)=-\E[\bXi^\top\bs_\kappa(\bu_\kappa(\bw)\mid\bx)]$ by Theorem~\ref{thm:duality}. By the Cauchy--Schwarz inequality,
$\norm{\widehat\bg_\kappa(\bw)-\nabla R_{\kappa,\pi}(\bw)}\le\E[\norm{\bXi}\norm{\widehat\bs_\kappa-\bs_\kappa}]\le C_\Xi\varepsilon_{\rm tar}$.
Since $\bw^\star$ minimizes $R_{\kappa,\pi}$ over $\R^M$ (Theorem~\ref{thm:ident}(b)), $\nabla R_{\kappa,\pi}(\bw^\star)=\bzero$, and strong convexity gives
$\mu\norm{\widetilde\bw-\bw^\star}^2\le\ip{\nabla R_{\kappa,\pi}(\widetilde\bw)}{\widetilde\bw-\bw^\star}
=\ip{\nabla R_{\kappa,\pi}(\widetilde\bw)-\widehat\bg_\kappa(\widetilde\bw)}{\widetilde\bw-\bw^\star}+\ip{\widehat\bg_\kappa(\widetilde\bw)}{\widetilde\bw-\bw^\star}$.
The last term is nonpositive by the variational inequality with $\bw=\bw^\star$, and the first is at most $C_\Xi\varepsilon_{\rm tar}\norm{\widetilde\bw-\bw^\star}$.
\end{proof}

\section{Proof of Proposition \ref{prop:local}}
\label{app_2}
\begin{proof}
\emph{(i)} $\E_q[W_q]=C_{\bw}$ and $\E_q[W_q^2]=C_{\bw}^2\int(q^\star)^2/q$, so $\Var_q[W_q]/C_{\bw}^2=\int(q^\star)^2/q-1=\chi^2$. Further, $2\TV=\int|q^\star/q-1|q\le(\int(q^\star/q-1)^2q)^{1/2}$ by Cauchy--Schwarz. \rev{If $q^\star(A)>0=q(A)$ for some $A$, the density ratio is infinite on a set of positive $q^\star$-measure and $\chi^2=\infty$.}

\emph{(ii)} Given the past, $\widehat\Gamma_{t+1}$ is the empirical measure of $J$ i.i.d.\ draws from $\widetilde\Gamma=\sum_j\widehat p_j\delta_{\bx^{(j)}}\otimes\widehat q_j$, hence $\E[\norm{\widehat\Gamma_{t+1}(\bm f)-\widetilde\Gamma(\bm f)}^2\mid\text{past}]\le J^{-1}$. Next, \rev{since $\norm{\mu(\bm f)-\nu(\bm f)}\le\int\norm{\bm f}d|\mu-\nu|\le2\TV(\mu,\nu)$ and the TV distance between two mixtures is at most the TV distance between the mixing weights plus the largest TV distance between components,} $\norm{\widetilde\Gamma(\bm f)-\Phi(\bm f)}\le2\TV(\widetilde\Gamma,\Phi)\le2\TV(\widehat{\bm p},\bm p)+2\max_j\TV(\widehat q_j,q^\star_j)$. With $\widehat p_j=p_j(1+\delta_j)/(1+\bar\delta)$, $|\delta_j|\le\varepsilon_C$, and $\bar\delta=\sum_jp_j\delta_j$, we obtain $\TV(\widehat{\bm p},\bm p)=\frac12\sum_jp_j|\delta_j-\bar\delta|/(1+\bar\delta)\le\varepsilon_C/(1-\varepsilon_C)$. \rev{Both bias terms are measurable with respect to the past, so} Minkowski's inequality concludes. \rev{No absolute continuity between $\widehat q_j$ and $q^\star_j$ is required.}

\rev{\emph{(iii)} Couple $\br\sim\nu$ with $T_D(\br)$; the two differ only on $\{\norm{\br}>D\}$, so $\TV\le\nu(\norm{\br}>D)$. Clipping acts on residuals and the affine map $\br\mapsto\bA(\bw)\bx^{(j)}+\br$ is a bijection, so it suffices to argue in residual coordinates. If the exact residual posterior $\nu^\star$ is supported in $\calB(\bzero,D)$, then $(T_D)_\#\nu^\star=\nu^\star$ and the data-processing inequality gives $\TV((T_D)_\#\widehat\nu,\nu^\star)\le\TV(\widehat\nu,\nu^\star)$. The kernel bound follows from the same coupling. Finally, $\E[\br_{t+1}\mid\bx_t]=\bzero$ and $\norm{T_D(\br)-\br}=(\norm{\br}-D)_+$.}
\end{proof}

\section{Proof of Proposition \ref{prop:PEmean}}
\label{app_3}
\rev{\begin{proof}
Let $\calG_t=\sigma(\by_{1:t},\bx_{t-1})$. As in Proposition~\ref{prop:PE}, $\E[\Cov(\bx_t\mid\by_{1:t})]\succeq\E[\Cov(\bx_t\mid\calG_t)]$. Given $\calF:=\sigma(\bx_{t-1},\by_{1:t-1})$, the residual $\br_t$ has prior density $\rho^\star(\cdot\mid\bx_{t-1})$, because it is independent of the past given $\bx_{t-1}$, and the new information is $\bz=\by_t-\bH_t\bA(\bw^\star)\bx_{t-1}=\bH_t\br_t+\bv_t$. Hence $\Cov(\bx_t\mid\calG_t)=\Cov(\br_t\mid\calF,\bz)$, and $\E[\Cov(\br_t\mid\calF,\bz)\mid\calF]$ is the minimum mean-square error matrix of $\br_t$ given $\bz$. The van Trees inequality~\cite{vantrees1968detection,gill1995applications}, applied to the scalar functional $\bm c^\top\br_t$ with prior Fisher information $\bm I_0(\bx_{t-1})$ and data Fisher information $\bH_t^\top\bR_t^{-1}\bH_t$, gives $\bm c^\top\E[\Cov(\br_t\mid\calF,\bz)\mid\calF]\bm c\ge(\bm c^\top(\bm I_0+\bH_t^\top\bR_t^{-1}\bH_t)\bm c)^{-1}\ge\lambda_0$ for unit $\bm c$, since $\lambda_{\max}(\bm I_0)\le\bar J$. Taking expectations proves~\eqref{eq:PEmean}. For the last claim, $Z=\bm c^\top\Cov(\bx_t\mid\by_{1:t})\bm c\in[0,R_x^2]$ and $\lambda_0\le\E Z\le\lambda_0/2+R_x^2\Prob(Z\ge\lambda_0/2)$.
\end{proof}}

\section{Proof of Theorem \ref{thm:closed}}
\label{app_F}
\begin{proof}
\emph{Step 1: decomposition.} Let $\calF_t$ be the $\sigma$-field of $\by_{1:t}$ and of all algorithmic randomness up to the computation of $\bw_t$, and let $\calF_t^+\supseteq\calF_t$ also contain $\by_{t+1}$ and the particles $\widehat\Gamma_{t+1}$. Write $\bDelta_t=\bw_t-\bw^\star$, $\bpsi_{\bw}(\bx,\bx')=-\E_{\beps}[\bXi(\bx)^\top\bs_\kappa(\bx'-\bA(\bw)\bx+\kappa\beps\mid\bx)]$, and $\widehat\bpsi_{\bw}$ for the same quantity with $\widehat\bs_\kappa$. Then
\begin{equation}
\widehat\bg_t=\nabla R_t(\bw_t)+\bm{\xi}_t+\bm m_t+\bb_t,
\end{equation}
with $\bm m_t=\widehat\bg_t-\widehat\Gamma_{t+1}(\widehat\bpsi_{\bw_t})$, $\bb_t=\widehat\Gamma_{t+1}(\widehat\bpsi_{\bw_t})-\Gamma^\star_{t+1}(\bpsi_{\bw_t})$, and $\bm{\xi}_t=\Gamma^\star_{t+1}(\bpsi_{\bw_t})-\pi^\star_t(\bphi_{\bw_t})$, using $\nabla R_t(\bw_t)=\pi^\star_t(\bphi_{\bw_t})$ by~\eqref{eq:phi}. Since the algorithm's randomness is conditionally independent of the hidden states given the observations, the law of $(\bx_t,\bx_{t+1})$ given $\calF_t$ is $\pi^\star_t(d\bx)P_{\bw^\star}(\bx,d\bx')$. Hence $\E[\Gamma^\star_{t+1}(\bpsi_{\bw_t})\mid\calF_t]=\pi^\star_t(\bphi_{\bw_t})$ and $\E[\bm{\xi}_t\mid\calF_t]=\bzero$, with $\E\norm{\bm{\xi}_t}^2\le\bar\Psi^2$. Likewise, $\E[\bm m_t\mid\calF_t^+]=\bzero$ and $\E[\norm{\bm m_t}^2\mid\calF_t^+]\le\sigma_m^2/J_g$ by Assumption~\ref{ass:score}, and $\bb_t$, $\bm{\xi}_t$ are $\calF_t^+$-measurable.

\emph{Step 2: bias bound.} We have $\norm{\widehat\Gamma_{t+1}(\widehat\bpsi-\bpsi)}\le\sup\norm{\bXi}\,\varepsilon_s\le G_\Xi\varepsilon_s$ by Jensen's inequality and Assumption~\ref{ass:score}. Since $\bm f=\bpsi_{\bw_t}/\bar\Psi$ is admissible ($\bw_t$ is $\calF_t$-measurable), $(\E\norm{(\widehat\Gamma_{t+1}-\Gamma^\star_{t+1})(\bpsi_{\bw_t})}^2)^{1/2}\le\bar\Psi d_{t+1}$. Therefore $\E\norm{\bb_t}^2\le2G_\Xi^2\varepsilon_s^2+2\bar\Psi^2d_{t+1}^2$.

\emph{Step 3: one-step topology recursion.} Since $\bw^\star\in\calW_s$, projection is nonexpansive, and $\norm{\bDelta_{t+1}}\le\norm{\bm a_t-\eta(\bm{\xi}_t+\bm m_t)}$ with $\bm a_t=\bDelta_t-\eta\nabla R_t(\bw_t)-\eta\bb_t$. Expanding the square, $\E\langle\bDelta_t-\eta\nabla R_t(\bw_t),\bm{\xi}_t+\bm m_t\rangle=0$ because the left factor is $\calF_t$-measurable, and $\E\langle\bb_t,\bm m_t\rangle=\E\langle\bm{\xi}_t,\bm m_t\rangle=0$ by $\calF_t^+$-measurability. Bounding $2\eta^2|\langle\bb_t,\bm{\xi}_t\rangle|\le\eta^2(\norm{\bb_t}^2+\norm{\bm{\xi}_t}^2)$ gives
\[
\E\norm{\bDelta_{t+1}}^2\le\E\norm{\bm a_t}^2+\eta^2\E\norm{\bb_t}^2+\eta^2\sigma_n^2 .
\]
Because $\nabla R_t(\bw^\star)=\bzero$ and $R_t$ is $\mu$-strongly convex and $L$-smooth on $\calW_s$, $\norm{\bDelta_t-\eta\nabla R_t(\bw_t)}^2\le(1-2\eta\mu+\eta^2L^2)\norm{\bDelta_t}^2\le(1-\eta\mu)\norm{\bDelta_t}^2$ for $\eta\le\mu/L^2$. Young's inequality with $\lambda=\eta\mu/(2(1-\eta\mu))$, and the trivial case if $\eta\mu=1$, yields $\norm{\bm a_t}^2\le(1-\frac{\eta\mu}2)\norm{\bDelta_t}^2+\frac{2\eta}{\mu}\norm{\bb_t}^2$. Since $\eta\mu\le\mu^2/L^2\le1$, $\eta^2\le\eta/\mu$, so with Step~2,
\begin{equation}
e_{t+1}^2\le\Big(1-\frac{\eta\mu}2\Big)e_t^2+K d_{t+1}^2+c_e .
\label{eq:erec}
\end{equation}

\emph{Step 4: closing the loop.} By Lemma~\ref{lem:prop}, $d_{t+1}\le Z_{t+1}=\varrho Z_t+C_0(\varepsilon_F+C_Le_t)$. Young's inequality with weights $(1/\varrho,1/(1-\varrho))$ gives $Z_{t+1}^2\le\varrho Z_t^2+c_1e_t^2+c_2$. Substituting into~\eqref{eq:erec} yields~\eqref{eq:closedrec}.

\emph{Step 5: stability and limits.} $\bG$ is entrywise nonnegative. A nonnegative $2\times2$ matrix has spectral radius below one if and only if $G_{11}<1$, $G_{22}<1$, and $\det(\bI-\bG)>0$. We have $\det(\bI-\bG)=\frac{\eta\mu}{2}(1-\varrho)-Kc_1=\eta\widetilde\Delta$, which is positive if and only if~\eqref{eq:smallgain} holds; then $Kc_1<\frac{\eta\mu}2(1-\varrho)$ also gives $G_{11}<1$. Iterating~\eqref{eq:closedrec} with nonnegative forcing $\bm c$ gives $\bm v_t\preceq\bG^t\bm v_0+(\bI-\bG)^{-1}\bm c$, where
\[
(\bI-\bG)^{-1}=\frac{1}{\eta\widetilde\Delta}\begin{bmatrix}1-\varrho & K\varrho\\ c_1 & \frac{\eta\mu}2-Kc_1\end{bmatrix}.
\]
The first component of $(\bI-\bG)^{-1}\bm c$ equals $[(1-\varrho)c_e+Kc_2]/(\eta\widetilde\Delta)$, and the second equals $[c_1c_e+\frac{\eta\mu}2c_2]/(\eta\widetilde\Delta)$. Substituting the constants gives~\eqref{eq:efloor}--\eqref{eq:dfloor}.
\end{proof}

\bibliographystyle{unsrt}
\bibliography{references_verified}

\end{document}